\documentclass[zamm,a4paper,fleqn,referee]{w-art}

\usepackage{times}
\usepackage{w-thm}
\usepackage{cite}
\usepackage[english]{babel}
\usepackage{hyperref,subfigure}
\usepackage{amsmath,amssymb}
\usepackage{graphicx,color}

\definecolor {darkred}{rgb}{0.6,0,0}
\definecolor {orange}{rgb}{1,0.5,0}
\definecolor {darkgreen}{rgb}{0,0.6,0.2}
\definecolor {mediumblue}{rgb}{0,0.5,0.9}
\definecolor {purple}{rgb}{0.5,0,1}
\definecolor {lightgray}{rgb}{0.45,0.45,0.45}

\newcommand {\Appx}[1]{Appendix~\ref{#1}}
\newcommand {\Eq}[1]{Eq.~(\ref{#1})}
\newcommand {\Fig}[1]{Fig.~\ref{#1}}
\newcommand {\Sect}[1]{Section~\ref{#1}}
\newcommand {\Tab}[1]{Table~\ref{#1}}

\newcommand {\ds}{\displaystyle}

\newcommand {\pa}[2]{\frac{\partial{#1}}{\partial{#2}}}

\newcommand {\grad}{\mathrm{grad}\,}

\newcommand {\fint}[1]{\mathbf{f}_{\mathrm{int}}#1}
\newcommand {\fext}[1]{\mathbf{f}_{\mathrm{ext}}#1}

\newcommand {\nel}{{n_\mathrm{el}}}
\newcommand {\nne}{{n_\mathrm{ne}}}
\newcommand {\nno}{{n_\mathrm{no}}}
\newcommand {\Dt}{\Delta t}

\newcommand {\mrd}{\mathrm{d}}

\newcommand {\mrp}{\mathrm{p}}
\newcommand {\mrq}{\mathrm{q}}

\newcommand {\mrL}{\mathrm{L}}

\newcommand {\mrT}{\mathrm{T}}

\newcommand {\mrW}{\mathrm{W}}

\newcommand {\vphi}{\varphi}

\newcommand {\bzero}{\boldsymbol{0}}

\newcommand {\bt}{\boldsymbol{t}}

\newcommand {\bv}{\boldsymbol{v}}

\newcommand {\bx}{\boldsymbol{x}}

\newcommand {\bB}{\boldsymbol{B}}
\newcommand {\bC}{\boldsymbol{C}}

\newcommand {\bI}{\boldsymbol{I}}

\newcommand {\bT}{\boldsymbol{T}}

\newcommand {\bX}{\boldsymbol{X}}

\newcommand {\bvphi}{\mbox{\boldmath$\varphi$}}

\newcommand {\bsig}{\mbox{\boldmath$\sigma$}}

\newcommand {\mf}{\mathbf{f}}

\newcommand {\mk}{\mathbf{k}}

\newcommand {\mm}{\mathbf{m}}

\newcommand {\mpp}{\mathbf{p}}
\newcommand {\mq}{\mathbf{q}}

\newcommand {\muu}{\mathbf{u}}
\newcommand {\mv}{\mathbf{v}}

\newcommand {\mx}{\mathbf{x}}

\newcommand {\mA}{\mathbf{A}}
\newcommand {\mB}{\mathbf{B}}

\newcommand {\mJ}{\mathbf{J}}

\newcommand {\mN}{\mathbf{N}}

\newcommand {\mX}{\mathbf{X}}

\newcommand {\sB}{\mathcal{B}}

\newcommand {\sS}{\mathcal{S}}

\newcommand {\sV}{\mathcal{V}}

\newcommand {\xz}{\hat{\mx}_0}
\newcommand {\xN}{\hat{\mx}_N}
\newcommand {\xn}{\hat{\mx}_n}
\newcommand {\xnpo}{\hat{\mx}_{n+1}}
\newcommand {\xnmo}{\hat{\mx}_{n-1}}

\newcommand {\un}{\hat{\muu}_n}
\newcommand {\unpo}{\hat{\muu}_{n+1}}

\newcommand {\vz}{\hat{\mv}_0}
\newcommand {\vN}{\hat{\mv}_N}
\newcommand {\vn}{\hat{\mv}_n}
\newcommand {\vnpo}{\hat{\mv}_{n+1}}

\newcommand {\pN}{\hat{\mpp}_N}
\newcommand {\qz}{\hat{\mq}_0}
\newcommand {\qN}{\hat{\mq}_N}

\newcommand {\pnp}{\hat{\mpp}^+_n}

\newcommand {\pnpop}{\hat{\mpp}^+_{n+1}}
\newcommand {\pnpope}{\hat{\mpp}^+_{e,n+1}}
\newcommand {\pnm}{\hat{\mpp}^-_n}
\newcommand {\pnme}{\hat{\mpp}^-_{e,n}}

\newcommand {\qnp}{\hat{\mq}^+_n}

\newcommand {\qnpop}{\hat{\mq}^+_{n+1}}
\newcommand {\qnpope}{\hat{\mq}^+_{e,n+1}}
\newcommand {\qnm}{\hat{\mq}^-_n}
\newcommand {\qnme}{\hat{\mq}^-_{e,n}}

\newcommand {\intTnpo}{\int_{t_n}^{t_{n+1}}}

\newcommand{\A}{\mathcal{S}}

\hypersetup{colorlinks=true, breaklinks=true, linkcolor=black,
	menucolor=black, citecolor=black, urlcolor=black}

\begin{document}

\DOIsuffix{zamm.DOIsuffix}

\Volume{90}\Month{01}\Year{2010}

\pagespan{1}{}

\Receiveddate{XXXX}\Reviseddate{XXXX}\Accepteddate{XXXX}\Dateposted{XXXX}

\keywords{Hamilton's law of varying action, Hermite interpolation,
	nonlinear elastodynamics, symplectic integration,
	variational integrators\vspace*{2ex}}




\title[$C^1$-continuous space-time discretization based on
	Hamilton's law]
	{$\bC^\mathbf{1}$-continuous space-time discretization based on \\
	Hamilton's law of varying action}

\author[J.C.~Mergel]{Janine C.~Mergel\inst{1}}
\author[R.A.~Sauer]{Roger A.~Sauer\inst{1,}\footnote{Corresponding author
	\quad E-mail:~\textsf{sauer@aices.rwth-aachen.de},
            Phone: +49\,241\,80\,99129,
            Fax: +49\,241\,80\,628498}}
\author[S. Ober-Bl{\"o}baum]{Sina Ober-Bl{\"o}baum\inst{2}}

\address[\inst{1}]{AICES Graduate School, RWTH Aachen University,
	Templergraben 55, 52056 Aachen, Germany}
\address[\inst{2}]{Department of Engineering Science, University of Oxford,
	Parks Road, Oxford, OX1 3PJ, United Kingdom}


\begin{abstract}
We develop a class of $C^1$-continuous time integration methods that are applicable to conservative problems in elastodynamics. These methods are based on Hamilton's law of varying action. From the action of the continuous system we derive a spatially and temporally weak form of the governing equilibrium equations. This expression is first discretized in space, considering standard finite elements. The resulting system is then discretized in time, approximating the displacement by piecewise cubic Hermite shape functions. Within the time domain we thus achieve $C^1$-continuity for the displacement field and $C^0$-continuity for the velocity field. From the discrete virtual action we finally construct a class of one-step schemes. These methods are examined both analytically and numerically. Here, we study both linear and nonlinear systems as well as inherently continuous and discrete structures. In the numerical examples we focus on one-dimensional applications. The provided theory, however, is general and valid also for problems in 2D or 3D. We show that the most favorable candidate --- denoted as p2-scheme --- converges with order four. Thus, especially if high accuracy of the numerical solution is required, this scheme can be more efficient than methods of lower order. It further exhibits, for linear simple problems, properties similar to variational integrators, such as symplecticity. While it remains to be investigated whether symplecticity holds for arbitrary systems, all our numerical results show an excellent long-term energy behavior.
\end{abstract}

\maketitle


\section{Introduction} \label{s:intro}



In this work we derive a class of time integration methods for the computational analysis of deformable solids. We consider a discrete version of Hamilton's law of varying action to obtain space-time discretization schemes based on cubic Hermite functions in time.


\subsection{Overview of existing methods}

One very common approach for the numerical analysis of elastodynamic problems is the application of so-called semi-discrete procedures: Here, the (spatially and temporally) continuous system is discretized in space and time separately. First, the mechanical equilibrium equations, describing the deformation of the body, are discretized in space by means of the finite element method (FEM). At this point we refer to any standard literature on nonlinear finite elements for solids; see e.g.~Ref.~\cite{belytschko00,zienkiewicz05,wriggers08}. The spatially discrete system is then discretized in time, using for instance a finite difference scheme or collocation based on Taylor series expansion. Discretization schemes of this type include many well-known methods, such as the Newmark algorithm~\cite{newmark59}, the HHT-$\alpha$ method~\cite{hilber77}, the WBZ-$\alpha$ (or Bossak-$\alpha$) method \cite{wood80}, and the generalized-$\alpha$ method~\cite{chung93}. Besides methods fulfilling the equilibrium equations at single time steps, there exist various approaches based on weighted residuals; these consider equilibrium in an weighted-average sense and go back to the publication by Zienkiewicz~\cite{zienkiewicz77}. A weighted residual approach based on cubic Hermite interpolation in time, for instance, is discussed in Ref.~\cite{fung96}. See also the generalized method proposed by Modak and Sotelino~\cite{modak02}. In addition, there exists a broad literature on methods combining finite elements (FE) in both space and time. In general, these solution schemes can be constructed by forming a (spatially and temporally) weak form of the equations of motion, and discretizing the resulting statement by means of finite elements. The first approaches accounting for finite elements in time go back to Ref.~\cite{argyris69,oden69,fried69}. Early publications on space-time FE methods include Ref.~\cite{hughes88,hulbert90,hulbert92}. A broader literature review on time integration methods in structural mechanics can be found e.g.~in Ref.~\cite{bauchau99,kuhl99a,betsch01,gross05}.



A special class of time integration schemes applied to mechanical systems is formed by geometrical integrators. Geometric integration enables the design of robust methods that provide both quantitatively and qualitatively accurate results. Since these methods preserve the geometric properties of the flow of a differential equation, they are able to exactly represent the main characteristic properties of the physical process~\cite{marsden01,hairer06,reich94}. Geometric integration methods can be mainly divided into two classes: 1)~energy-momentum integrators and 2) symplectic momentum-conserving integrators. The first class of methods fulfills the conservation laws of energy and momentum automatically; for methods of this type we refer e.g.~to Simo and Tarnow~\cite{simo92a}, Simo et al.~\cite{simo92b}, Gonzalez~\cite{gonzalez00}, Betsch and Steinmann~\cite{betsch01}, Gro{\ss} et al.~\cite{gross05}, Leyendecker~et~al.~\cite{leyendecker06}, Hesch and Betsch~\cite{hesch09,hesch10}, Gautam and Sauer~\cite{gautam13} and Krenk~\cite{krenk14}, Betsch and Janz~\cite{betsch16}, and the references therein. See also the generalized energy-momentum method discussed in Ref.~\cite{kuhl99a,kuhl99b,kuhl00}. The second class preserves both the symplectic form and --- in the presence of symmetries --- momentum maps; it additionally shows excellent long-term energy behavior. Symplectic-momentum integrators can be represented by the class of variational integrators~\cite{marsden01,suris90}. For conservative systems, these methods are constructed by forming a discrete version of Hamilton's principle, choosing both a finite-dimensional function space and a suitable numerical quadrature. See Ref.~\cite{lew04a,lew04b} for an overview. For dissipative or controlled mechanical systems, they can be derived from a discrete version of the Lagrange-d'Alembert principle~\cite{kane00,oberbloebaum11}. Within the last years variational integrators have been extended towards constrained~\cite{leyendecker08,leyendecker10,kobilarov10,cortes01}, non-smooth~\cite{fetecau03,johnson14}, stochastic~\cite{bourabee08}, multirate and multiscale~\cite{tao10,leyendecker13} systems as well as to electric circuits~\cite{oberbloebaum13}. For variational integrators in combination with spatial discretization we refer to Ref.~\cite{lew03,wolff13,demoures15} and the references therein. Besides these semi-discrete approaches, there exists a covariant space-time discretization method by Marsden et al.~\cite{marsden98}. This multi-symplectic scheme allows for symplecticity in both space and time.



In most of the previously mentioned work the solution is approximated by using piecewise Lagrange interpolation. For a mechanical system, this leads to a smooth approximation of the position, but to discontinuities in the velocity at the discrete time steps. Besides this approach, Leok and Shingel~\cite{leok12} have developed a variational integrator based on piecewise Hermite interpolation. In their prolongation-collocation approach, not only the solution of the discrete Euler-Lagrange equations, but also its time derivatives are approximated with sufficient accuracy. This leads to a globally smooth approximation of the solution. Note that Ref.~\cite{leok12} does not include the combined discretization in both space and time.



To incorporate any initial conditions of the mechanical system explicitly, early publications on structural dynamics have considered the so-called Hamilton's law of varying action; see e.g.~Argyris and Scharpf~\cite{argyris69}, Fried~\cite{fried69}, Bailey~\cite{bailey75a,bailey75b,bailey76a,bailey76b}, Simkins~\cite{simkins78,simkins81}, and Borri et al.~\cite{borri85}. This law can be regarded as a generalization of Hamilton's principle: It accounts for any initial and final velocities by considering non-zero variations in the displacement at the boundaries of the time domain. Some of the studies mentioned above also include cubic Hermite interpolation in time for the displacement. Based on this law, a family of methods has been proposed~\cite{baruch82,riff84b} that combines different zero-variations of the displacement or velocity at initial and final time. In addition, G{\'e}radin~\cite{geradin74} has constructed from a subsequent application of Hamilton's law a time integration method based on Hermite interpolation. While variational integrators are automatically symplectic (due to the discretized action integral serving as generating function~\cite{hairer06}), it is not clear if the same properties hold for integration schemes constructed from Hamilton's law.


\subsection{Objectives}

In this paper we derive a class of space-time integration methods that are based on piecewise cubic Hermite interpolation in time. To this end, we consider Hamilton's law of varying action. We thus directly incorporate the additional boundary terms (arising from the non-zero variations) into our time integration method. Using a semi-discrete approach, we first discretize our resulting equilibrium equation in space and then in time. Instead of deriving additional conditions on the time derivatives of the approximated solution (as it is done in Ref.~\cite{leok12}) we consider independent variations of the position and velocities.

In general, one could then construct a variational integrator by 1) varying the action of the entire temporal domain, and 2) deriving from the variation a set of discrete Euler-Lagrange equations. Since for a cubic Hermite approximation, however, this would lead to an (unconditionally) unstable numerical method, we pursue a different approach: We vary the action for each discrete time interval individually, which leads to an overdetermined system of four equations. By choosing different combinations of equations, we derive a family of six different one-step methods. One of these schemes coincides with the method proposed by G{\'e}radin~\cite{geradin74}.

In fact, our time integration methods are not variational in the sense that they are not derived from the virtual action of the total time domain. We will demonstrate numerically, however, that the most favorable of our schemes --- denoted in the following as $\mrp2$-scheme --- shows similar properties like true variational integrators: an excellent long-term behavior and, for a simple harmonic oscillator, symplecticity. Interestingly, this is not the case for the variant discussed in Ref.~\cite{geradin74}.

We emphasize that the aim of this work is to both present the construction of our time integration methods, and to demonstrate their most important features by means of various numerical examples. A further analytical investigation (including e.g.~the proof of the order of convergence) goes beyond the scope of this paper; instead, this should be addressed in future work.

Note that like other methods based on $C^1$-continuous approximations in time, our integration schemes are not favorable for the simulation of discontinuous changes (such as shock waves) in mechanical systems. Instead, we apply the $\mrp2$-scheme to temporally smooth examples; these include both linear and nonlinear as well as intrinsically discrete and spatially continuous problems. Here, we focus on one-dimensional applications; the theory, however, is also valid for conservative systems in 2D and 3D. Compared to both the formulation of G{\'e}radin and the method of Leok and Shingel based on cubic Hermite interpolation, our $\mrp2$-scheme exhibits a higher order of convergence. In addition, based on the desired accuracy, it may be more efficient than classical methods like the Newmark algorithm.


\subsection{Outline}

The remainder of this paper is structured as follows. \Sect{s:hamilton} introduces the action integral of a continuous body deforming over time. From Hamilton's law of varying action a (spatially and temporally) weak form of the mechanical equilibrium equation is derived. \Sect{s:dAh} briefly outlines its spatial discretization by means of standard finite elements. The temporal discretization is discussed in \Sect{s:tempdiscr}, providing a solution strategy that leads to a class of different integration methods. These methods are then related to other approaches from the literature. In \Sect{s:props}, we study the main characteristic properties of our integration schemes, such as symplecticity and the long-term and convergence behavior. The most favorable scheme is then applied to investigate several numerical examples (\Sect{s:results}). \Sect{s:concl} finally concludes this paper.


\section{Hamilton's law of varying action} \label{s:hamilton}

In this section we summarize the governing equations describing a body undergoing finite motion and deformation. For the general theory of continuum mechanics, the reader is referred to text books~\cite{chadwick99,holzapfel00}. Consider a body deforming within the time domain $[0,T]$. In the initial configuration, at $t = 0$, the body is denoted by $\sB_0$; its boundary is denoted by $\partial\sB_0$. The body can be subjected to volumetric loads,~$\bar{\bB}$ (applied in $\sB_0$), deformations, $\bar{\bvphi}$ (prescribed on $\partial_{\boldsymbol{\varphi}}\sB_0 \subseteq \partial\sB_0$), and surface loads, $\bar{\bT}$ (applied on~$\partial_{\bt}\sB_0 \subseteq \sB_0$).\footnote{We assume that $\partial_{\boldsymbol{\varphi}}\sB_0 \cup \partial_{\bt}\sB_0 = \partial\sB_0$ and $\partial_{\boldsymbol{\varphi}}\sB_0 \cap \partial_{\bt}\sB_0 = \emptyset$.} At any time $t \in (0,T]$, the deformation of the body is characterized by a unique mapping of a material point, $\bX\in\sB_0$, to its current position, $\bx = \bvphi(\bX,t) \in \sB$. The material time derivative $\bv := \partial\bx/\partial t$ corresponds to the velocity of a material particle located at $\bx$; in short we will also write $\bv = \dot{\bx}$.


\subsection{Action of the continuous system} \label{s:A}

In the following we assume 1)~conservation of mass, 2)~hyperelastic material behavior, and 3)~that the external forces do not depend on the deformation. We start with the action integral of the continuous system, defined as
\begin{equation}
	\A = \int_0^T L(\bx,\bv)~\mrd t. \label{e:A}
\end{equation}
The integrand corresponds to the Lagrangian of our system, given by
\begin{equation}
	L(\bx,\bv) = K(\bv) - \Pi(\bx).
\end{equation}
Here, $K$ is the kinetic energy, and $\Pi$ is the potential energy due to both internal strains and external forces;
\begin{align}
	K(\bv)		&	= \frac{1}{2}\int_{\sB_0}\rho_0\,\bv\cdot\bv~\mrd V, \\
	\Pi(\bx)	&	= \int_{\sB_0} W(\bx)~\mrd V
							-	\int_{\sB_0} \bx\cdot\rho_0\,\bar{\bB}~\mrd V
							- \int_{\partial_{\bt}\sB_0} \bx\cdot\bar{\bT}~\mrd A.
							\label{e:Pi}
\end{align}
The terms $\rho_0$ and $W(\bx)$ respectively denote the initial material density and an energy density function characterizing the material behavior; a detailed description of $W(\bx)$ for different material models can be found e.g.~in Ref.~\cite{belytschko00,zienkiewicz05,wriggers08}.


\subsection{Variation of the action} \label{s:dA}

We now consider an admissible variation of the deformation, $\delta\bx\in\sV$,
\begin{equation}
	\sV  = \Big\{ \delta\bx:\ \sB_0\times [0,T] \rightarrow \left.{
	\mathbb{R}^d}\right|\ \delta\bx(\bX,t)|_{\partial
	_{\boldsymbol{\vphi}}\sB_0} = \bzero \Big\}. \label{e:spaceV}
\end{equation}
Here, $d$ is the dimension of Euclidean space. Varying the action integral \eqref{e:A} yields
\begin{equation}
	\delta\A = \int_0^T \delta L(\bx,\bv)~\mrd t, \qquad
	\delta L(\bx,\bv) = \delta K(\bv) - \delta\Pi(\bx), \label{e:dAdef}
\end{equation}
where the variations of the energy terms are given by
\begin{align}
	\delta K	(\bv)	&	=	\int_{\sB_0} \delta\bv \cdot \rho_0 \, \bv~\mrd V,
										\label{e:dK} \\
	\delta\Pi(\bx)	&	=	\int_{\sB}\grad{\delta\bx}:\bsig~\mrd v
										- \int_{\sB_0} \delta\bx\cdot\rho_0\,\bar{\bB}~\mrd V
										- \int_{\partial_{\bt}\sB_0} \delta\bx
										\cdot \bar{\bT}~\mrd A. \label{e:dPi}
\end{align}
The tensor $\bsig$ denotes the Cauchy stress, which is derived from the strain energy density function,~$W$, appearing in \Eq{e:Pi}.

If we enforce the deformation at $t = 0$ and $t = T$ to be fixed (s.t.~the variations $\delta\bx(\bX,0)$ and $\delta\bx(\bX,T)$ become zero), we arrive at the classical Hamilton's principle,
\begin{equation}
	\delta\A = 0 \qquad \forall\ \delta\bx \in \Big\{ \delta\bx\in\sV\,\big|\ 
	\delta\bx(\bX,0) = \bzero,\ \delta\bx(\bX,T) = \bzero \Big\}, \label{e:dA0}
\end{equation}
see e.g.~the book by Lanczos~\cite{lanczos70}. Instead, however, we leave the variations $\delta\bx(\bX,0)$ and $\delta\bx(\bX,T)$ arbitrary for now. In this case, $\delta\A$ is equal to the following boundary term evaluated at $t = 0$ and $t = T$:
\begin{equation}
	\delta\A = \left.{\left\langle{\delta\bx \, , \,
	\frac{\partial L (\bx,\bv)} {\partial\bv}}\right\rangle}\right|_0^T
	\qquad \forall\ \delta\bx \in\sV. \label{e:dA}
\end{equation}
The term in angle brackets corresponds to the scalar product of the variation $\delta\bx$ and the linear momentum of the body, i.e.
\begin{equation}
	\left\langle{\delta\bx \, , \, \frac{\partial L (\bx,\bv)}
	{\partial\bv}}\right\rangle = \int_{\sB_0} \delta\bx \cdot \rho_0 \,
	\bv~\mrd V. \label{e:scalprod}
\end{equation}
A more detailed discussion of arbitrary variations at the boundaries can be found in Ref.~\cite{lanczos70} (see Chapter~V.3 there). Due to the scalar product, \Eq{e:dA} depends on the system's initial momentum (and thus on its initial velocity) explicitly. This is the reason why expression~\eqref{e:dA} has been discussed in several early publications studying initial value problems in structural dynamics, such as~Ref.~\cite{argyris69,bailey75a,bailey75b,simkins78,simkins81,baruch82}; an overview can also be found in Ref.~\cite{peters88}. Following the terminology used in several of these papers, we refer to \Eq{e:dA} as ``Hamilton's law of varying action''.\footnote{Alternatively, this equation has been referred to as ``Hamilton's weak principle'' (HWP); see e.g.~Ref.~\cite{hodges91}.} A very interesting comment on its origin can be found in a paper by Bailey (Ref.~\cite{bailey75b}, p.~434): ``When copies of Hamilton's original papers \cite{hamilton34,hamilton35} were obtained, it was found that Hamilton had furnished what he called the `law of varying action'. He did not furnish what is now known as `Hamilton's principle'. Evidently, in the latter part of the 19th century, application of the concepts of the variational calculus of Euler and Lagrange reduced Hamilton's law to Hamilton's principle.''

We will explain later why the two additional boundary terms (caused by the non-zero variations) are important for the derivation of our final integration schemes; see \Sect{s:sol}. Equation~\eqref{e:dA} finally represents the (spatially and temporally) weak form of the governing equilibrium equations. Note that this expression is general and valid for elastodynamic problems in $\mathbb{R}^d$, $d \in \{1,2,3\}$. Nevertheless, since we focus on the development and analysis of a new time integration method, we numerically investigate only one-dimensional problems with $\bar{\bB} = \bar{\bT} = \bzero$ in this paper. A detailed study of problems in both 2D or 3D may be the subject of future work.


\section{Spatial discretization} \label{s:dAh}

We now briefly outline the spatial discretization by means of the finite element method. Regarding nonlinear finite elements for solids we refer to text books such as Ref.~\cite{belytschko00,zienkiewicz05,wriggers08}. We spatially discretize \Eq{e:dA} by using $\nel$ finite elements; for each element,~$\Omega^e$, the initial position,~$\bX$, the deformation, $\bx(\bX,t)$, and the velocity, $\bv(\bX,t)$, are approximated by
\begin{equation}
	\bX^h(\bX) = \mN_e(\bX) \, \mX_e, \qquad
	\bx^h(\bX,t) = \mN_e(\bX) \,\,\mx_e(t), \qquad
	\bv^h(\bX,t) = \mN_e(\bX) \,\,\dot{\mx}_e(t), \label{e:xh}
\end{equation}
where the dot indicates the derivative with respect to time. The vectors
\begin{equation}
	\mX_e = \begin{bmatrix} \mX_1 \\ \vdots \\ \mX_\nne
	\end{bmatrix}, \qquad \mx_e(t) = \begin{bmatrix} \mx_1(t) \\ \vdots \\
	\mx_\nne(t) \end{bmatrix}, \qquad \dot{\mx}_e(t) = \begin{bmatrix}
	\dot{\mx}_1(t) \\ \vdots \\ \dot{\mx}_\nne(t) \end{bmatrix}
\end{equation}
contain the initial and current positions as well as the velocities of those $\nne$~nodes belonging to element~$\Omega^e$. These quantities are still continuous with respect to time. The array $\mN_e = \big[ N_1 \, \bI_d \,,\, \dots \,,\, N_{n_\mathrm{ne}} \, \bI_d \big]$ contains the nodal shape functions $N_1$ -- $N_{n_\mathrm{ne}}$ associated with~$\Omega^e$. Using an isoparametric concept, we discretize the variations by means of the same shape functions,
\begin{equation}
	\delta\bx^h(\bX,t) = \mN_e(\bX) \, \delta\mx_e(t), \qquad
	\delta\bv^h(\bX,t) = \mN_e(\bX) \, \delta\dot{\mx}_e(t).
\end{equation}
We can further write
\begin{equation}
	\grad{\delta\bx}^h = \mB_e(\bX) \, \delta\mx_e(t)
\end{equation}
for a suitable definition of the strain operator $\mB_e$. By inserting these relations into \Eq{e:dA} and \eqref{e:scalprod}, we obtain
\begin{equation}
	\delta\A^h = \sum_{e=1}^\nel \, {\Big[{ \delta\mx_e^\mrT \,
	\mm_e \, \dot{\mx}_e} \Big]\Big|}_0^T, \label{e:dAh}
\end{equation}
where $\delta\A^h$ can be computed from
\begin{equation}
	\delta\A^h = \int_0^T \delta L^h (\mx,\dot{\mx})~\mrd t, \qquad \delta L^h
	(\mx,\dot{\mx}) = \sum_{e=1}^\nel \left[{\delta \dot \mx_e^\mrT \,
	\mm_e \, \dot\mx_e - \delta\mx_e^\mrT \, \mf^e}\right]; \label{e:dAhdef}
\end{equation}
cf.~\Eq{e:dAdef} -- \eqref{e:dPi}. Here, the vectors $\mx(t)$ and $\dot{\mx}(t)$ denote the deformation and velocity at all spatial FE nodes. The elemental mass matrices, $\mm_e$, and the force vectors, $\mf^e := \fint^e - \fext^e$, are computed through
\begin{align}
	\mm_e		&	:=	\int_{\Omega^e_0} \rho_0 \, \mN^\mrT_e\,\mN_e~\mrd V,
							\label{e:me} \\
	\fint^e	&	:=	 \int_{\Omega^e} \mB_e^\mrT \, \bsig~\mrd v, \label{e:finte} \\
	\fext^e	&	:= \int_{\Omega^e_0} \rho_0 \, \mN^\mrT_e\,\bar{\bB}~\mrd V
							+ \int_{\Gamma_{0\bt}^e} \mN^\mrT_e \, \bar{\bT}~\mrd A,
\end{align}
introducing $\Gamma^e_{0\bt} = \Omega^e_0 \cap \partial_{\bt}\sB^h_0$. To obtain a shorter notation we will later refer to the global mass matrix,~$\mm$, and the global force vectors, $\mf := \fint - \fext$, assembled from the elemental contributions. Equation~\eqref{e:dAh} finally corresponds to the spatially discrete version of Hamilton's law of varying action.


\section{Temporal discretization} \label{s:tempdiscr}

We now discretize the (spatially discrete) virtual action, $\delta\A^h$, in time. To achieve temporal $C^1$-continuity, we approximate the nodal deformation of element $\Omega^e$ by cubic Hermite shape functions; $\mx_e(t) \approx \mx_e^t(t)$ with
\begin{equation}
	\mx^t_e(t) = R_1(t) \, \xn^e + R_2(t) \, \xnpo^e + H_1(t) \, \vn^e
	+ H_2(t) \, \vnpo^e, \qquad t \in [t_n,t_{n+1}], \label{e:xt}
\end{equation}
where $n = 0,\dots,N\!-\!1$, $t_0 = 0$, and $t_N = T$. The vectors $\xn^e$ and $\xnpo^e$ contain the nodal deformations at $t_n$ and $t_{n+1}$; the vectors $\vn^e$ and $\vnpo^e$ are the corresponding nodal velocities. See \Appx{a:shpfct} for the definition of the shape functions $R_1$, $R_2$, $H_1$, and $H_2$. We can further write
\begin{align}
	\dot{\mx}^t_e(t)	&	= \dot{R}_1(t) \ \ \,\xn^e + \dot{R}_2(t) \ \ \,\xnpo^e
											+ \dot{H}_1(t) \ \ \,\vn^e +\dot{H}_2(t) \ \ \,\vnpo^e,
											\label{e:vt} \\
	\delta\mx^t_e(t)	&	= R_1(t) \, \delta\xn^e + R_2(t) \, \delta\xnpo^e
											+ H_1(t) \, \delta\vn^e + H_2(t) \, \delta\vnpo^e, \\
	\delta\dot{\mx}_e^t(t)
							&	= \dot{R}_1(t) \, \delta\xn^e + \dot{R}_2(t) \, \delta\xnpo^e
								+ \dot{H}_1(t) \, \delta\vn^e + \dot{H}_2(t) \, \delta\vnpo^e.
									\label{e:dvt}
\end{align}
In order to improve readability, we will also use the assembled counterparts accounting for all spatial FE nodes at once. We will denote them e.g.~by $\xn$ instead of $\xn^e$; see the analogy to $\mx(t)$ / $\mx_e(t)$ introduced in the previous section.


\subsection{Virtual action for a single time interval} \label{s:dAht}

Approximations \eqref{e:xt} -- \eqref{e:dvt} are now inserted into the virtual action for a single time interval, $[t_n,t_{n+1}]$; this yields
\begin{equation}
	\delta\A^{ht}_{n+1} = \intTnpo \delta L^{ht} \big(\mx^t,\dot{\mx}^t\big)
	~\mrd t, \qquad \delta L^{ht} \big(\mx^t,\dot{\mx}^t\big)
	= {\left( {\delta\dot\mx^t} \right)}^\mrT \, \mm \, \dot\mx^t - {\left(
	{\delta\mx^t} \right)}^\mrT \, \mf\big(\mx^t\big). \label{e:dAhtdef}
\end{equation}
The increment $\delta\A^{ht}_{n+1}$ depends on four variables: $\xn$, $\xnpo$, $\vn$, and $\vnpo$; it can thus be reformulated to
\begin{equation}
	\delta\A^{ht}_{n+1} = \delta\xn \cdot \pa{\A^{ht}_{n+1}}{\xn} + \delta\xnpo
	\cdot \pa{\A^{ht}_{n+1}}{\xnpo} + \delta\vn \cdot \pa{\A^{ht}_{n+1}}{\vn}
	+ \delta\vnpo \cdot \pa{\A^{ht}_{n+1}}{\vnpo}. \label{e:dAhtdetail}
\end{equation}
Following the terminology of Marsden and West~\cite{marsden01}, we now define the discrete momenta
\begin{equation}
	\pnm := -\,\pa{\A^{ht}_{n+1}}{\xn}, \qquad
	\pnpop := \pa{\A^{ht}_{n+1}}{\xnpo}. \label{e:defp}
\end{equation}
We further introduce two analogous variables that appear due to the Hermite discretization~\eqref{e:xt},
\begin{equation}
	\qnm := -\,\pa{\A^{ht}_{n+1}}{\vn}, \qquad
	\qnpop := \pa{\A^{ht}_{n+1}}{\vnpo}. \label{e:defq}
\end{equation}
Since $\qnm$ and $\qnpop$ have the unit ``momentum $\times$ time'', we refer to them as discrete ``pseudo-momenta''. The four terms can be computed from the discrete action, $\sS^{ht}_{n+1}$, given in \Appx{a:Aht}; this results in
\begin{align}
	\pnm		&	= - \intTnpo \! \Big[ \dot{R}_1(t) \, \mm \, \dot{\mx}^t
						- R_1(t) \, \mf \big(\mx^t\big)\Big]~\mrd t, \quad\ \ 
	\, \pnpop =	 \intTnpo \! \Big[ \dot{R}_2(t) \, \mm \, \dot{\mx}^t
						- R_2(t) \, \mf \big(\mx^t\big) \Big]~\mrd t, \label{e:pnmp} \\
	\qnm		&	= - \intTnpo \! \Big[ \dot{H}_1(t) \, \mm \, \dot{\mx}^t
						- H_1(t) \, \mf \big(\mx^t\big) \Big]~\mrd t, \quad\ \ 
	\qnpop	 = \intTnpo \! \Big[ \dot{H}_2(t) \, \mm \, \dot{\mx}^t
						- H_2(t) \, \mf \big(\mx^t\big) \Big]~\mrd t. \label{e:qnmp}
\end{align}
By inserting these (pseudo-)momenta, we can simplify \Eq{e:dAhtdetail} to
\begin{equation}
	\delta\A^{ht}_{n+1} = - \delta\xn \cdot \pnm + \delta\xnpo \cdot \pnpop
	- \delta\vn \cdot \qnm + \delta\vnpo \cdot \qnpop. \label{e:dAhtabbr}
\end{equation}
In addition, since the variations $\delta\xn$ and~$\delta\xnpo$ remain arbitrary, $\delta\A^{ht}_{n+1}$ must fulfill Hamilton's law of varying action applied to the interval $[t_n,t_{n+1}]$,
\begin{equation}
	\delta\A^{ht}_{n+1} = \delta\xnpo \cdot \mm \, \vnpo
	- \delta\xn \cdot \mm \, \vn. \label{e:dAht}
\end{equation}
After inserting \Eq{e:dAhtabbr}, we finally obtain the spatially and temporally discrete version of \Eq{e:dA} for a single time interval:
\begin{equation}
	\delta\xn \cdot \Big[{ \mm\,\vn - \pnm }\Big] + \delta\xnpo \cdot
	\Big[{ \pnpop - \mm\,\vnpo }\Big] + \delta\vn \cdot \Big[{ -\qnm }\Big]
	+ \delta\vnpo \cdot \Big[{ \qnpop }\Big] = 0 \label{e:dAhtfinal}
\end{equation}
for arbitrary variations $\delta\xn$, $\delta\xnpo$, $\delta\vn$, and $\delta\vnpo$. Note that we derive \Eq{e:dAhtfinal} by first varying the continuous action integral, and then discretizing its variation in space and time. We would, however, obtain the same expression if we first discretized the action itself, and then varied the discrete action for a single time interval.


\subsection{Solution strategy} \label{s:sol}

\begin{figure}[ht]
	\centering
	\includegraphics[width=0.85\textwidth]{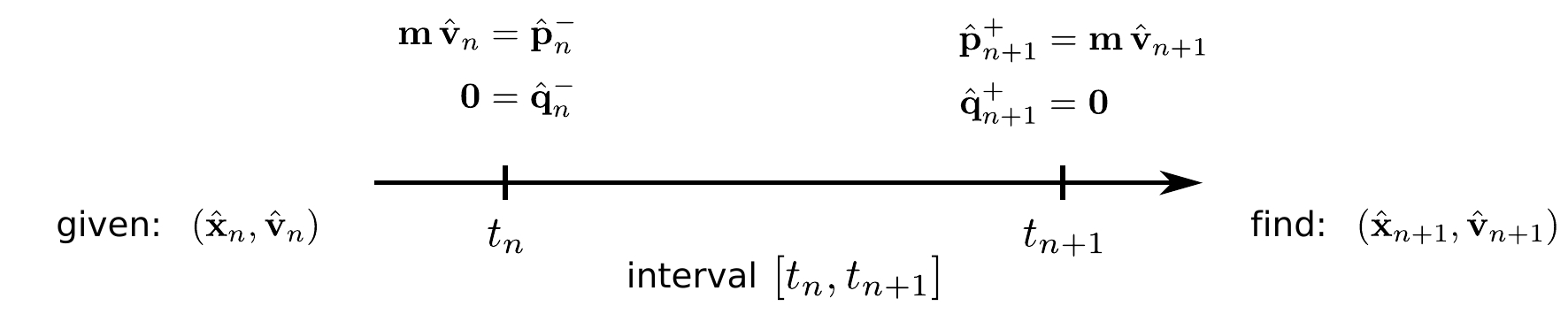}
	\caption{Equilibrium equations for time interval $[t_n,t_{n+1}]$.}
	\label{f:eqs}
\end{figure}

\begin{figure}[ht]
	\centering
	\hspace*{-5ex}\subfigure[{$\mrp2$-scheme}]{
		\includegraphics[width=0.48\textwidth]{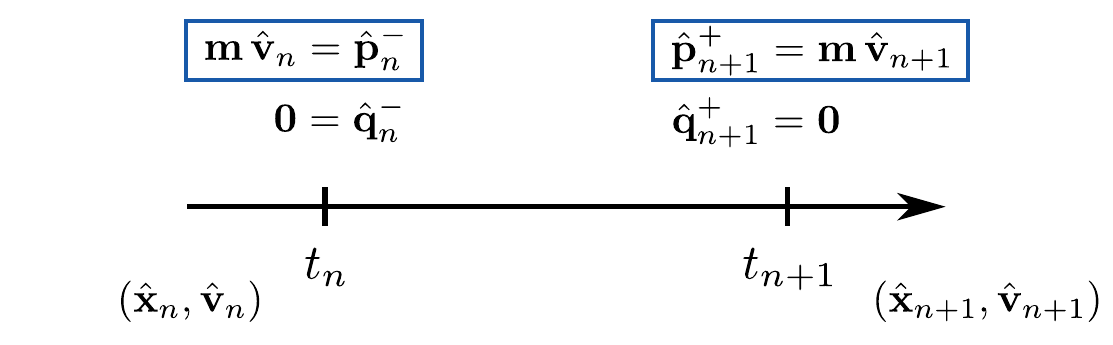}
	}\subfigure[{$\mrq2$-scheme}]{
		\includegraphics[width=0.48\textwidth]{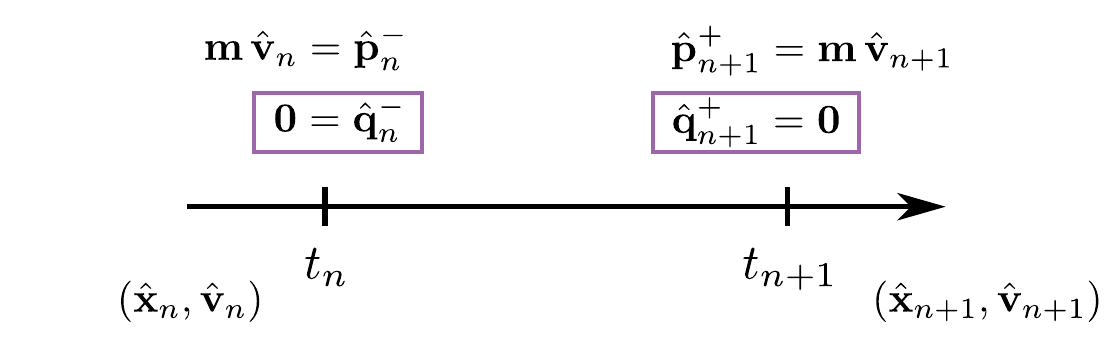}}
	
	\vspace*{1ex}
	\hspace*{-5ex}\subfigure[{$\mrp^+\mrq^-$-scheme}]{
		\includegraphics[width=0.48\textwidth]{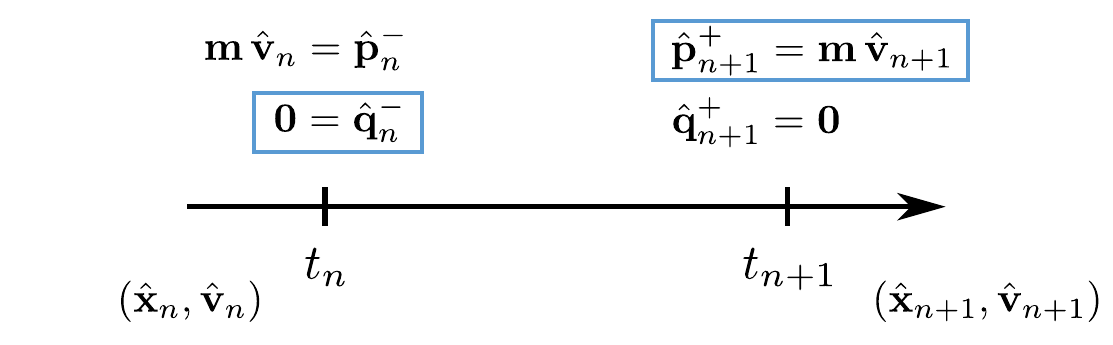}
	}\subfigure[{$\mrp^+\mrq^+$-scheme}]{
		\includegraphics[width=0.48\textwidth]{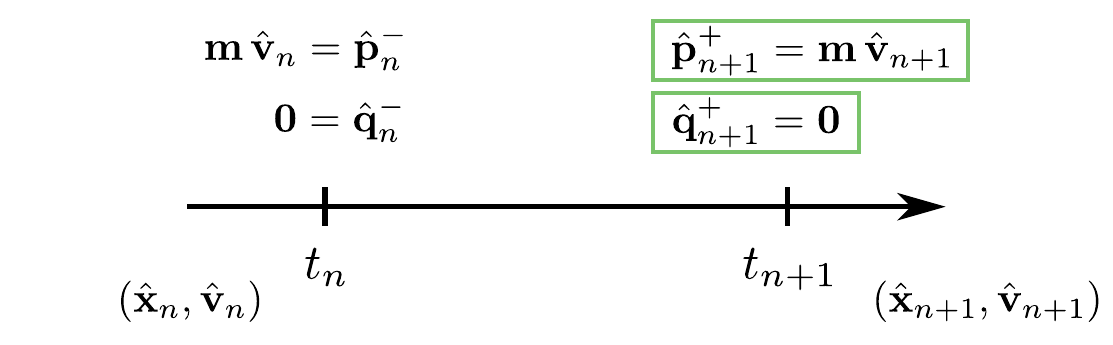}}
	
	\vspace*{1ex}
	\hspace*{-5ex}\subfigure[{$\mrp^-\mrq^-$-scheme}]{
		\includegraphics[width=0.48\textwidth]{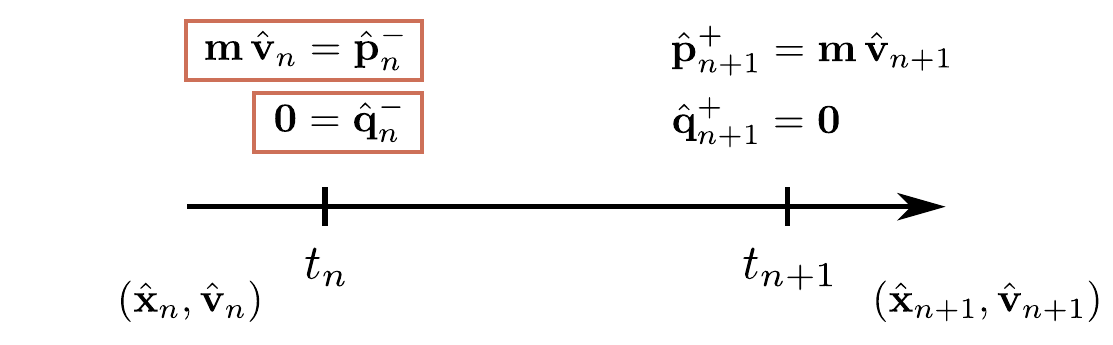}
	}\subfigure[{$\mrp^-\mrq^+$-scheme}]{
		\includegraphics[width=0.48\textwidth]{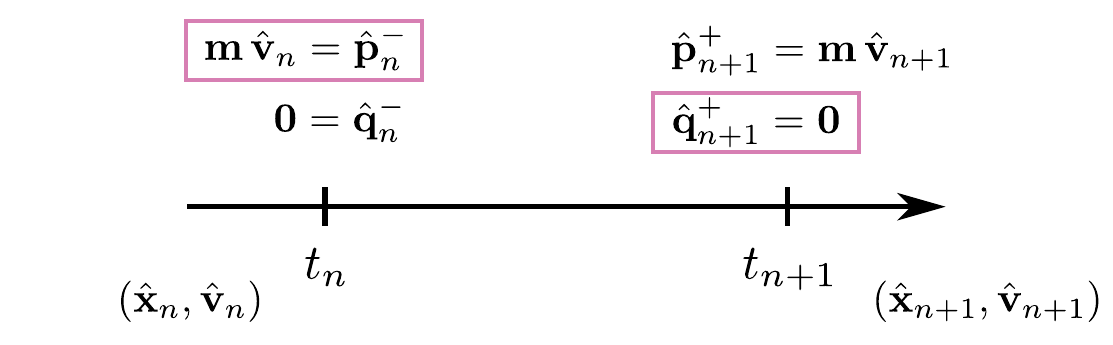}}
	\caption{Illustration of the six integration schemes.}
	\label{f:schemes}
\end{figure}

In general, variational integrators are constructed by 1)~summing up the discrete action for all time intervals, 2)~taking its variation, and 3)~re-arranging the summands. Doing so for \Eq{e:dAhtfinal}, we would arrive at the following statement:
\begin{align}
	&	\delta\xz \cdot \Big[ \mm\,\vz - \hat{\mpp}_0^- \Big]
		+ \sum_{n=1}^{N-1} \delta\xn \cdot \Big[ \pnp - \pnm \Big]
		+ \delta\xN \cdot \Big[ \mm\,\vN - \pN^+ \Big] \nonumber \\
	& \hspace*{2.5cm} + \delta\vz \cdot \Big[ -\qz^- \Big]
		+ \sum_{n=1}^{N-1} \delta\vn \cdot \Big[ \qnp - \qnm \Big]
		+ \delta\vN \cdot \Big[ \qN^+ \Big] = 0. \label{e:dAhtinst}
\end{align}
This expression is equivalent to \Eq{e:dAhtfinal} summed up over the entire time domain. Solving \Eq{e:dAhtinst} subsequently, however, results in a two-step method that is unconditionally unstable (i.e.~for which the spectral radius is larger than one). This observation is also discussed in Ref.~\cite{riff84a}. Instead, we develop a class of one-step methods arising from the virtual action for the individual time interval $[t_n,t_{n+1}]$. If the virtual displacements and velocities are presumed to be arbitrary, \Eq{e:dAhtfinal} provides us with $(4\cdot d\,\nno)$ equations,
\begin{align}
	\pnm			&	=	\mm \, \vn, \label{e:eqpm} \\
	\pnpop		&	=	\mm \, \vnpo, \label{e:eqpp} \\
	\qnm			&	=	\bzero, \label{e:eqqm} \\
	\qnpop		&	=	\bzero, \label{e:eqqp}
\end{align}
where $\nno$ is the number of finite element nodes. Physically, the first two equations relate the discrete momenta, $\pnm$ and~$\pnpop$, to the linear momenta at $t_n$ and $t_{n+1}$; see \Fig{f:eqs}. The second two equations arise from the chosen Hermite approach.
Assuming that the displacement and velocity of the previous time step, $\xn$ and $\vn$, are given, we need only $(2\cdot d\,\nno)$ equations to determine the new state, $\xnpo$ and $\vnpo$. The system (\ref{e:eqpm}) -- (\ref{e:eqqp}) is thus over-determined. For this reason, we set two of the (so far arbitrary) variations, $\delta\xn$, $\delta\xnpo$, $\delta\vn$, and $\delta\vnpo$, to zero; this approach is further motivated in the following. The new deformation and velocity are then computed from the remaining two equations. We finally obtain six methods, which are illustrated in \Fig{f:schemes}.
\vspace*{1ex}
\begin{description}
	\item[$\mpp\mathbf{2}$-scheme] \quad $\delta\vn = \delta\vnpo
		= \bzero$
		\begin{equation}
			\pnm = \mm \, \vn, \qquad \pnpop = \mm \, \vnpo. \label{e:pp}
		\end{equation}
		This seems to be the most promising approach: From its definition follows
		that it enforces a matching of the momenta at the discrete time
		steps, i.e.~that $\pnm = \pnp = \mm \, \vn$.
		\vspace*{1ex}
	\item[$\mq\mathbf{2}$-scheme] \quad $\delta\xn = \delta\xnpo = \bzero$
		\begin{equation}
			\qnm = \bzero, \qquad \qnpop = \bzero. \label{e:qq}
		\end{equation}
		This method can be seen as the counterpart of the $\mrp2$-scheme.
		\vspace*{1ex}
	\item[$\mpp^+\mq^-$-scheme] \quad $\delta\xn = \delta\vnpo = \bzero$
		\begin{equation}
			\pnpop = \mm \, \vnpo, \qquad \qnm = \bzero. \label{e:ppqm}
		\end{equation}
		This is one of four mixed methods, varying once the displacement and
		once the velocity.
		\vspace*{1ex}
	\item[$\mpp^+\mq^+$-scheme] \quad $\delta\xn = \delta\vn = \bzero$
		\begin{equation}
			\pnpop = \mm \, \vnpo, \qquad \qnpop = \bzero. \label{e:ppqp}
		\end{equation}
		This scheme corresponds to the formulation proposed by
		G{\'e}radin~\cite{geradin74}.
		\vspace*{1ex}
	\item[$\mpp^-\mq^-$-scheme] \quad $\delta\xnpo = \delta\vnpo = \bzero$
		\begin{equation}
			\pnm = \mm \, \vn, \qquad \qnm = \bzero. \label{e:pmqm}
		\end{equation}
	\item[$\mpp^-\mq^+$-scheme] \quad $\delta\xnpo = \delta\vn = \bzero$
		\begin{equation}
			\pnm = \mm \, \vn, \qquad \qnpop = \bzero. \label{e:pmqp}
		\end{equation}
\end{description}
Interestingly, the resulting six methods have completely different characteristics; it is not surprising that some of them are more favorable than others. In \Sect{s:props} we numerically investigate the properties of each scheme in terms of preservation of energy and convergence behavior. Furthermore, for linear systems we analyze the stability and symplecticity of the schemes. For several reasons we especially focus on the first, i.e.~on the $\mrp2$-scheme. Here, the discrete velocities, $\vn$ and $\vnpo$ --- appearing in the Hermite ansatz~\eqref{e:xt} --- are connected to the displacements by setting the linear momenta, $\mm\,\vn$ and $\mm\,\vnpo$, equal to the discrete momenta, $\pnm$ and $\pnpop$. Note that this approach is only possible because we explicitly account for the boundary terms appearing in Hamilton's law of varying action. Second, since from definition~(\ref{e:pp}) directly follows that $\pnp = \pnm$, the discrete Euler-Lagrange equations, $\pnp - \pnm = \bzero$, are fulfilled automatically (cf.~Chapter VI, Eq.~(6.7) of Hairer et al.~\cite{hairer06}). Interestingly, the $\mrp2$-scheme satisfies balance of linear momentum (a generalization of conservation of momentum) averaged over the time step~(\Appx{a:ppreform}):
\begin{equation}
	\intTnpo \mm \, \ddot{\mx}^t + \mf \big( \mx^t \big)~\mrd t = 0;
	\label{e:ppreform}
\end{equation}

Nevertheless, we must point out that our six integrators are not variational; this becomes apparent especially for the four mixed methods, where we use one zero-variation for the displacement, and one for the velocity. We will show, however, that --- at least for simple linear problems --- the first two methods \eqref{e:pp} and \eqref{e:qq} have similar properties like variational integrators, such as symplecticity. In addition, we numerically demonstrate that these methods show a very good energy-preserving behavior even for nonlinear problems with multiple degrees of freedom. In future work it would be interesting to compare these schemes with variational integrators in more detail, and to investigate symplecticity for arbitrary systems. It would be further interesting to examine whether the $\mrp2$-schemes conserves momentum maps associated with symmetries of the Lagrangian.


\subsection{Implementation} \label{s:implement}

In general, equations~\eqref{e:eqpm} -- \eqref{e:eqqp} are nonlinear. They thus must be linearized by using e.g.~Newton's method; this provides a system of linear equations that is iteratively solved for the new positions and velocities of the finite element nodes. For linearization, the derivatives of the discrete (pseudo-)momenta are required; see \Appx{a:implempq}. In analogy to the force vectors and the mass matrix, the terms $\pnm$, $\pnpop$, $\qnm$, and $\qnpop$ can be computed by assembling the contributions of each spatial element, denoted by $\pnme$, $\pnpope$, $\qnme$, and $\qnpope$. Where possible, the integrals should be computed analytically; this can be done for the contributions due to 1)~the kinetic energy and due to 2)~any linear elastic internal energy (\Appx{a:implempq}). The remaining integrals are evaluated by Gaussian quadrature, choosing a sufficient number of quadrature points. For our time integration schemes, neither the kinetic/potential energy nor the total energy of the system must be evaluated explicitly. Since we want to investigate them in our numerical examples, however, we discuss these quantities in \Appx{a:Aht}.


\subsection{Relation to other methods} \label{s:others}

The idea of applying Hamilton's law of varying action to initial value problems in structural dynamics goes back to the first approaches using finite elements in both space and time; see e.g.~Ref.~\cite{argyris69,fried69,bailey75a}. Instead of zero-variations of the displacement at initial and final time (as required for Hamilton's principle), these publications account for zero-variations of both the initial displacement and velocity: $\delta\mx(0) = \delta\dot{\mx}(0) = \bzero$.

This idea has motivated Baruch and Riff~\cite{baruch82} to combine different zero-variations of either the displacement or velocity at both $t = 0$ and $t = T$. Their approach results in six different methods that can be related --- with several important differences --- to our schemes. Since the same authors have discovered in a previous work~\cite{riff84a} the instability of the solution scheme given by \Eq{e:dAhtinst}, they propose a modified discretization of the virtual displacements in Ref.~\cite{riff84b}. In their approach, $\delta\mx(t)$ is discretized by considering the second derivatives of the shape functions, $\ddot{R}_\bullet(t)$ and $\ddot{H}_\bullet(t)$; the variation of the displacement is thus approximated by a linear (instead of a cubic) function in time. This modification leads to different partial derivatives of the action, and therefore to a different integration method. An even more important difference to our schemes concerns the zero-variations at the boundaries: While Riff and Baruch~\cite{riff84b} consider the boundaries of the entire time domain ($t = 0$ and $t = T$), we derive our schemes from zero-variations within each time interval,~$[t_n,t_{n+1}]$. Our approach results in six different one-step methods solving the equations subsequently. In contrast, this is possible only for the so-called F4-method of Ref.~\cite{baruch82,riff84b}, where $\delta\mx(T) = \delta\dot{\mx}(T) = \bzero$. For the remaining formulations in Ref.~\cite{baruch82,riff84b}, all equations would have to be solved simultaneously. In summary, one could loosely relate our six schemes to a subsequent application of the methods by Riff and Baruch~\cite{riff84b} for each time interval, $[t_n,t_{n+1}]$. The underlying equations, however, are approximated differently. Besides that, the references mentioned above discuss only linear dynamic systems (where the forces depend on the displacement linearly).

Recently, Leok and Shingel~\cite{leok12} have proposed a variational integrator based on Hermite finite elements in time. Their formulation is derived from a prolongation-collocation approach: In addition to the discrete Euler-Lagrange equations this method accounts for the system's equation of motion in strong form,
\begin{equation}
	\mm \, \ddot{\mx}^t(t_\bullet) + \mf \big[\mx^t(t_\bullet)\big]
	= \bzero, \qquad \bullet \in \{n,n+1\}. \label{e:eqmotion}
\end{equation}
For cubic Hermite shape functions --- as they are used in our schemes --- the velocities, $\vn$ and $\vnpo$, are computed from \Eq{e:eqmotion}, using
\begin{equation}
	\ddot{\mx}^t(t) = \ddot{R}_1(t) \, \xn + \ddot{R}_2(t) \, \xnpo
	+ \ddot{H}_1(t) \, \vn + \ddot{H}_2(t) \, \vnpo;
\end{equation}
cf.~\Eq{e:vt}. These expressions are inserted into the temporally discrete action for one time interval,~$\A^t_{n+1}$, which then depends only on the displacements, $\xn$ and $\xnpo$. The final time integration method of Leok and Shingel is obtained by 1)~varying the incremental action, $\A^t_{n+1}$, with respect to the displacements, and 2)~setting the total virtual action to zero. Compared to our six Hermite formulations, the resulting method requires only half the number of unknowns to be solved within each time step. Its rate of convergence, however, is lower than the best of our schemes; see \Sect{s:conv:single}. Note that the combination of both spatial and temporal discretizations is not discussed in Ref.~\cite{leok12}.


\section{Properties of the six schemes} \label{s:props}

We now investigate the different properties of the six formulations, first focusing on a linear problem with a single degree of freedom.


\subsection{Long-term behavior} \label{s:longterm}

Consider a simple harmonic oscillator (i.e.~a spring pendulum) with mass $m$ and stiffness $k$. For an initial elongation of $u_0 = u(0)$, the displacement and velocity of the oscillator can be computed analytically: $u_\mathrm{an}(t) = u_0 \, \cos(\omega\,t)$ and $v_\mathrm{an}(t) = - \omega \, u_0 \, \sin(\omega\,t)$. The frequency of oscillation is given by $\omega = \sqrt{k/m}$; the period length (i.e.~the duration of one oscillation) is determined through $T_0 = 2\pi/\omega$. In the following, the numerical results are normalized by $u_0$, $\omega$, $T_0$, and by the initial energy of the system, $E_0 = \frac{1}{2} \, k \, u_0^2$.

We now compare our Hermite schemes with the implicit Newmark algorithm~\cite{newmark59}, choosing the Newmark parameters as $\beta = 1/4$ and $\gamma = 1/2$. Regarding linear systems, this method is then not only unconditionally stable; it can further conserve the energy (see e.g.~Ref.~\cite{hughes76} or \cite{krenk06}). Besides, it is discussed in Ref.~\cite{kane00} that for $\gamma = 1/2$, the Newmark is variational. In addition to the Newmark algorithm, we consider a variational integrator based on linear finite elements in time (see \Appx{a:varintlin}). We will refer to this method as L1-integrator. 

\begin{figure}[ht]
	\subfigure[{Displacement}]{
		\includegraphics[width=0.49\textwidth]{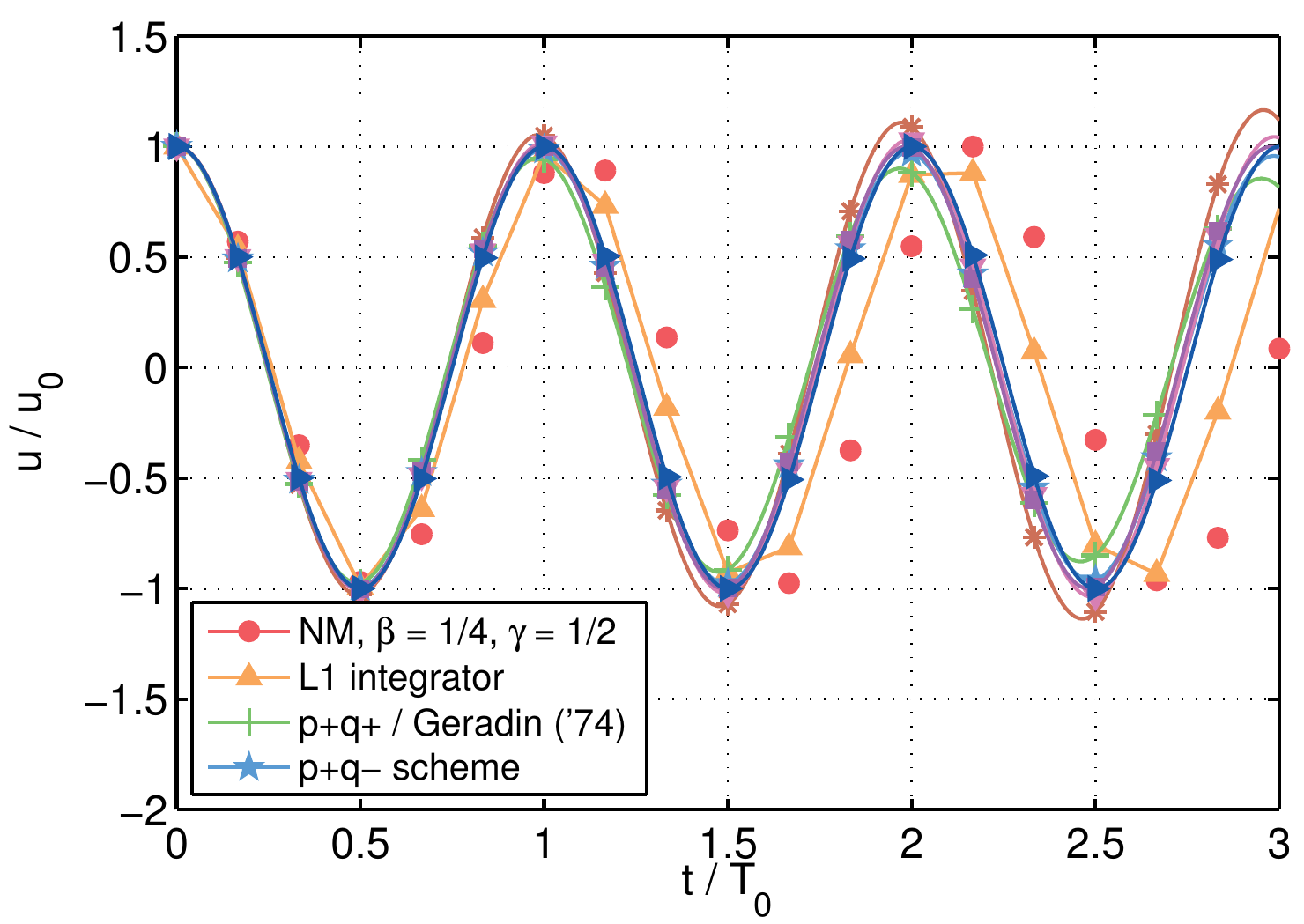}
		\label{f:single:u}}\subfigure[{Velocity}]{
		\includegraphics[width=0.49\textwidth]{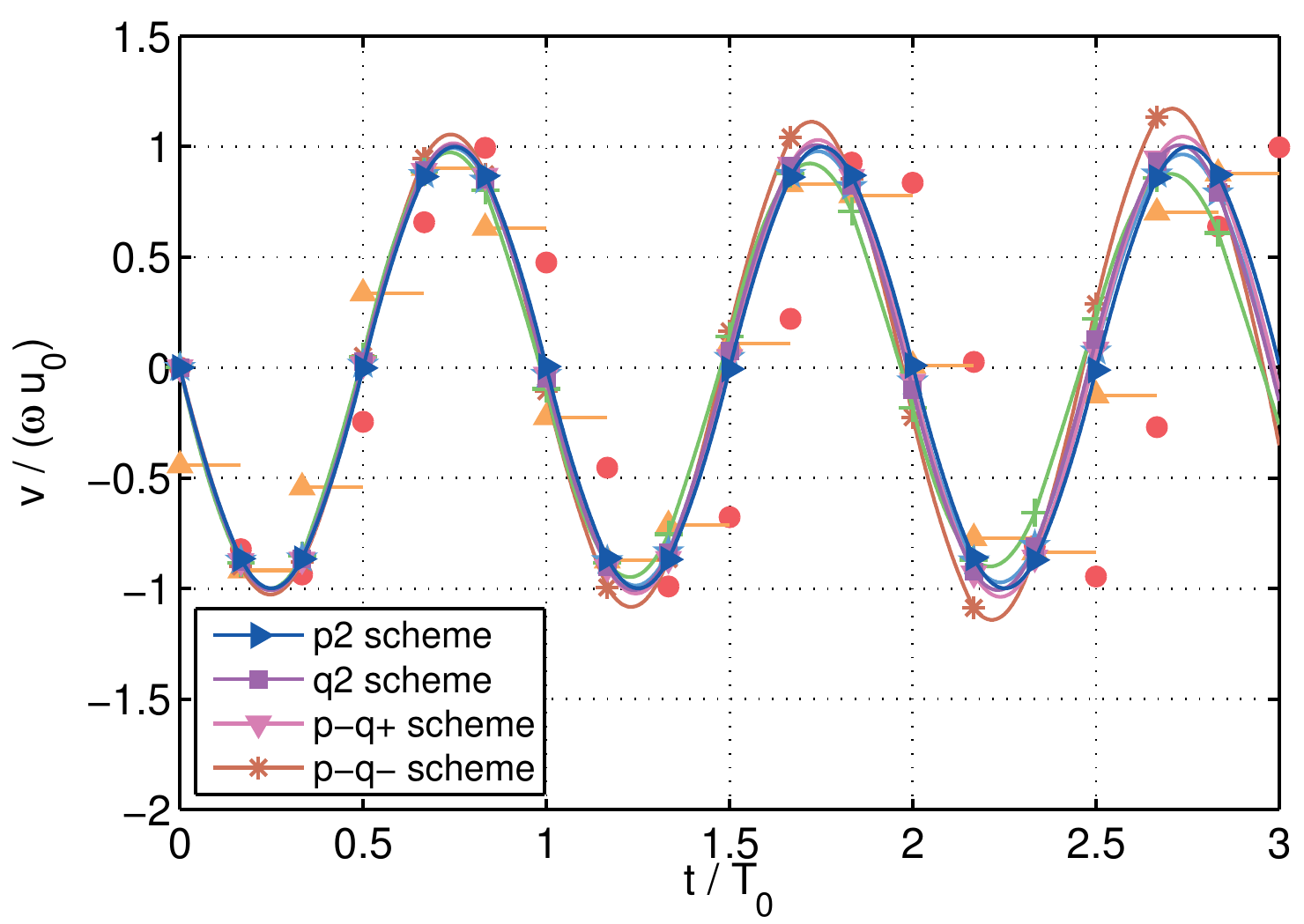}
		\label{f:single:v}}
	\caption{Harmonic oscillator: Displacement and velocity for three
		periods of oscillation; the six Hermite schemes are compared with the
		Newmark algorithm (NM) and with the L1-integrator; $\Dt = T_0/6$.}
	\label{f:single:uv}
\end{figure}

\Fig{f:single:uv} shows the displacement and the velocity of the oscillator for three periods and a very coarse time discretization. As expected, for our six schemes the displacement (velocity) is $C^1$-continuous ($C^0$-continuous) at the discrete time steps. In contrast, the L1-integrator approximates the velocity as a constant along each time interval; this leads to discontinuities at the interval boundaries. For the Newmark algorithm, the displacement and velocity are evaluated only at the discrete time steps, $t_n$. We observe that the oscillation period increases for both Newmark's method and the L1-integrator.

Regarding the maximum displacement, for two of our mixed methods ($\mrp^-\mrq^+$ and $\mrp^-\mrq^-$) the amplitude of oscillation seems to increase remarkably~(\Fig{f:single:u}). This indicates that these methods may be unstable. For the remaining mixed schemes ($\mrp^+\mrq^+$/G{\'e}radin and $\mrp^+\mrq^-$) the amplitudes in both the displacement and velocity decrease. In contrast, both the amplitude and the period of oscillation are well-preserved for our $\mrp2$- and $\mrq2$-schemes.
\begin{figure}[ht]
	\subfigure[{Total energy, $E = K + \Pi$}]{
		\includegraphics[width=0.49\textwidth]{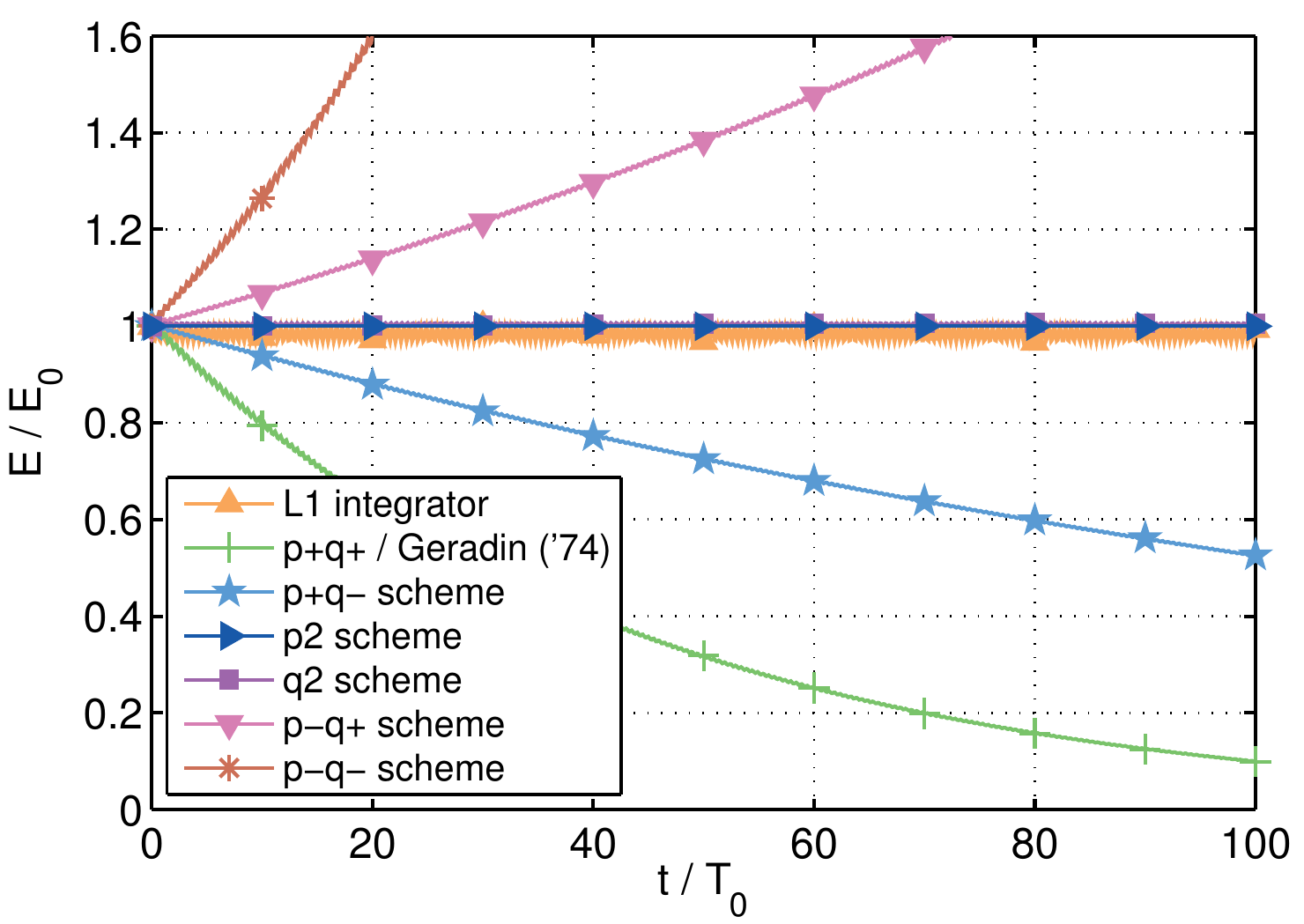}
		\label{f:longtime:E}
	}\subfigure[{Zoom of (a)}]{
		\includegraphics[width=0.49\textwidth]{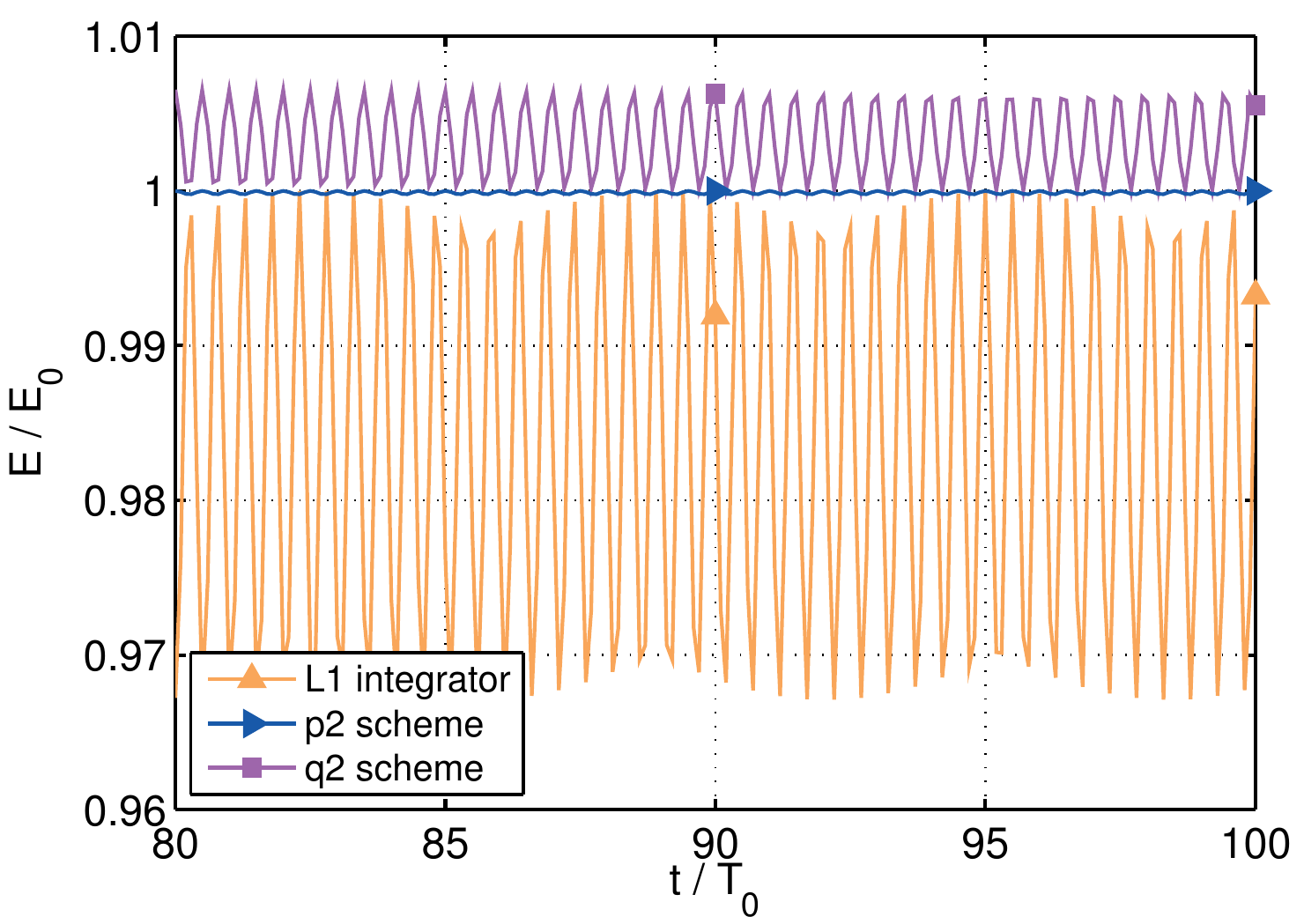}
		\label{f:longtime:Ez}
	}
	\caption{Harmonic oscillator: Long-term energy behavior for
		100 periods of oscillation; the six Hermite schemes are compared with
		the L1-integrator; $\Dt = T_0 / 10$.}
	\label{f:longtime}
\end{figure}

\Fig{f:longtime} shows the total energy of the system over 100 oscillation periods. We observe that the schemes $\mrp^-\mrq^+$ and $\mrp^-\mrq^-$ are unstable, while the schemes $\mrp^+\mrq^+$ and $\mrp^+\mrq^-$ are strongly dissipative. This agrees with the results shown in \Fig{f:single:uv}. As expected, for both the remaining two schemes, $\mrp2$ and $\mrq2$, and for the L1-integrator, the total energy is qualitatively preserved. Interestingly, the $\mrp2$- and the $\mrq2$-scheme are more accurate; this is indicated by smaller amplitudes of oscillation in \Fig{f:longtime:Ez}. Compared to the linear integrator, the maximum relative errors are smaller by one order of magnitude for the $\mrq2$-scheme, and even by two orders for the $\mrp2$-scheme; see \Tab{t:single:eEmax}.
\begin{vchtable}[ht]
	\vchcaption{Harmonic oscillator: Maximum error in the total energy
		for the methods shown in \Fig{f:longtime:Ez}.}
	\label{t:single:eEmax}
	\begin{tabular}{@{}cccc@{}}
	\hline\noalign{\smallskip}
		&	$\mrp2$	&	$\mrq2$	&	$\mrL 1$ \\
	\noalign{\smallskip}\hline\noalign{\smallskip}
		$e_E^\mathrm{max}$	& 0.023\,\%	&	0.653\,\%	&	3.29\,\%	\\
	\hline
	\end{tabular}
\end{vchtable}


\subsection{Stability} \label{s:stab}

\begin{figure}[ht]
	\subfigure[{Spectral radius of $\mA$}]{
		\includegraphics[width=0.485\textwidth]{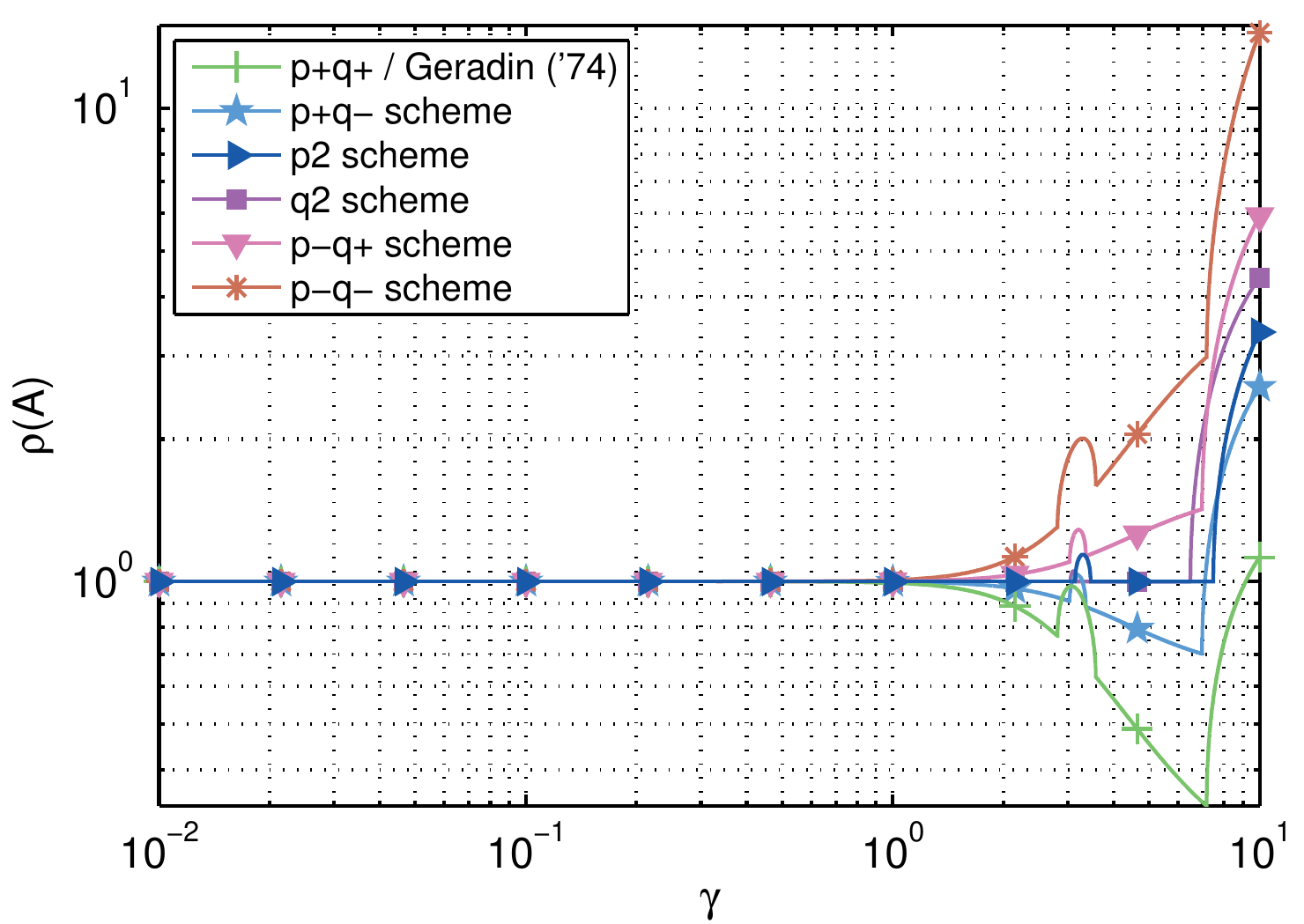}
		\label{f:spectr:all}
	}\subfigure[{Zoom of (a); the purple line is hidden behind the blue one}]{
		\includegraphics[width=0.49\textwidth]{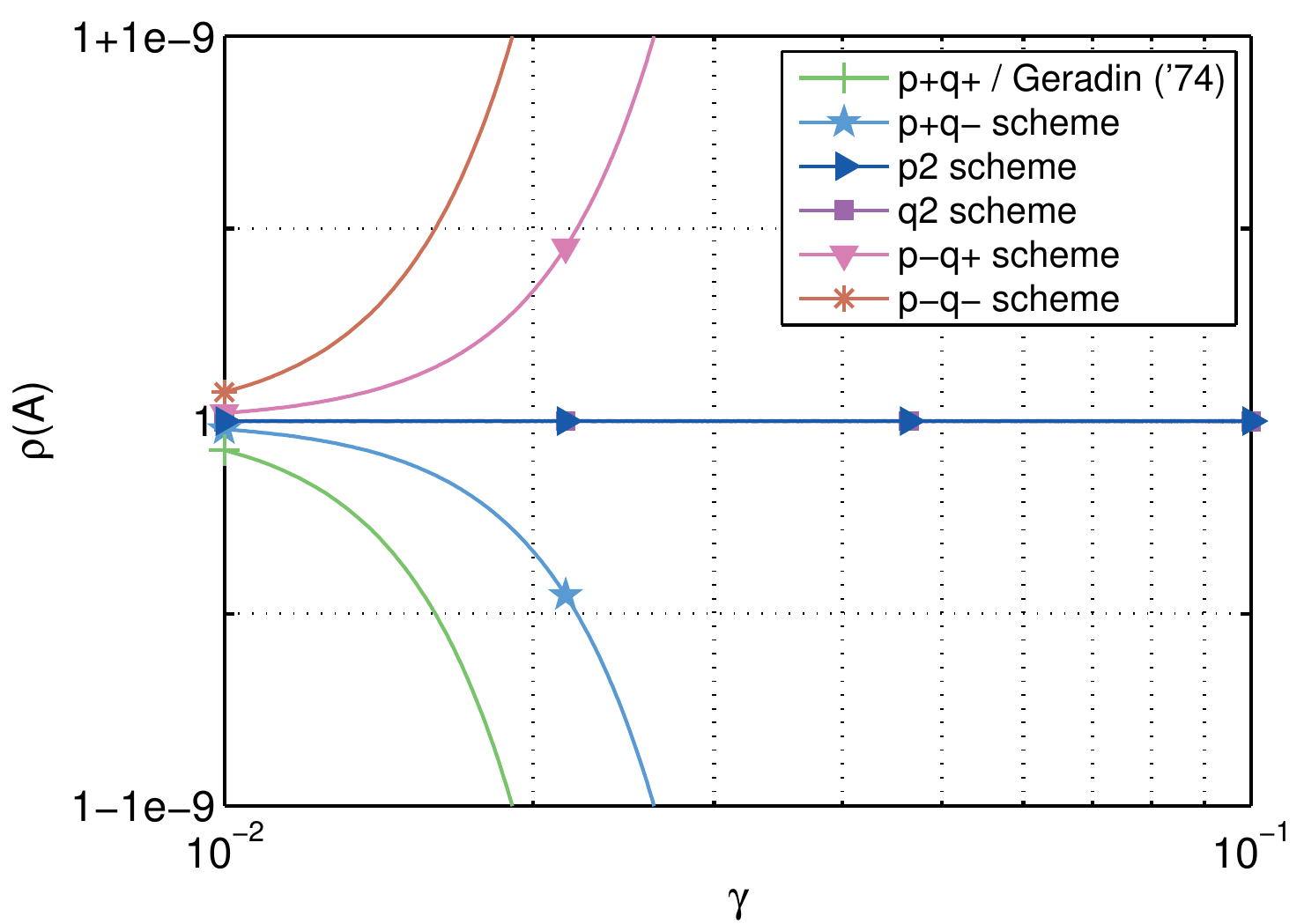}
		\label{f:spectr:z}
	}
	\caption{Harmonic oscillator: Spectral radius for the six Hermite schemes
		as a function of the normalized time step, $\gamma$.}
	\label{f:spectr}
\end{figure}

We now investigate the stability of the six schemes by means of the harmonic oscillator. For this purpose, we introduce the normalized time step $\gamma := \omega \, \Dt$ and insert it into \Eq{e:pp} -- (\ref{e:pmqp}). Following Ref.~\cite{leimkuhler05,oberbloebaum15}, the six schemes can be expressed in the form
\begin{equation}
	\begin{bmatrix} v_{n+1} \\ \omega \, u_{n+1} \end{bmatrix} = \mA
	\begin{bmatrix} v_n \\ \omega \, u_n \end{bmatrix}, \label{e:ampmat}
\end{equation}
where $\mA$ is the amplification matrix given in \Appx{a:ampmat}. The terms $u_n$ and $v_n$ denote the displacement and the velocity at time step $t_n$. \Fig{f:spectr} shows the spectral radius, $\rho(\mA)$, for each of the schemes. \Tab{t:spectr} shows the maximum permitted time step, $\Delta t_\mathrm{stab}$, for which the schemes are stable, i.e.~for which $\rho(\mA) \le 1$. Both the table and \Fig{f:spectr:z} show that the last two schemes, $\mrp^-\mrq^-$ and $\mrp^-\mrq^+$, are unstable even for very small time steps. In contrast, the schemes $\mrp^+\mrq^-$ and $\mrp^+\mrq^+$ seem to be stable for large time steps; nevertheless, these methods are numerically dissipative because of $\rho(\mA) < 1$ for $\gamma < \gamma_\mathrm{stab}$ (see \Fig{f:spectr:z}). The most promising methods seem to be the $\mrp2$-scheme and the $\mrq2$-scheme; they show both excellent stability and energy preservation.
\begin{vchtable}[ht]
	\vchcaption{Harmonic oscillator: Maximum time steps, $\gamma_\mathrm{stab}$
		and $\Dt_\mathrm{stab}$, for which the schemes are stable.}
	\label{t:spectr}
	\begin{tabular}{@{}rcccccc@{}}
	\hline\noalign{\smallskip}
	&	$\mrp2$				&	$\mrq2$				&	$\mrp^+\mrq^-$
	&	$\mrp^+\mrq^+$	&	$\mrp^-\mrq^-$	&	$\mrp^-\mrq^+$ \\
	\noalign{\smallskip}\hline\noalign{\smallskip}
	$\gamma_\mathrm{stab}$	$[-]$	&	3.144	&	3.055	&	3.083	&	9.165
					&	--	&	-- \\[0.3ex]
	$\Dt_\mathrm{stab}$ $[T_0]$	&	0.500	&	0.486	&	0.491	&	1.459
					&	--	&	-- \\
	\hline
	\end{tabular}
\end{vchtable}


\subsection{Symplecticity} \label{s:symp}

As discussed in \Sect{s:longterm}, for the linear oscillator both the $\mrp 2$-scheme and the $\mrq 2$-scheme preserve the energy of the system well. This motivates us to investigate whether these methods are generally symplectic. One way to prove symplecticity is to investigate the derivatives of the phase state, $(\hat{\mpp}_{n+1},\xnpo)$, w.r.t.~the previous state, $(\hat{\mpp}_n,\xn)$. Here, $\hat{\mpp}_\bullet$ is the linear momentum, given by $\hat{\mpp}_\bullet = \mm \, \hat{\mv}_\bullet$. This results in the Jacobian
\begin{align}
	&	\nonumber \\[-3ex]
	\mB_{n+1}	&	:= \begin{bmatrix}
			\ds \pa{\hat{\mpp}_{n+1}}{\hat{\mpp}_n} \,
		&	\ds \, \pa{\hat{\mpp}_{n+1}}{\xn} \\[2ex]
			\ds \pa{\xnpo}{\hat{\mpp}_n} \,
		& \ds \, \pa{\xnpo}{\xn}
	\end{bmatrix} = \begin{bmatrix}
			\ds \mm \, \pa{\vnpo}{\vn} \, \mm^{-1} \,
		&	\ds \, \mm \, \pa{\vnpo}{\xn} \\[2ex]
			\ds \pa{\xnpo}{\vn} \, \mm^{-1} \,
		&	\ds \, \pa{\xnpo}{\xn}
	\end{bmatrix}. \\[-3ex]
	&	\nonumber
\end{align}
According to Ref.~\cite{hairer06}, the mapping $(\hat{\mpp}_n,\xn) \mapsto (\hat{\mpp}_{n+1},\xnpo)$ is symplectic if $\mB_{n+1}$ is symplectic, i.e.~if
\begin{equation}
	\mB^\mrT_{n+1} \, \mJ \, \mB_{n+1} = \mJ, \qquad
	\mJ = \begin{bmatrix}
		\bzero 		& \bI_{d\cdot\nno} \\
		-\bI_{d\cdot\nno} & \bzero
	\end{bmatrix}, \qquad n = 0, \dots, N\!-\!1, \label{e:defsymp}
\end{equation}
where $\bI_{d\cdot\nno}$ is the identity matrix of dimension $d\cdot\nno$. For the harmonic oscillator, the Jacobian reduces to
\begin{equation}
	\mB_{n+1} = \begin{bmatrix}
	\ds \pa{v_{n+1}}{v_n}								&	\ds m \, \pa{v_{n+1}}{u_n} \\[2ex]
	\ds \frac{1}{m} \, \pa{u_{n+1}}{v_n}	&	\ds \pa{u_{n+1}}{u_n} \end{bmatrix}.
\end{equation}
In this case one can show that the determinant of $\mB_{n+1}$ is equal to the determinant of the amplification matrix introduced in \Eq{e:ampmat}: $\det\,(\mB_{n+1}) = \det\,(\mA)$. With this relation the condition for symplecticity~\eqref{e:defsymp} is fulfilled if the determinant of~$\mA$ is equal to one,
\begin{equation}
	\begin{bmatrix} 0 & \det\,(\mA) \\ -\det\,(\mA) & 0 \end{bmatrix}
	= \begin{bmatrix} \ \ 0 & 1 \\ -1 & 0 \end{bmatrix} \qquad
	\Rightarrow \qquad \det\,(\mA) = 1.
\end{equation}
By means of \Appx{a:ampmat} one can show that both our $\mrp 2$-scheme and $\mrq 2$-scheme fulfill $\det\,(\mA) = 1\ \ \forall\ \gamma$. This implies that --- at least for the harmonic oscillator --- these schemes are symplectic. It remains to be subject of further investigation whether this is also true for arbitrary systems with multiple degrees of freedom.


\subsection{Convergence behavior} \label{s:conv:single}

\begin{figure}[ht]
	\subfigure[Maximum error in the displacement]{
		\includegraphics[width=0.49\textwidth]
		{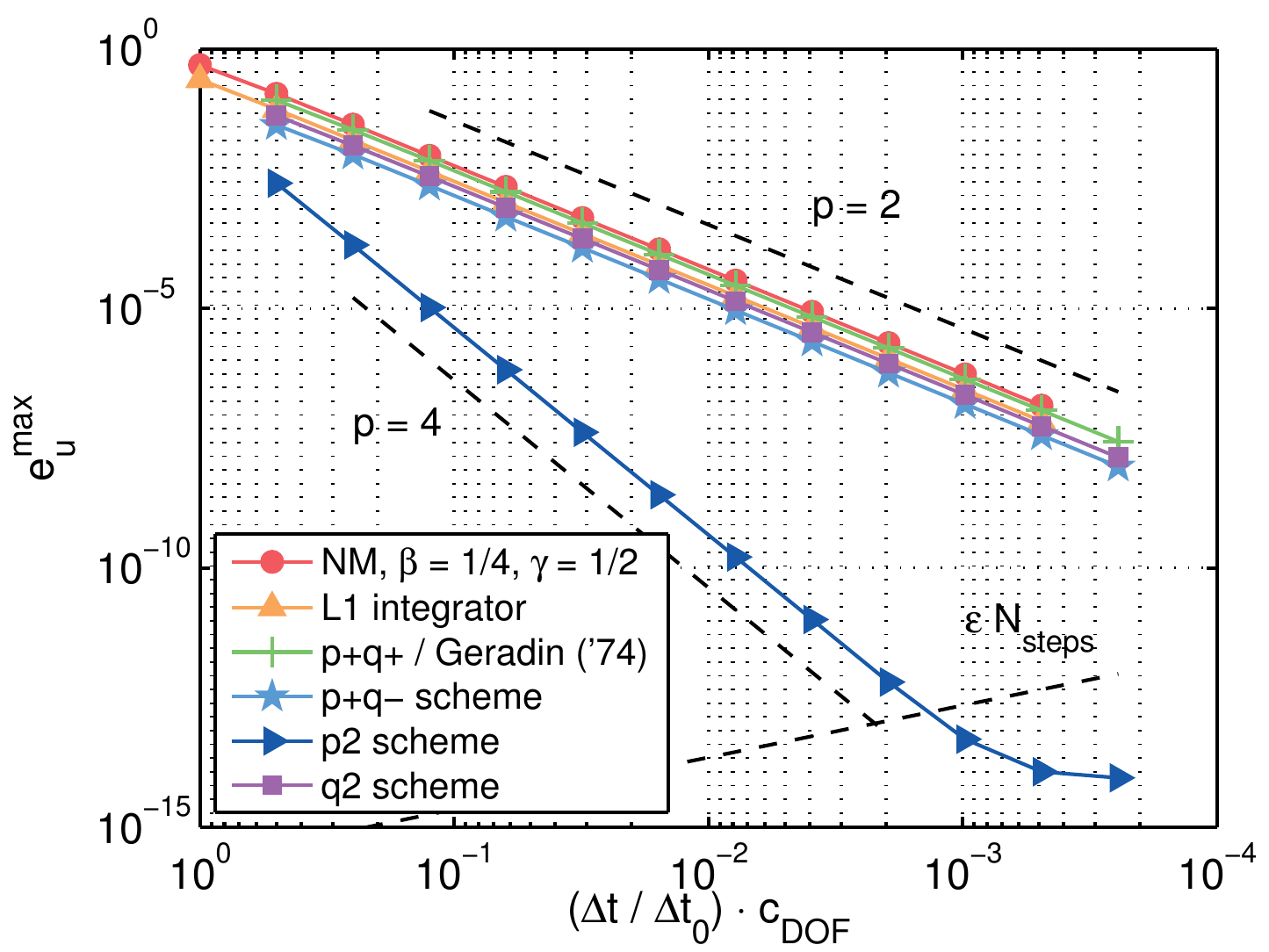}\label{f:single:eu_max}}
	\subfigure[Maximum error in the velocity]{
		\includegraphics[width=0.49\textwidth]
		{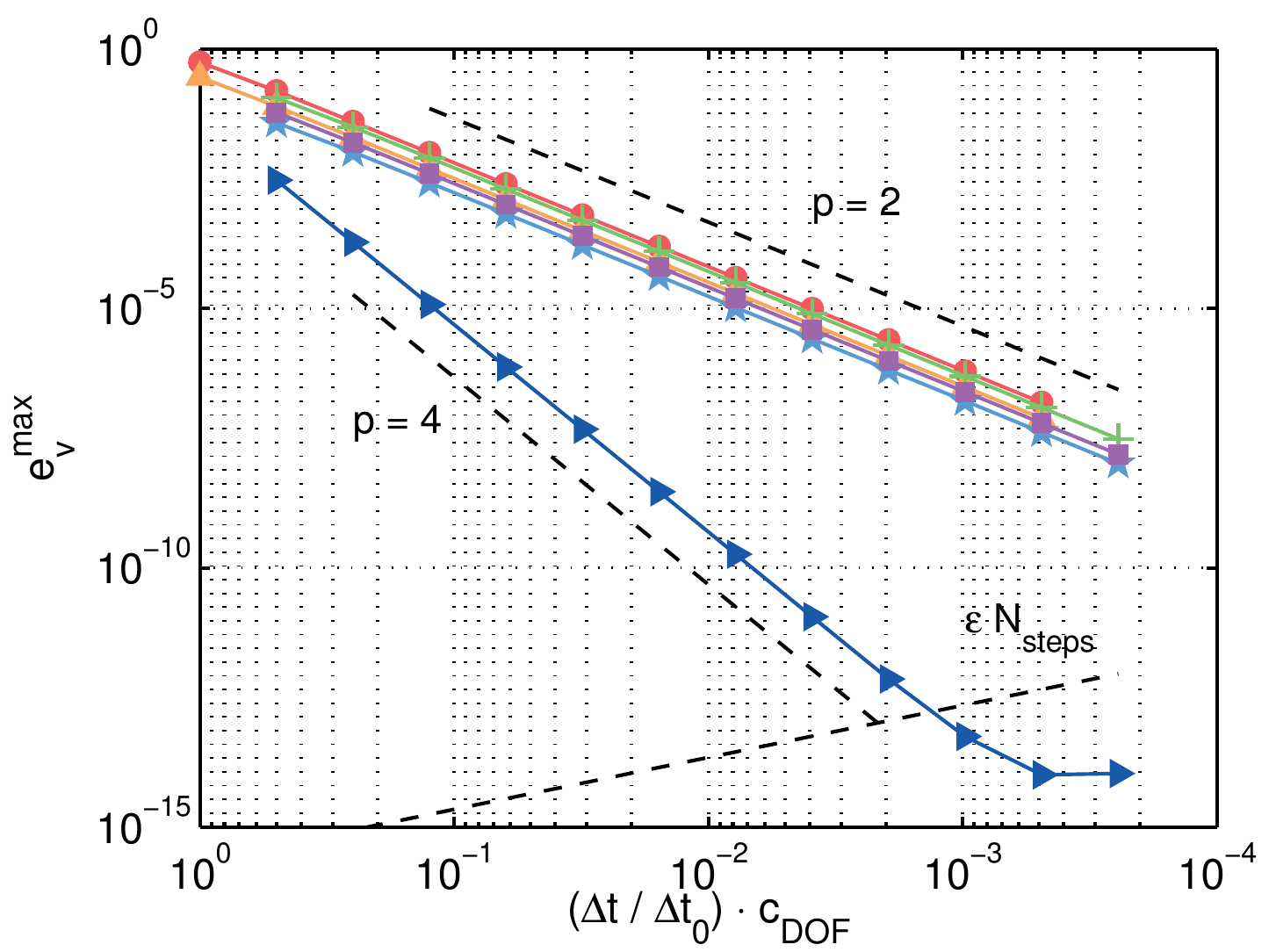}\label{f:single:ev_max}}
	\subfigure[Maximum error in the energy]{
		\includegraphics[width=0.49\textwidth]
		{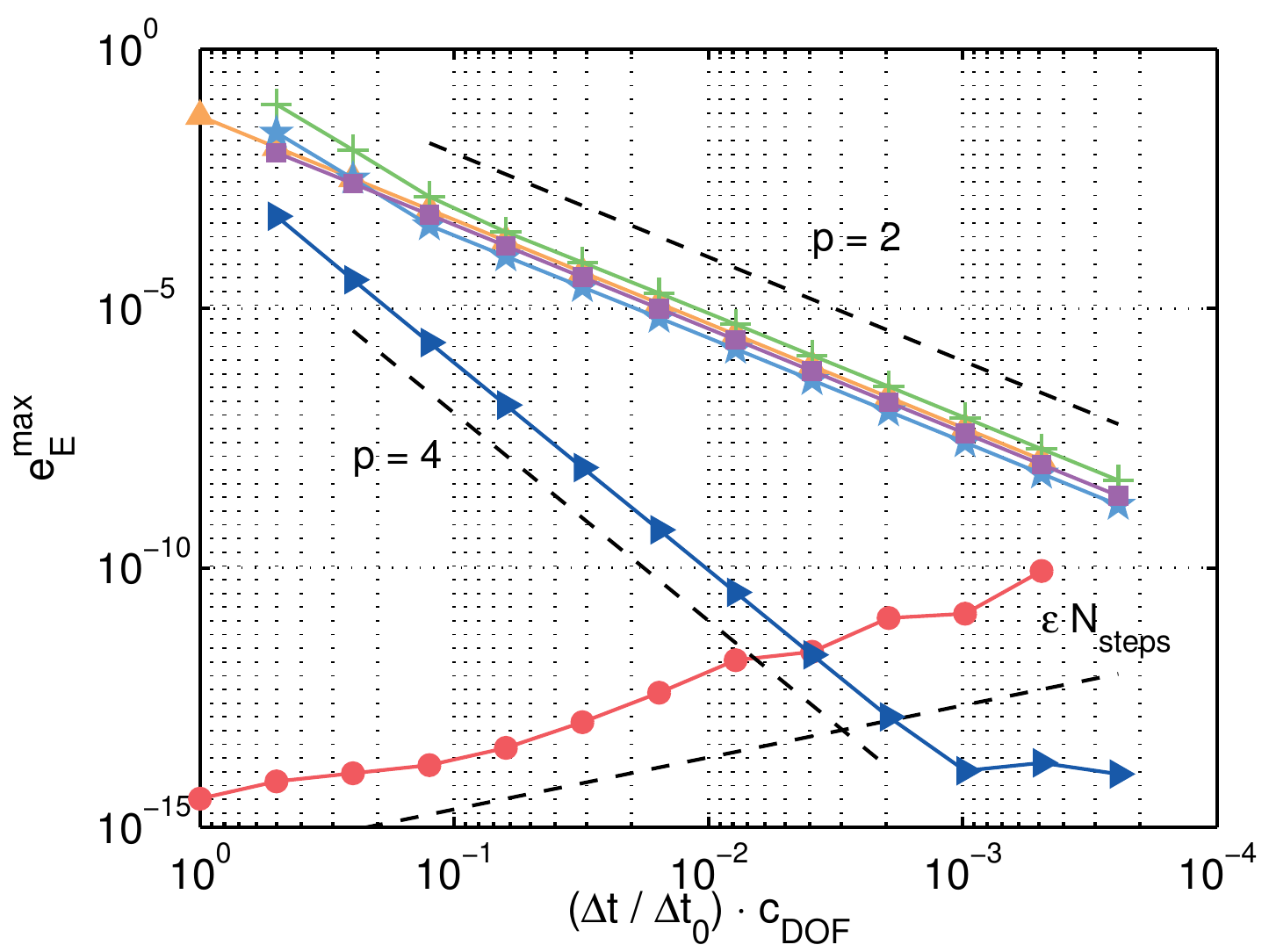}\label{f:single:eE_max}}
	\caption{Harmonic oscillator: Convergence behavior for
		$T = 2 \, T_0$, $\Dt_0 = T_0 / 8$; $c_\mathrm{DOF} = 1/2$ for our
		Hermite schemes, otherwise $c_\mathrm{DOF} = 1$; the dashed line labeled
		with $\varepsilon \, N_\mathrm{steps}$ denotes the estimated machine
		precision multiplied with the number of computed steps.}
	\label{f:single:euvE}
\end{figure}

We now focus on the four stable schemes: $\mrp2$, $\mrq2$, $\mrp^+\mrq^-$, and $\mrp^+\mrq^+$ (which is equivalent to the method of G\'eradin~\cite{geradin74}). In order to study convergence for the harmonic oscillator we consider the maximum errors of the displacement, velocity, and total energy at the discrete time steps; these are given by
\begin{equation}
	e_\bullet^\mathrm{max} = \max_{n = 0,\dots, N} | e_\bullet(t_n) |,
	\qquad \bullet \in \{u,v,E\},
\end{equation}
where
\begin{equation}
	e_u(t) = |u(t) - u_\mathrm{an}(t)| \, / \, |u_0|, \qquad 
	e_v(t) = |v(t) - v_\mathrm{an}(t)| \, / \, |\omega \, u_0|, \qquad
	e_E(t) = |E(t) - E_0| \, / \, E_0. \label{e:eEconv}
\end{equation}
\Fig{f:single:euvE} shows the convergence behavior of these errors for our stable schemes, the Newmark algorithm, and the L1-integrator. In addition, we account for the results discussed in Ref.~\cite{leok12} for cubic Hermite interpolation. Note that compared to the last three methods, our time integration schemes must account for twice the number of unknowns in each step: the nodal displacements and the nodal velocities. This has been considered in the scaling of the abscissae by introducing the factor $c_\mathrm{DOF}$ ($c_\mathrm{DOF} = 1/2$ for our schemes, otherwise $c_\mathrm{DOF} = 1$).

As already shown in \Fig{f:single:v}, for the L1-integrator the piecewise approximation of the velocity is discontinuous at the discrete time steps. Nevertheless, we can determine the discrete velocity at $t_n$ by computing the momentum from a discrete Legendre transformation~\cite{hairer06}. This approach is also discussed by Ober-Bl{\"o}baum et al.~\cite{oberbloebaum11} considering boundary conditions for the velocity. For the L1-integrator the maximum errors in both the velocity and in the energy (Fig.~\ref{f:single:ev_max} and \ref{f:single:eE_max}) are finally obtained by
\begin{equation}
	\vn = \mm^{-1} \, \hat{\mpp}_n, \qquad \hat{\mpp}_n :=
	- \pa{\A^{ht}_{n+1}(\xn,\xnpo)}{\xn} = \pa{\A^{ht}_n(\xnmo,\xn)}{\xn}.
\end{equation}

As expected, Newmark's method can conserve the energy of the (linear) system; the corresponding error thus lies in the range of machine precision. The errors in the displacement and in the velocity, however, are of order~$p = 2$. The same rate of convergence can be observed for the L1-integrator. This agrees with the discussion in Ref.~\cite{oberbloebaum15} that integrators interpolating the displacement linearly in time can be at most of second order.

Our schemes $\mrp^+\mrq^+$~\cite{geradin74}, $\mrp^+\mrq^-$, and $\mrq2$ also show an order of $p = 2$. Since their computational effort is higher than for both Newmark's method and the L1-integrator, these methods are not favorable. In comparison, Leok and Shingel prove convergence with order three for their integrator based on cubic Hermite interpolation (i.e.~$n = 2$ in Theorem~2 of Ref.~\cite{leok12}). An even better rate of convergence is achieved with our favorite candidate, the $\mrp2$-scheme: All of the errors in the displacement, the velocity, and the energy converge with order $p = 4$.


\section{Numerical results for 1D elastodynamics} \label{s:results}

As shown in the previous sections, our $\mrp2$-scheme is symplectic for the harmonic oscillator; it further possesses the highest rate of convergence. So far we have studied only linear problems with a single degree of freedom. In the following section we apply the $\mrp2$-scheme to spatially continuous problems.


\subsection{Axial vibration of a linear elastic bar} \label{s:freevibr}

We now discuss the free axial vibration of a linear elastic bar  (\Appx{a:bar}). For the spatial discretization we either use linear Lagrange or cubic Hermite finite elements; see~\Appx{a:shpfct} and text books such as Ref.~\cite{wriggers08}. The second type of element yields a $C^1$-continuous approximation of the displacement also in space.

Consider the bar vibrating in the first (i.e.~the lowest) characteristic eigenmode. For this test case, the corresponding displacement, velocity, and energy can be analytically computed from the one-dimensional wave equation; this yields
\begin{equation}
	u_\mathrm{an}(X,t) = u_0 \cdot \cos \,\bigg(\frac{\pi\,X}{L}\bigg)
	\cdot \cos\,(\omega_\mathrm{an}\,t), \qquad
	v_\mathrm{an}(X,t) = -u_0 \cdot \cos\,\bigg(\frac{\pi \, X}{L}\bigg)
	\cdot \omega_\mathrm{an} \, \sin\left({\omega_\mathrm{an}\,t}\right),
	\label{e:bar:uvan}
\end{equation}
and $E_0 = EA/(4\,L) \cdot (\pi\,u_0)^2$. Here, $X \in [0,L]$, $u_0$ is the amplitude of oscillation, and $\omega_\mathrm{an} = \pi \sqrt{ E / (\rho_0 \, L^2)}$ is the first natural frequency. The deformation of the bar vibrating in the first mode is shown in \Fig{f:L1bar:deform}. Here, six linear elements are used for the spatial discretization. As expected, the bar performs sinusoidal oscillations. Due to the coarse finite element mesh, however, the structure oscillates with a frequency slightly higher than the analytical solution: $|\omega - \omega_\mathrm{an}| \, / \, \omega_\mathrm{an} \approx 1.15\,\%$. The oscillation period of the discrete system, $T_0$, is thus smaller than the analytical solution, $T_0^\mathrm{an} = 2\pi / \omega_\mathrm{an}$.
\begin{vchfigure}[ht]
	\includegraphics[width=0.48\textwidth]{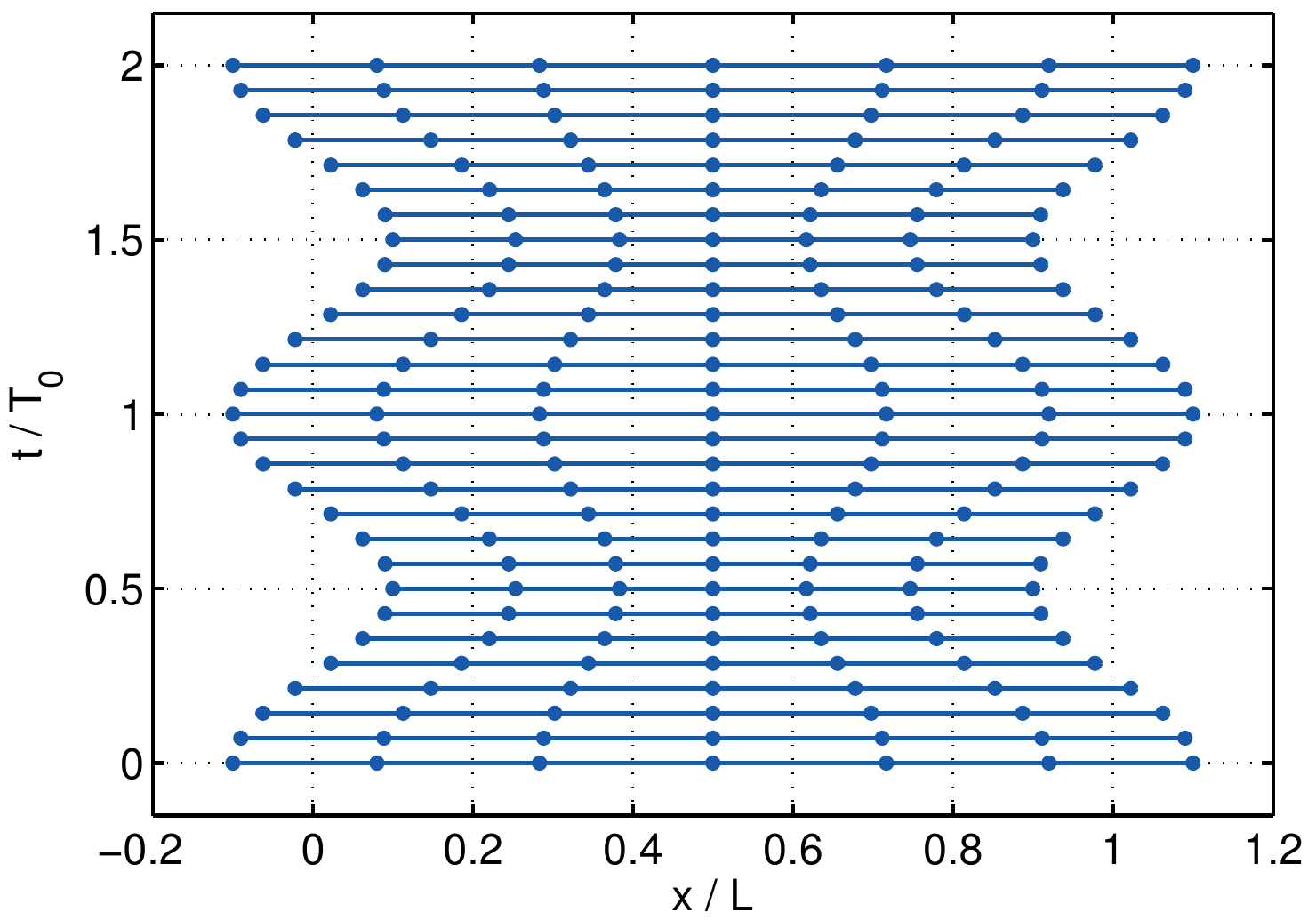}
	\vchcaption{Linear elastic bar: Two oscillations in the
		first natural frequency using the $\mrp2$-scheme and linear FE.}
	\label{f:L1bar:deform}
\end{vchfigure}

Like variational integrators, our $\mrp2$-scheme does not incorporate any numerical dissipation to damp spurious oscillations; we thus must carefully adjust the size of the time step, $\Dt$, to the spatial FE element mesh. Since a stability analysis including spatial discretization can be quite tedious, we roughly estimate the maximum permitted time step for either linear or Hermite elements. We therefore consider the Courant-Friedrichs-Lewy (CFL) condition for one-dimensional problems,
\begin{equation}
	C_\mathrm{CFL} \le C^\mathrm{max}, \qquad C_\mathrm{CFL} :=
	\frac{c_\mathrm{W} \cdot \Dt}{\Delta L}, \label{e:cflcond}
\end{equation}
where $c_\mrW = \sqrt{E/\rho_0}$ is the velocity of wave propagation, and $\Delta L$ is the characteristic discretization length. For a linearly interpolated element applies $\Delta L := L_e$; for a Hermite finite element, we choose $\Delta L := L_e/2$ to take into account that it has twice the number of unknowns (and thus higher accuracy).

We now vary the CFL number for a bar vibrating for at least 1000 oscillations. The estimated maximum values are useful to choose appropriate parameters for the following numerical examples. Of course we cannot ensure, however, stability for arbitrary CFL numbers smaller than these estimates. As shown in the stability analysis for a single degree of freedom (\Sect{s:stab}), the methods may also become unstable for small ranges of parameters. This becomes apparent in \Fig{f:spectr:all}, where the spectral radius of the $\mrp2$-scheme exceeds the limit~(one) for a small range of time steps, while being stable again for larger steps. Apart from that, the CFL condition does not serve as sufficient condition for stability. For linear finite elements we estimate $C^\mathrm{max}_{\mrL 1} \approx 1.00$ and $C^\mathrm{max}_{\mrp2} \approx 0.90$; this implies that the time step should fulfill $\Dt \le 1.00 \, \Delta L/c_\mrW$ for the L1-integrator, and $\Dt \le 0.90 \, \Delta L/c_\mrW$ for the $\mrp2$-scheme. For Hermite finite elements we obtain $C^\mathrm{max}_{\mrL 1} \approx 0.72$ and $C^\mathrm{max}_{\mrp2} \approx 0.96$. In the case of linear problems Newmark's method is unconditionally stable if $\beta = 1/4$ and $\gamma = 1/2$.


\subsection{Convergence for the linear bar} \label{s:freevibr:conv}

\begin{figure}[ht]
	\vspace*{-1ex}
	\subfigure[{Error in the displacement for spatially linear FE}]{
		\includegraphics[width=0.49\textwidth]
			{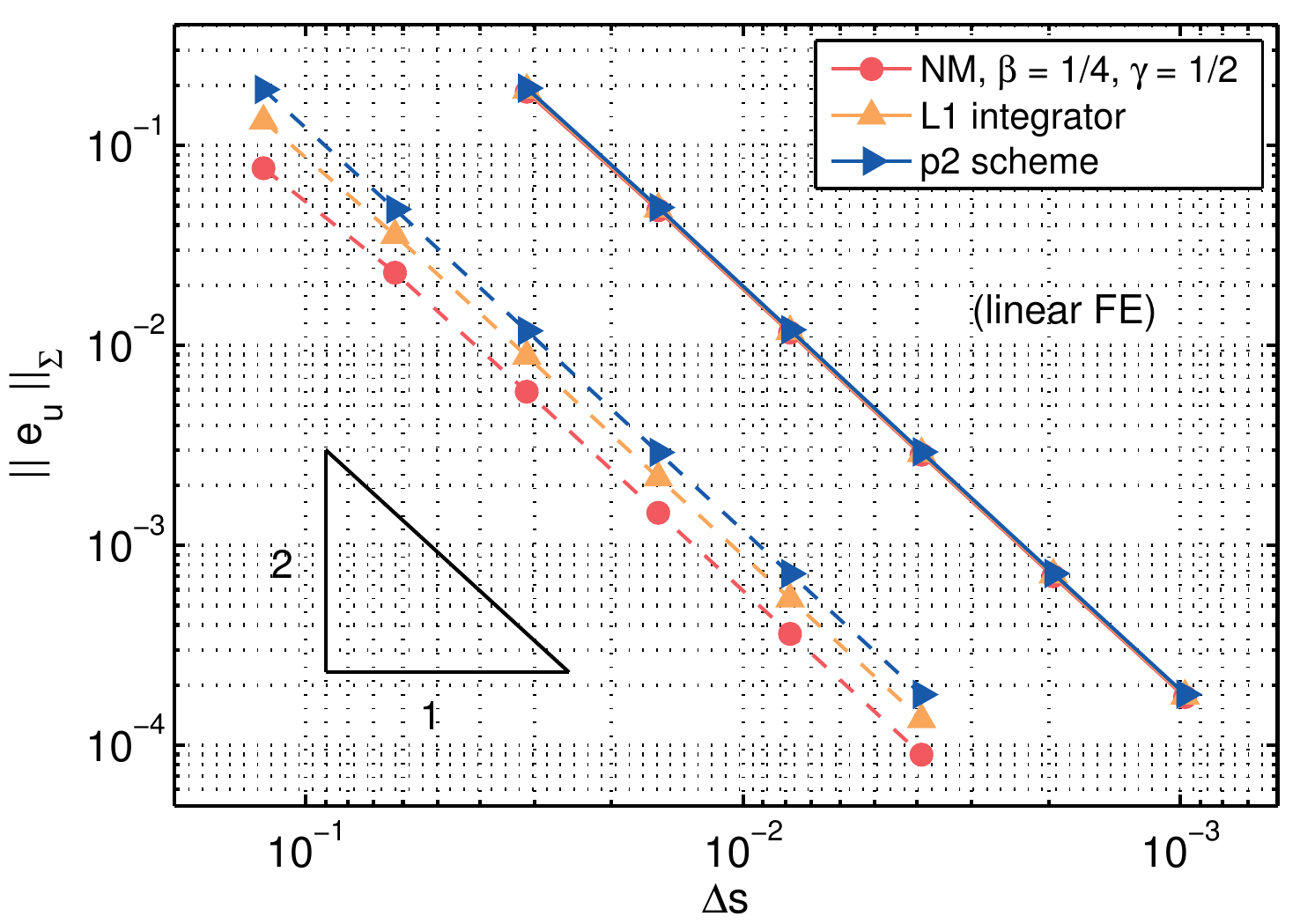}\label{f:L1bar:eu_CFL}}
	\subfigure[{Error in the displacement for spatial Hermite FE}]{
		\includegraphics[width=0.49\textwidth]
			{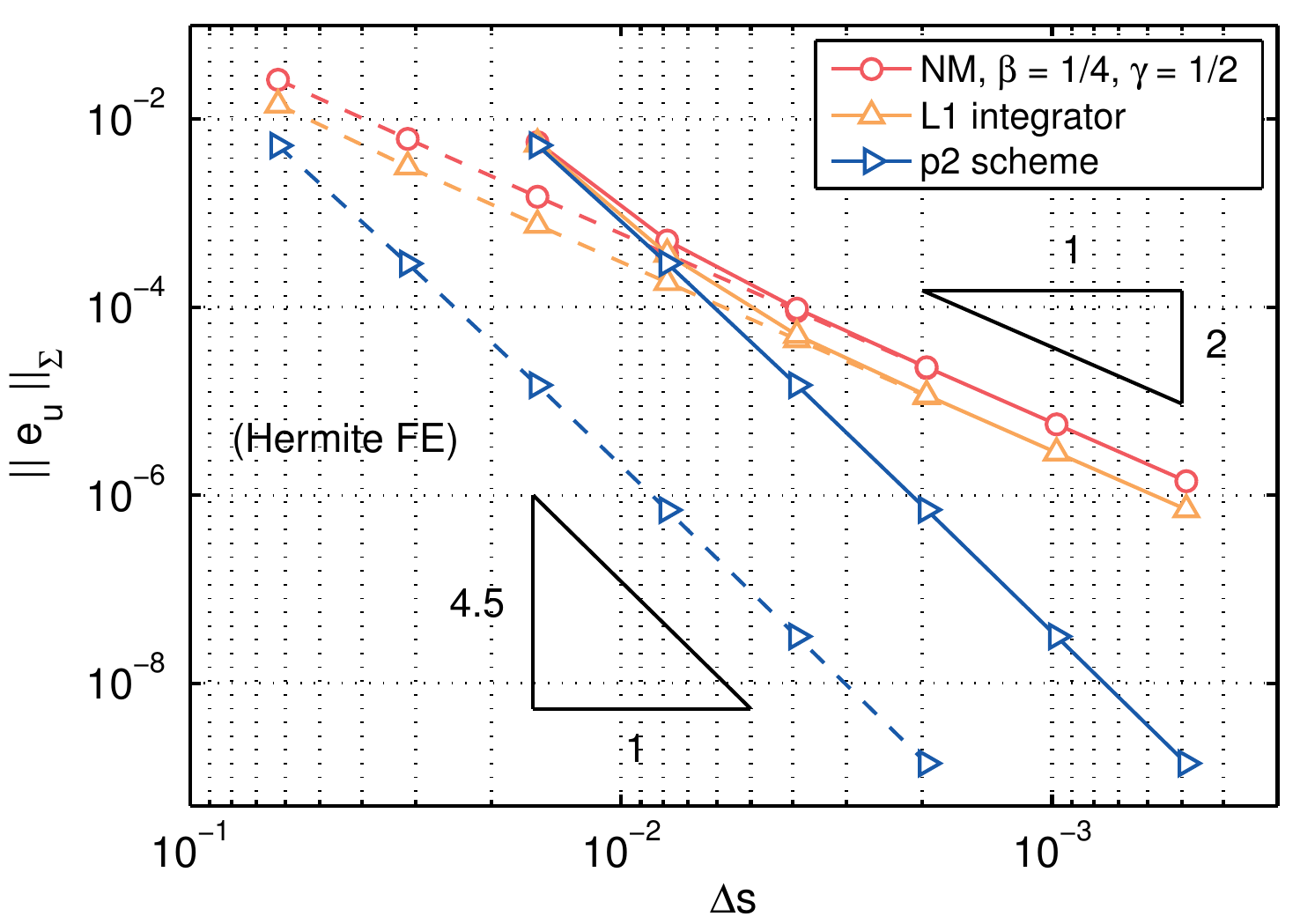}\label{f:H3bar:eu_CFL}}
		\subfigure[{Error in the velocity for spatially linear FE}]{
		\includegraphics[width=0.49\textwidth]
			{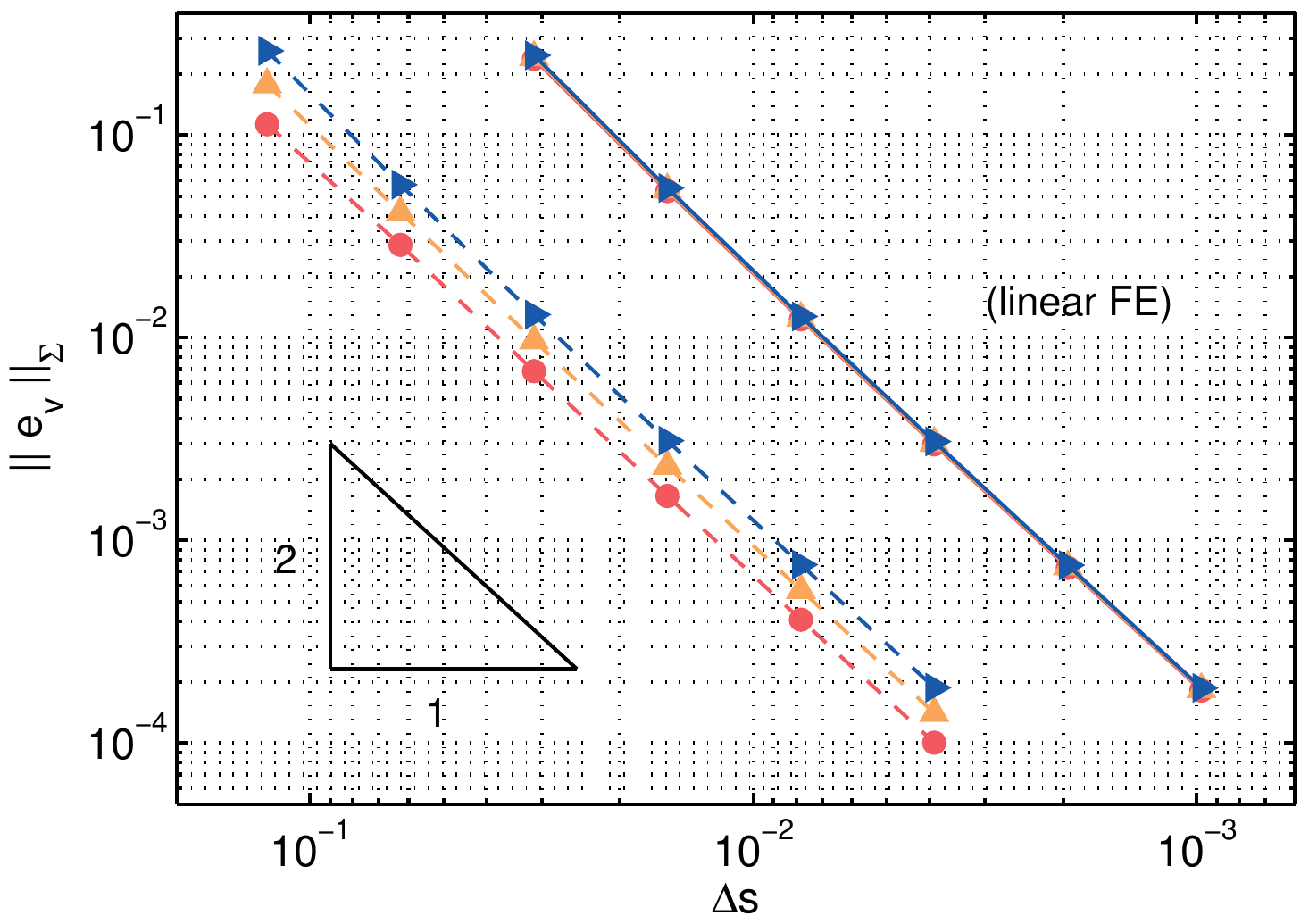}\label{f:L1bar:ev_CFL}}
	\subfigure[{Error in the velocity for spatial Hermite FE}]{
		\includegraphics[width=0.49\textwidth]
			{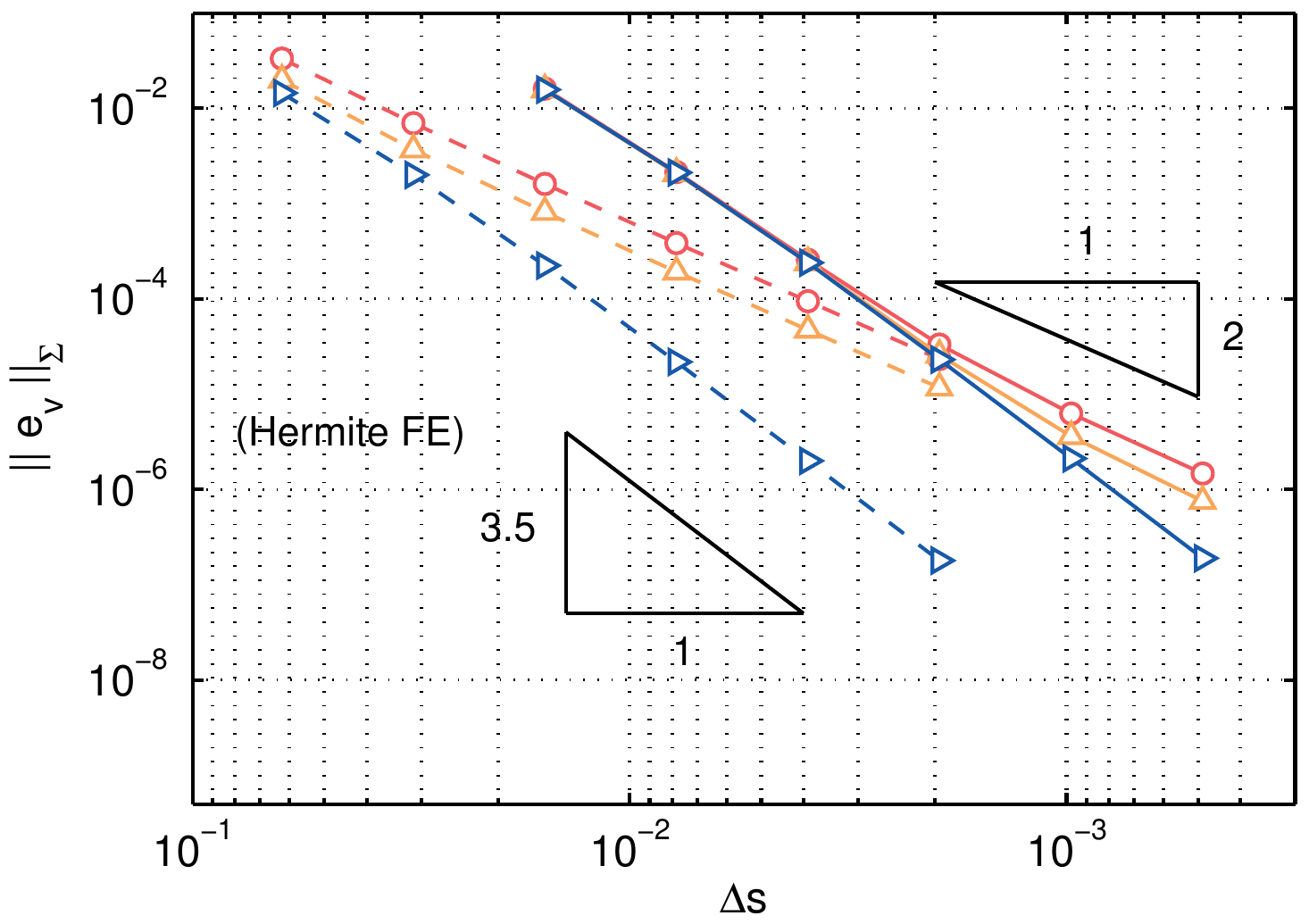}\label{f:H3bar:ev_CFL}}
	\subfigure[{Error in the energy for spatially linear FE}]{
		\includegraphics[width=0.49\textwidth]
		{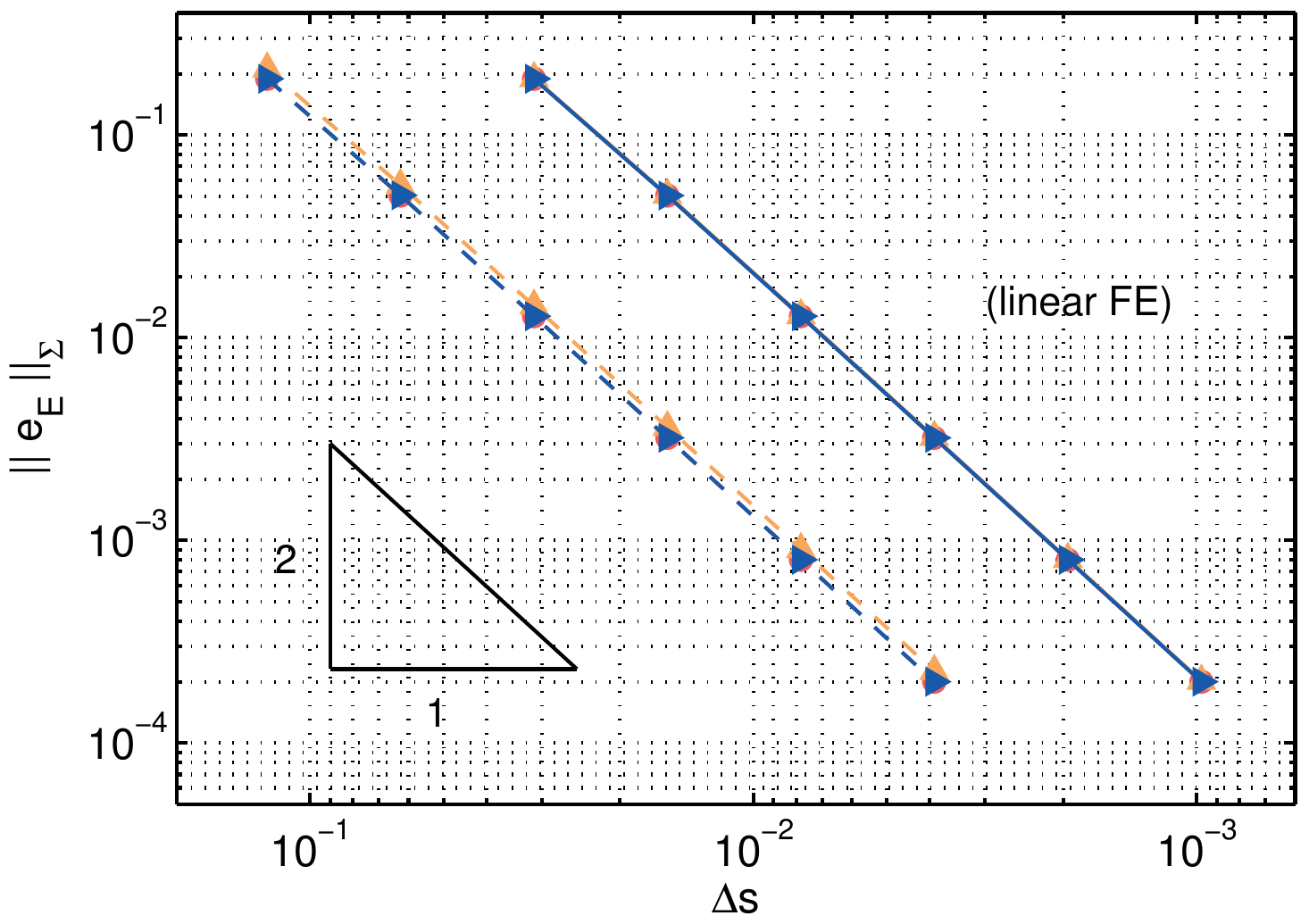}\label{f:L1bar:eE_CFL}}
	\subfigure[{Error in the energy for spatial Hermite FE}]{
		\includegraphics[width=0.49\textwidth]
		{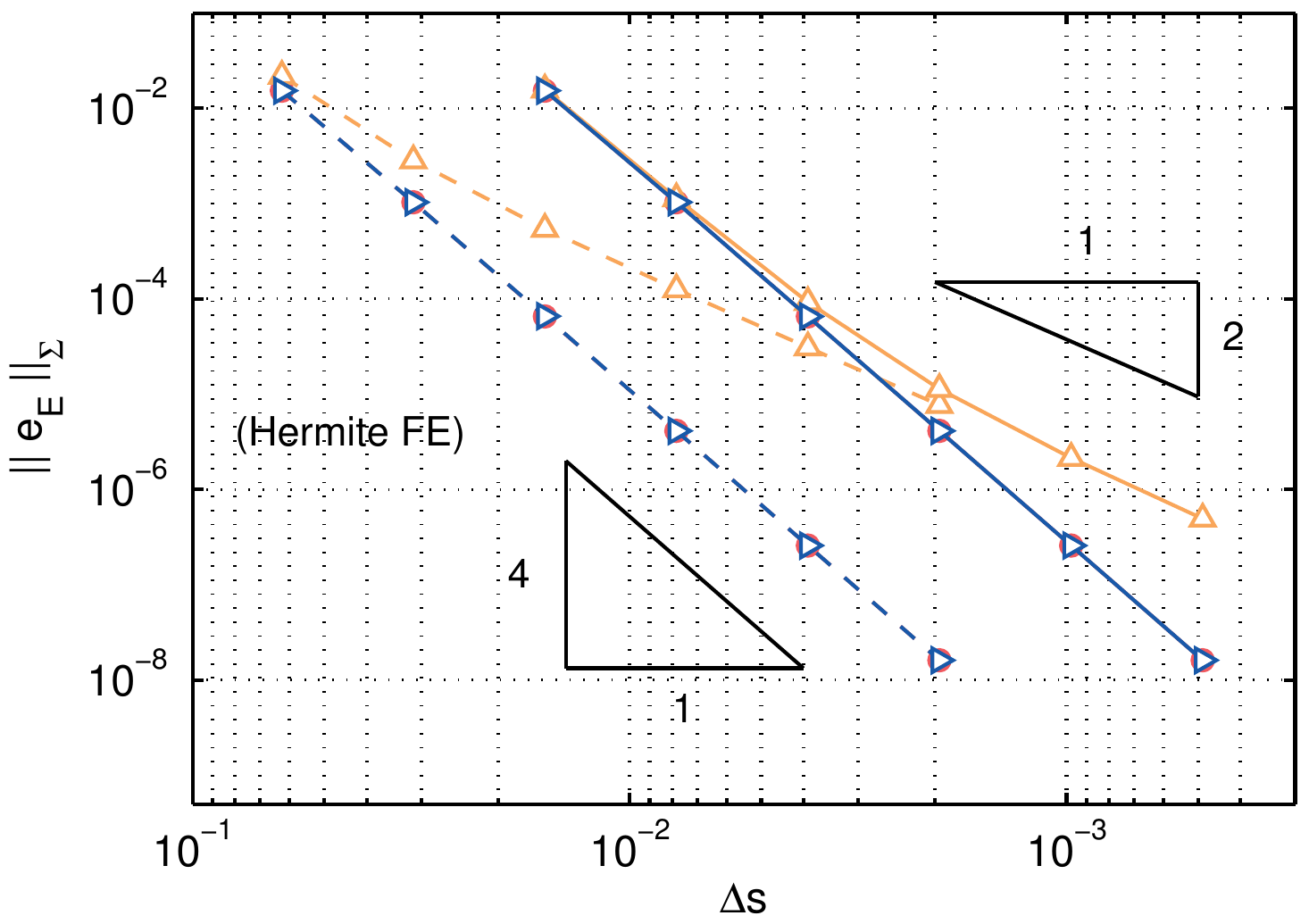}\label{f:H3bar:eE_CFF}}
	\vspace*{-1ex}
	\caption{Linear elastic bar: Convergence behavior of the displacement,
		the velocity, and the energy refining both the mesh and time step;
		the parameter $\Delta s$ is given by $\Delta s
		= \Delta t / T_0^\mathrm{an} = C_\mathrm{CFL} \cdot \Delta L/(2 \, L)$;
		dashed line: $C_\mathrm{CFL} = 0.5$, solid line:
		$C_\mathrm{CFL} = 0.125$; $T = 1 \, T_0^\mathrm{an}$.}
	\label{f:CFL:euvE}
\end{figure}

Reconsider the axial vibration discussed in the previous section. Following Ref.~\cite{demoures15}, we introduce for the displacement and velocity discrete $L^2$-norms including the relative errors at all time steps and finite element nodes:
\begin{equation}
	{||e_\bullet||}_{\Sigma} := \sqrt{ \sum_{n=0}^N \, \sum_{I=1}^\nno \,
	\frac{{|e_\bullet(X_I,t_n)|}^2}{(N\!+\!1) \cdot \nno} \, }, \qquad
	\bullet \in \{u,v\}. \label{e:L2discr}
\end{equation}
This corresponds to the Frobenius norm of the arrays $e_u(X_I,t_n)$ and $e_v(X_I,t_n)$ normalized by (the square roots of) the numbers of nodes and time steps. The relative errors are defined as
\begin{equation}
	e_u(X,t) = |u(X,t) - u_\mathrm{an}(X,t)| \,/\, |u_0|, \qquad
	e_v(X,t) = |v(X,t) - v_\mathrm{an}(X,t)| \,/\, |\omega_\mathrm{an}\,u_0|.
	\label{e:euvdiscr}
\end{equation}
In analogy to \Eq{e:L2discr}, we define a discrete $L^2$-norm for the error in the energy,
\begin{equation}
	{||e_E||}_{\Sigma} := \sqrt{ \sum_{n=0}^N \,
	\frac{{|e_E(t_n)|}^2}{N\!+\!1} \, },
\end{equation}
using $e_E(t)$ from \Eq{e:eEconv}. \Fig{f:CFL:euvE} shows the convergence behavior of the displacement, the velocity, and the energy. For the spatial discretization, we either use linear and Hermite elements. We further consider two fixed CFL numbers, while refining both the mesh and the time step simultaneously.

For a linear finite element mesh (left column of \Fig{f:CFL:euvE}), the three time discretization methods converge with the same order. This indicates that for this specific problem, the error caused by the spatial discretization predominates. In contrast, the error due to the spatial Hermite discretization (right column of \Fig{f:CFL:euvE}) carries considerably less weight. The $\mrp2$-scheme yields a significantly higher convergence than for both the Newmark algorithm and the L1-integrator.

\begin{figure}[ht]
	\subfigure[{Temporal error in the displacement}]{
		\includegraphics[width=0.49\textwidth]
			{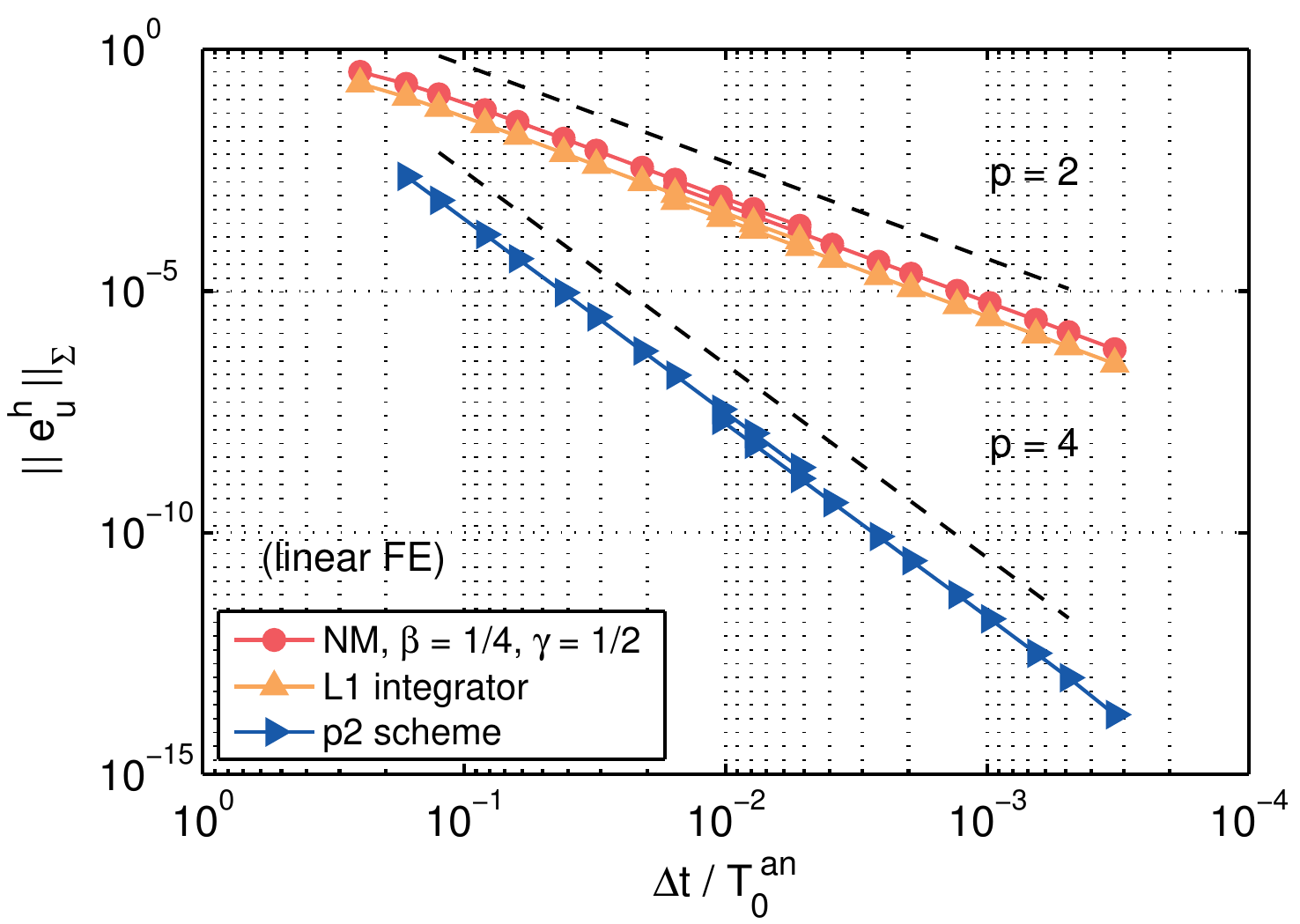}\label{f:L1bar:euh_F}}
	\subfigure[{Temporal error in the velocity}]{
		\includegraphics[width=0.49\textwidth]
			{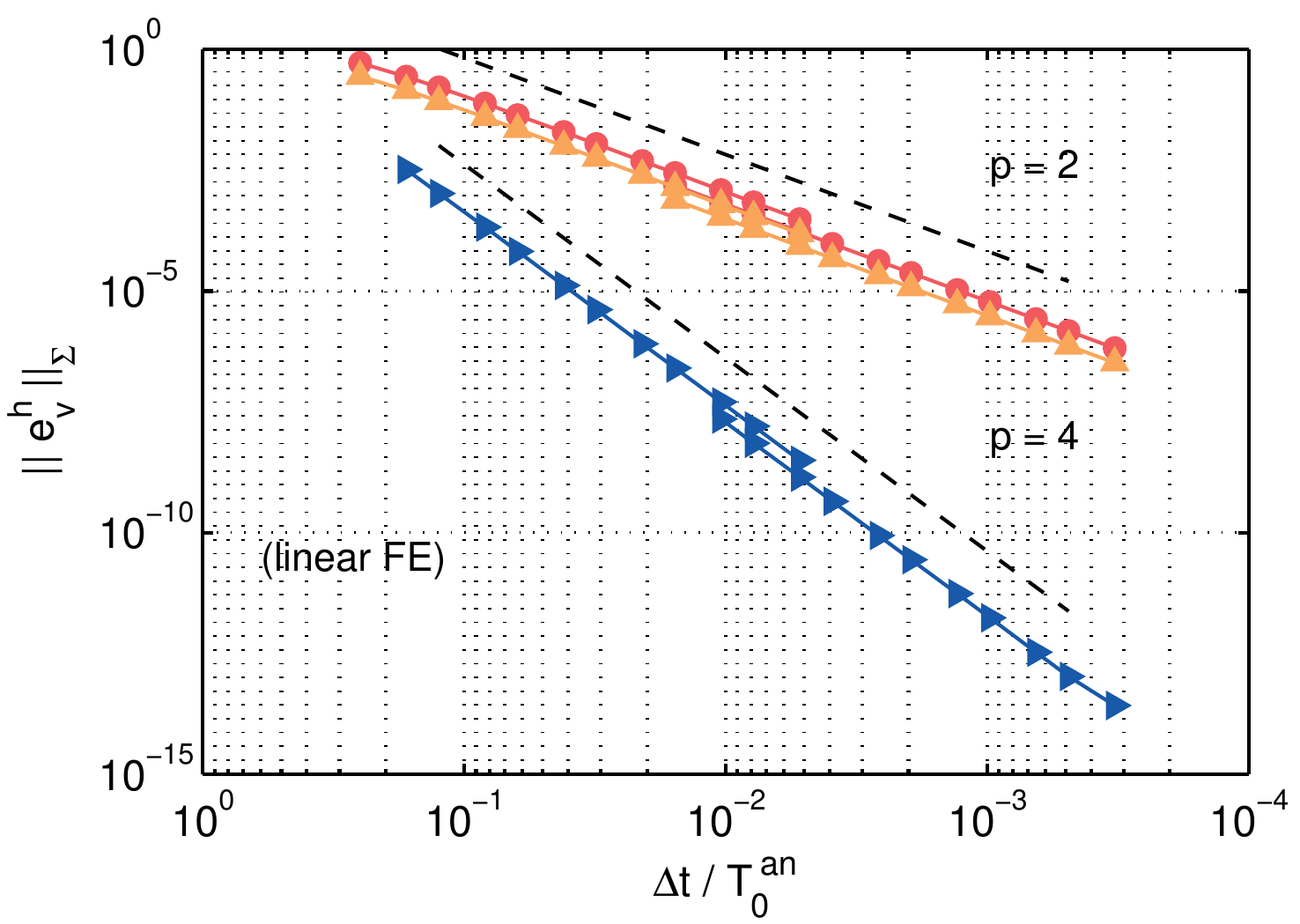}\label{f:L1bar:evh_F}}
	\caption{Linear elastic bar: Convergence of the temporal errors in
		the displacement and the velocity for linear Lagrange FE;
		$T = 1 \, T_0^\mathrm{an}$, $\Dt_0 = T_0^\mathrm{an} / 4$.}
	\label{f:L1bar:euvE}
\end{figure}

Note that if the bar is discretized with linear elements, the resulting system can be treated as a naturally discrete spring-mass system consisting of linear springs. For such a spring-mass system, the (temporally) analytical solution is
\begin{equation}
	u_\mathrm{an}^h ( X_I,t \big) = u_0 \cdot \cos\left({
	\frac{\pi \, X_I}{L}}\right) \cdot \cos \, (\omega\,t), \qquad
	v_\mathrm{an}^h ( X_I,t \big) = -u_0 \cdot \cos\left({\frac{\pi \, X_I}{L}
	}\right) \cdot \omega \, \sin \, (\omega\,t), \label{e:euvanh}
\end{equation}
$I = 1, \dots, \nno$; cf.~\Eq{e:bar:uvan}. The natural frequency, $\omega$, can be determined by analyzing the eigenmodes of the discrete system. \Fig{f:L1bar:euh_F} and \ref{f:L1bar:evh_F} show the maximum errors in the displacement and velocity arising from the temporal discretization. Here, the errors $e_u^h$ and $e_v^h$ are determined from \Eq{e:euvdiscr}, inserting the analytical solutions given by \Eq{e:euvanh}. As expected, the orders of convergence agree with those studied for a single degree of freedom (\Sect{s:conv:single}).

\Tab{t:H3bar:runtime} compares the computational cost of the three methods for those results shown in the right column of \Fig{f:CFL:euvE}, $C_\mathrm{CFL} = 0.5$. Here, we measure the computation time that is required to obtain an error in the displacement smaller either than $1\,\%$, $0.1\,\%$, or $0.01\,\%$; see \Fig{f:H3bar:eu_CFL}. Although the $\mrp2$-scheme must account for twice the number of unknowns within each time step, it takes --- due to its higher order of convergence --- less computation time than both the L1-integrator and Newmark's method. Note that for this linear example, the integrals in the discrete momenta can be evaluated analytically, i.e.~without numerical quadrature. For an example requiring quadrature also in time, we refer to the following section.
\begin{vchtable}[ht]
	\vchcaption{Linear elastic bar: Step size, $\Delta s$, and computation time,
		$T_\mathrm{ct}$, of the test cases from \Fig{f:H3bar:eu_CFL}
		($C_\mathrm{CFL} = 0.5$), for which the error in the displacement is
		smaller than $1\,\%$, $0.1\,\%$, and $0.01\,\%$; $T_\mathrm{ct}$ denotes
		the time for one oscillation.}
	\label{t:H3bar:runtime}
	\begin{tabular}{@{}llrlrlr@{}}
	\hline\noalign{\smallskip}
	& \multicolumn{2}{c}{Newmark}
	& \multicolumn{2}{c}{L1-integrator}
	& \multicolumn{2}{c@{}}{$\mrp2$-scheme} \\[0.3ex]
	& \ \,$\Delta s$ & $T_\mathrm{ct}$ [ms]\ \, 
			& \ \,$\Delta s$ & $T_\mathrm{ct}$ [ms]\ \,
			& \ \,$\Delta s$ & $T_\mathrm{ct}$ [ms] \\
		\hline\noalign{\smallskip}
			${||e_u||}_\Sigma < 1\,\%$	& \ \,1/32 & 79.6\ \,	
					& \ \,1/32 & 105.0\ \,	&\ \,1/16 & 44.9 \\[0.5ex]
			${||e_u||}_\Sigma < 0.1\,\%$	& \ \,1/128 & 815.8\ \,
					&\ \,1/64 & 416.3\ \,		&\ \,1/32 & 135.7 \\[0.5ex]
			${||e_u||}_\Sigma < 0.01\,\%$	& \ \,1/256 & 2,919.4\ \,
					&\ \,1/256 & 5,381.3\ \,		&\ \,1/64 & 574,3 \\
		\noalign{\smallskip}\hline
	\end{tabular}
\end{vchtable}


\subsection{Vibration of a nonlinear bar} \label{s:nonlvibr}

\begin{figure}[ht]
	\centering
	\subfigure[{Deformation for $t\,/\,T_0^\mathrm{an} \in [0,3]$}]{
		\includegraphics[width=0.48\textwidth]
			{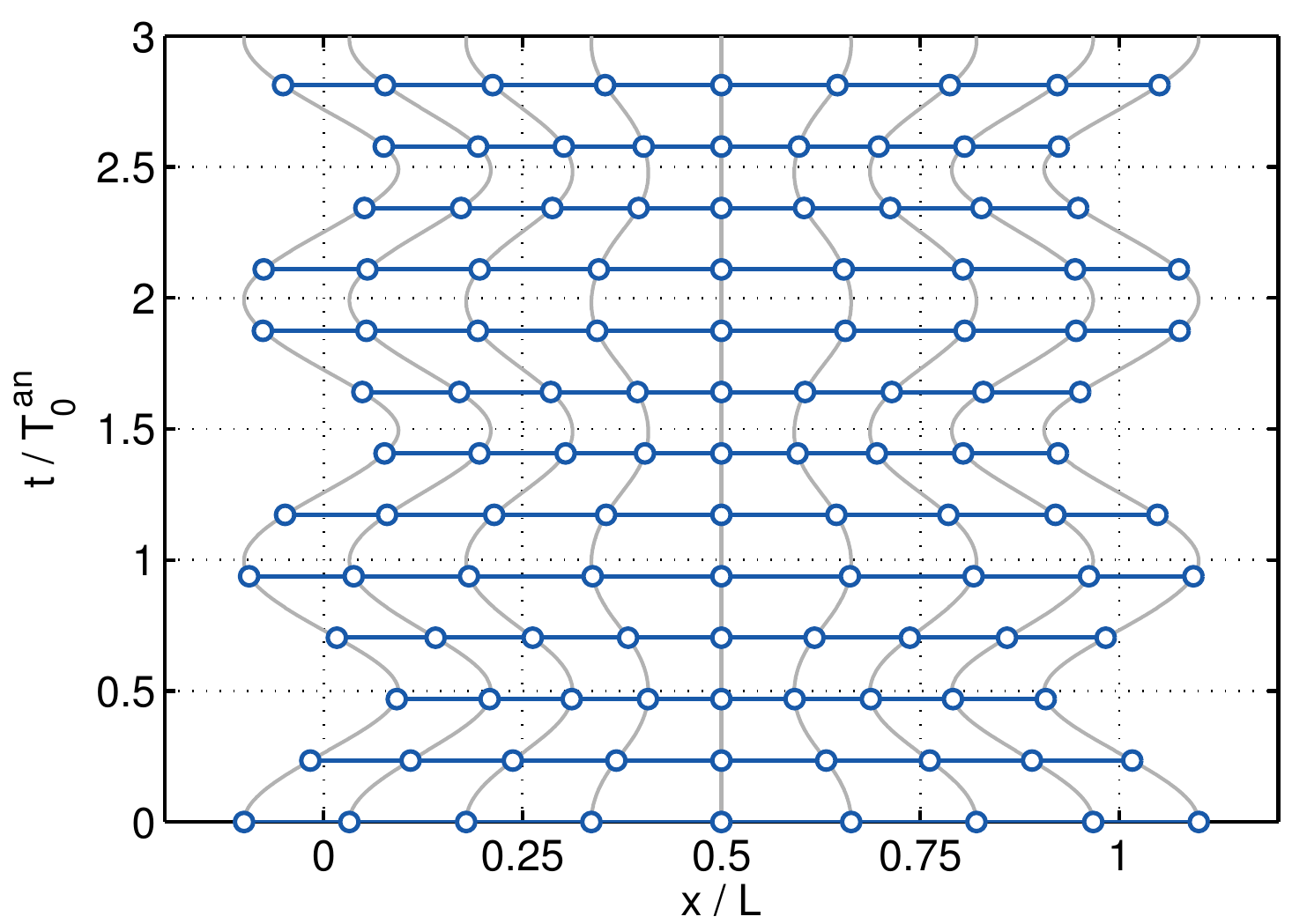}\label{f:nonlbar:deform1}}
	\subfigure[{Deformation for $t\,/\,T_0^\mathrm{an} \in [197,200]$}]{
		\includegraphics[width=0.48\textwidth]
			{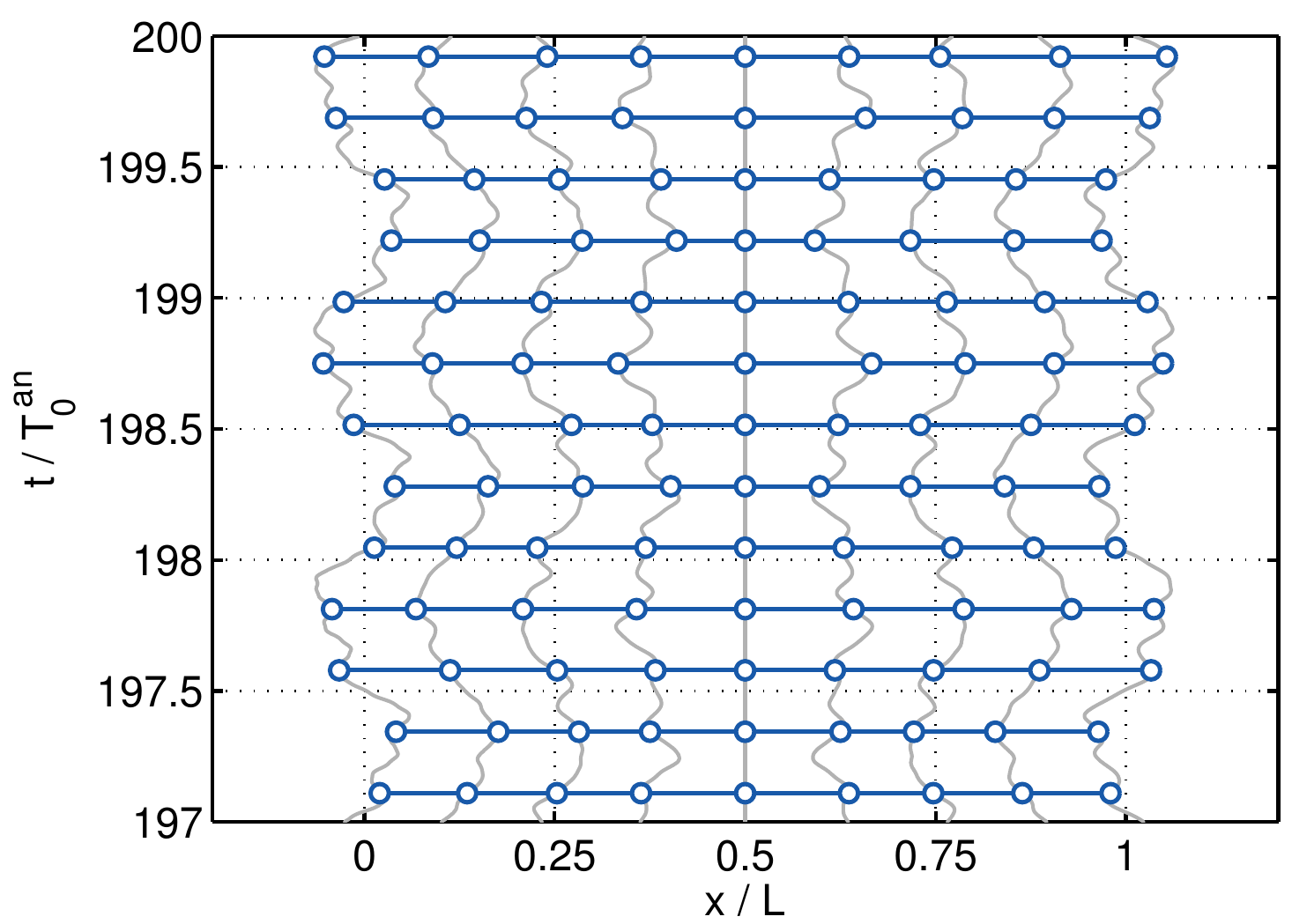}\label{f:nonlbar:deform2}}
	\caption{Nonlinear bar: Deformation of the bar (a)~at the beginning of
		the oscillation, and (b)~after a long time period, using the
		$\mrp2$-scheme and Hermite FE; $T = 200 \, T_0^\mathrm{an}$,
		$\Dt = T_0^\mathrm{an} / 512$, $L_e = L/8$.}
	\label{f:nonlbar:deform}
\end{figure}

The numerical examples discussed in the previous sections cover both naturally discrete and (spatially discretized) continuum systems. So far only linear problems (for which the internal forces depend on the displacement linearly) have been investigated. We therefore consider a nonlinear Neo-Hooke material behavior, which is described in \Appx{a:bar}. The bar is initially deformed by prescribing the same displacement and velocity as in \Eq{e:bar:uvan}. \Fig{f:nonlbar:deform} shows the deformation of the bar (a)~at the very beginning and (b)~after a long period of oscillations.
Here, a Hermite finite element mesh with eight elements is chosen. For a better comparison with the results from the previous section, the time step is normalized by the period length of the first eigenmode, $T_0^\mathrm{an}$. Since the mechanical response of the system differs from the linear case, however, the initially sinusoidal oscillations turn into a set of different interfering oscillations (\Fig{f:nonlbar:deform2}).

\begin{figure}[ht]
	\centering
	\includegraphics[width=\textwidth]{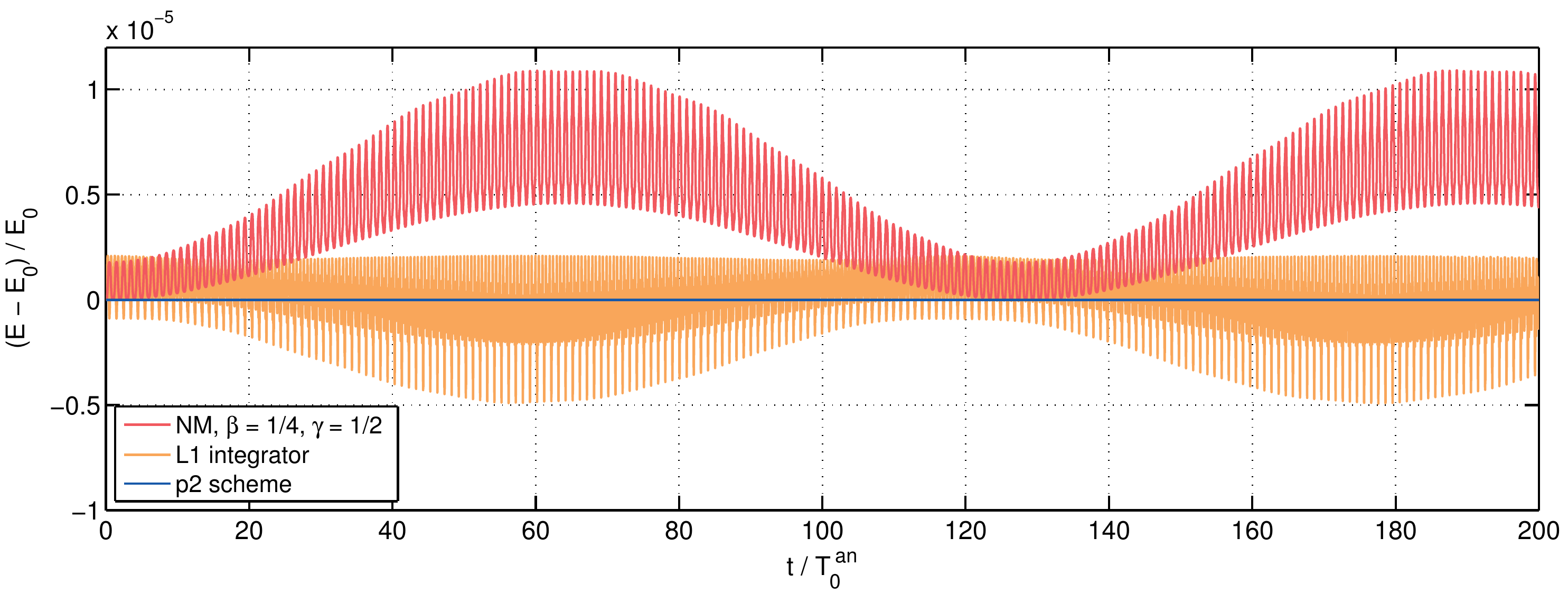}
	\caption{Nonlinear bar: Long-term energy behavior for approximately 200
		periods of oscillation using Hermite FE; the maximum relative error 
		of the $\mrp2$-scheme is $3.7 \cdot 10^{-10}$;
		$\Dt = T_0^\mathrm{an} / 256$, $L_e = L / 4$.}
	\label{f:nonlbar:E}
\end{figure}

\Fig{f:nonlbar:E} shows the long-term behavior for a nonlinear bar, comparing the system's total energy for the $\mrp2$-scheme, the Newmark algorithm, and the L1-integrator. Like for the harmonic oscillator from \Sect{s:longterm}, the energy oscillates while being qualitatively preserved. In comparison with the other methods the relative error of the $\mrp2$-scheme is smaller by five orders of magnitude.

\begin{figure}[ht]
	\centering
	\subfigure[{Error in the displacement (w.r.t.~reference solution)}]{
		\includegraphics[width=0.48\textwidth]
			{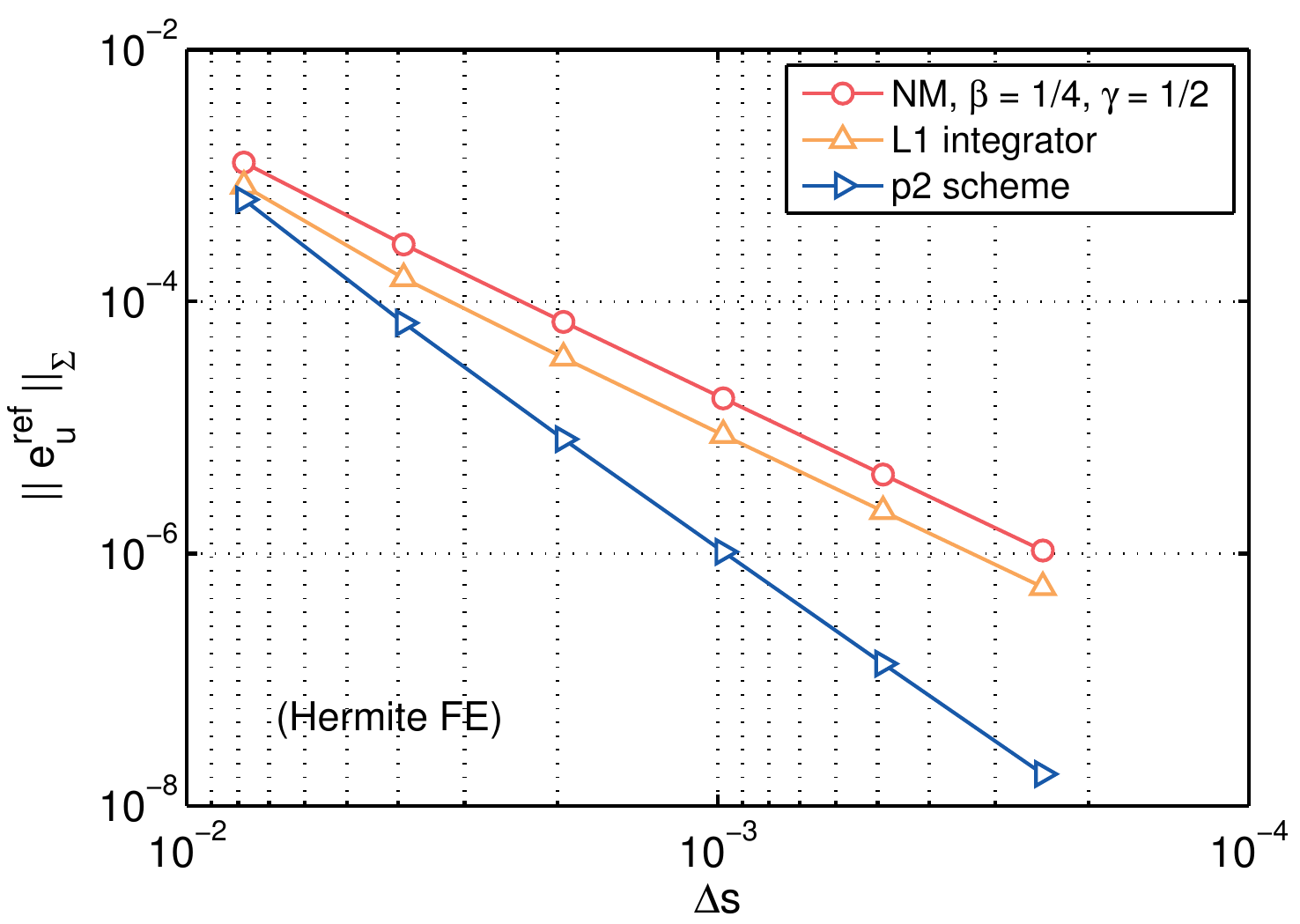}\label{f:nonlbar:eurel}}
	\subfigure[{Error in the velocity (w.r.t.~reference solution)}]{
		\includegraphics[width=0.48\textwidth]
			{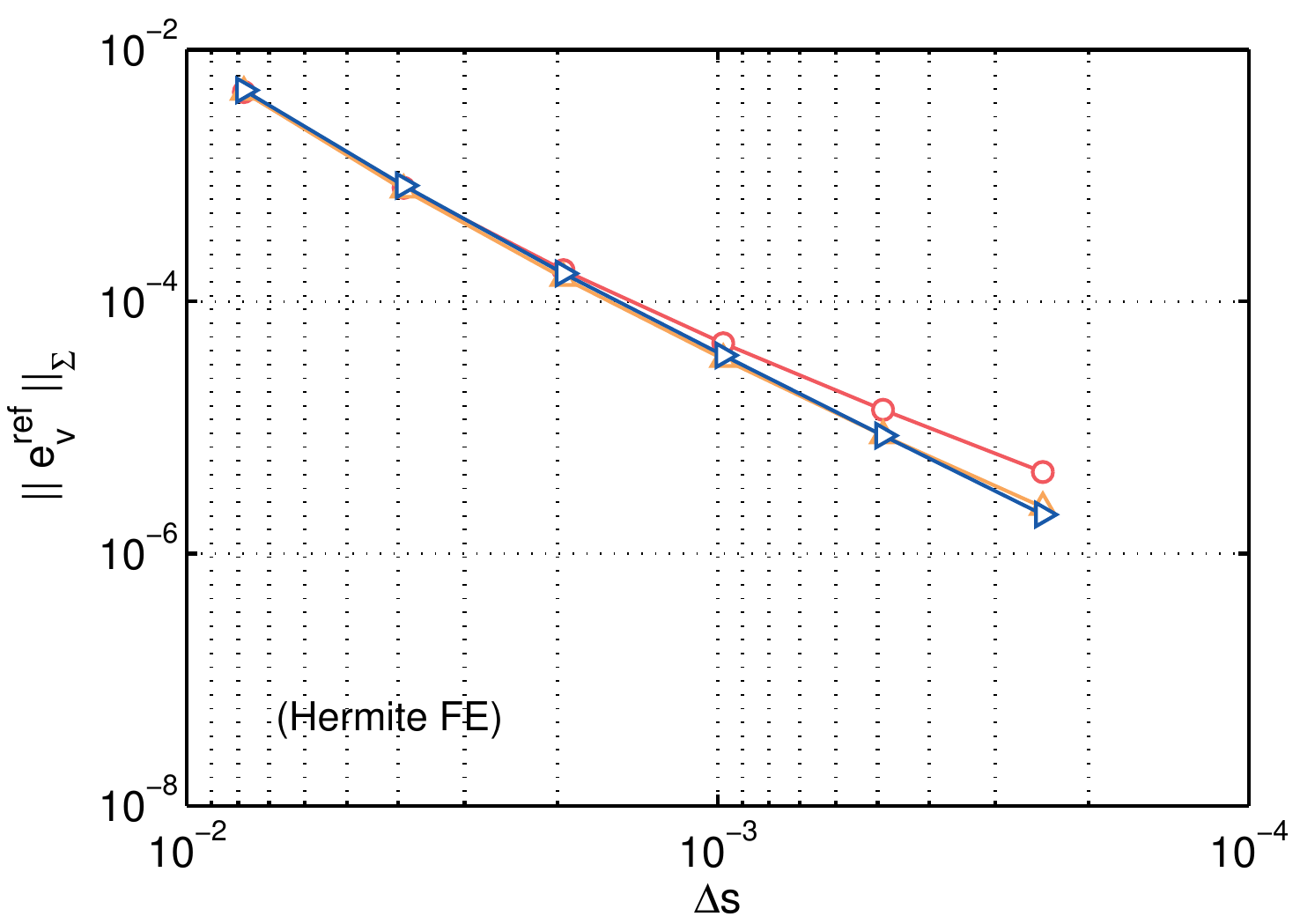}\label{f:nonlbar:evrel}}
		\subfigure[{Error in the total energy (w.r.t.~reference solution)}]{
		\includegraphics[width=0.48\textwidth]
			{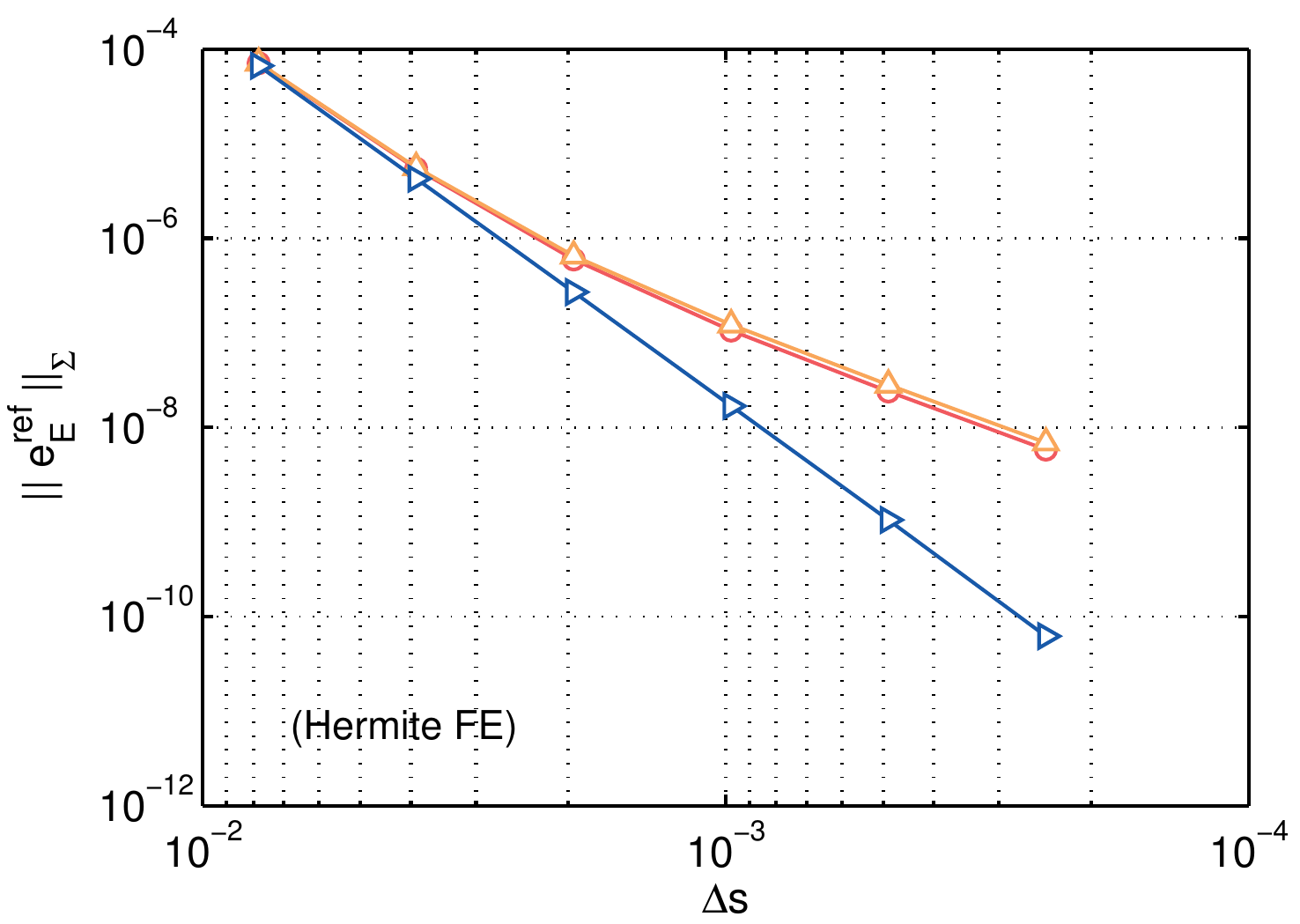}\label{f:nonlbar:eErel}}
	\subfigure[{Error in the total energy (w.r.t.~$E_0$)}]{
		\includegraphics[width=0.48\textwidth]
			{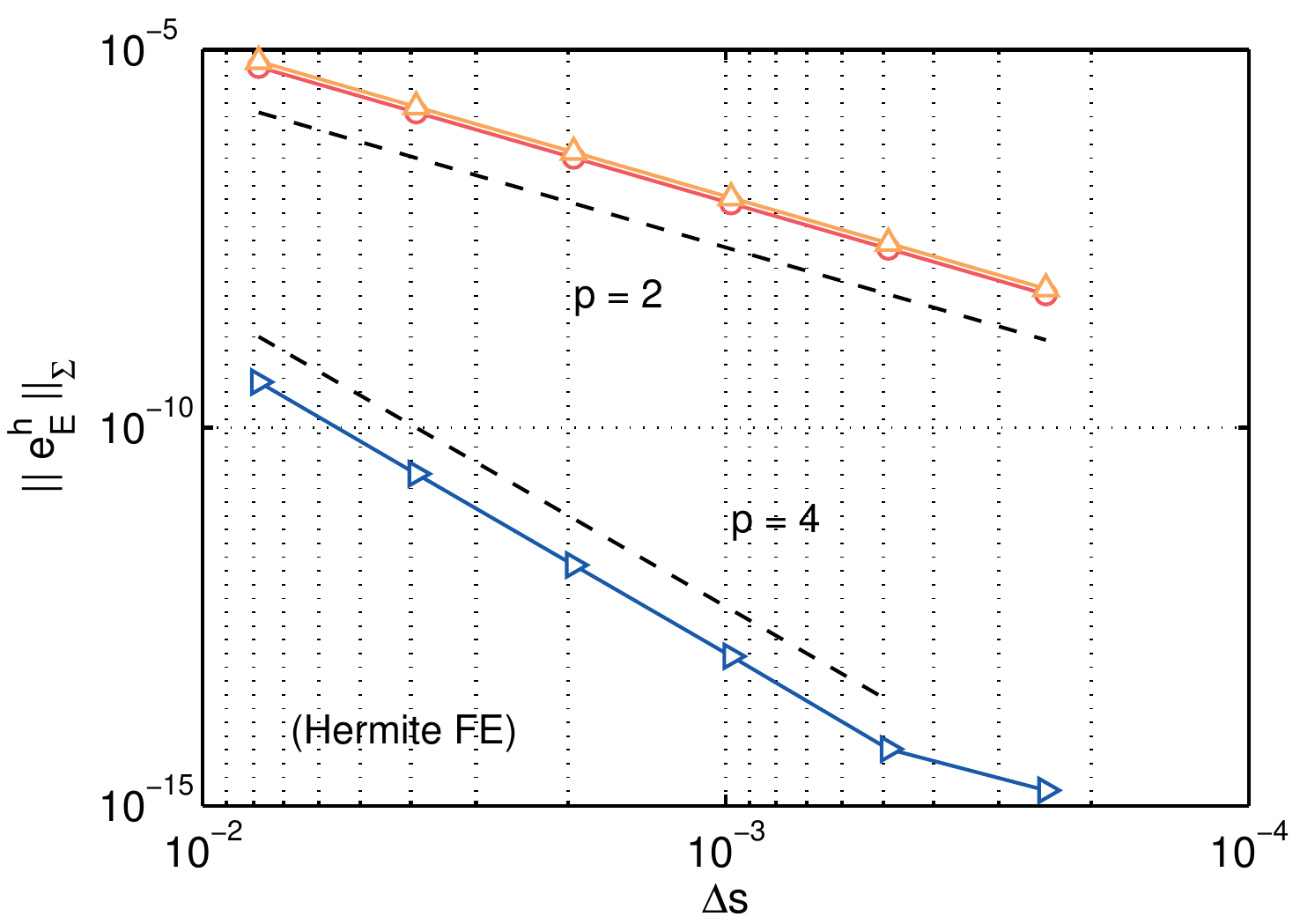}\label{f:nonlbar:eEh}}
	\caption{Nonlinear bar: Relative errors with respect to (a) -- (c) a very
		fine reference solution and (d)~to the initial energy of the spatially
		discrete system, refining both the mesh and time step; $\Delta s
		= \Dt/T_0^\mathrm{an} = \Delta L / (16\,L)$; shown is $C_\mathrm{CFL}
		= 0.125$ for Hermite FE and $T = 2 \, T_0^\mathrm{an}$; the reference
		solution is computed with the $\mrp2$-scheme using $L_e^\mathrm{ref}
		= L / 256$ and $\Dt^\mathrm{ref} = T_0^\mathrm{an} / 8192$.}
	\label{f:nonlbar:euvE}
\end{figure}

In addition, we investigate the convergence behavior for a nonlinear bar, considering Hermite finite elements in space. Since for this case the deformation cannot be computed analytically, we compare our results with a fine reference solution, using both a very fine FE mesh and a small time step. Like for the examples shown in the previous section, we refine both discretizations simultaneously. The results are shown in \Fig{f:nonlbar:eurel} -- \ref{f:nonlbar:eErel}. The accuracy of our scheme becomes most apparent for the displacement and the total energy; for these quantities we observe a significantly higher convergence. Note that caused by the nonlinear material law, both the $\mrp2$-scheme and the L1-integrator require numerical quadrature to evaluate the time integral over the internal forces; see also \Appx{a:implempq}. If we demand, however, a sufficient high accuracy w.r.t.~the fine solution, the computational cost may still be lower for the $\mrp2$-scheme than for the Newmark algorithm. An error in the displacement smaller than~$10^{-5}$ (\Fig{f:nonlbar:eurel}), for instance, requires a step size of $\Delta s = 1/2048$ for Newmark's method. In contrast, the $\mrp2$-scheme achieves this accuracy already for $\Delta s = 1/512$. Therefore, the measured computation time is significantly lower: $48.7\%$ of the time required for the Newmark algorithm.

We finally investigate how well the (spatially discrete) initial energy is preserved over time; see \Fig{f:nonlbar:eEh}. The error plotted here thus arises only from the temporal discretization. A comparison with \Fig{f:single:eE_max} shows that --- for both the linear oscillator and the nonlinear bar --- we achieve with the $\mrp2$-scheme the same high order of convergence: $p = 4$.


\section{Conclusion} \label{s:concl}

In this work we derive a class of $C^1$-continuous time integration methods for conservative elastodynamic problems. Using piecewise Hermite interpolation, we approximate the displacement of a deformable solid by $C^1$-continuous functions in time. The velocity of the body is thus $C^0$-continuous in the entire time domain. To explicitly account for any initial velocities in the system, we consider a generalization of Hamilton's principle, referred to as Hamilton's law of varying action. From the action integral of the continuous system we derive the spatially and temporally weak form of the equilibrium equation in elastodynamics. This expression is first discretized in space, using a standard (Galerkin) finite element method. Afterwards, the spatially discrete system is discretized in time by 1)~subdividing the temporal domain into a set of smaller time intervals, and 2)~approximating the displacement by cubic Hermite shape functions.

Generally, methods belonging to the class of variational integrators can be constructed by varying the action integral over the entire time domain. From this equation one can derive the discrete Euler-Lagrange equations in order to develop a subsequent time integration method. For cubic Hermite interpolation, however, the resulting variational integrator is unconditionally unstable; this issue has been discussed also by Riff and Baruch~\cite{riff84a}. Therefore, we first vary the action of each time interval individually, and then derive different one-step methods to solve for the new (unknown) displacement and velocity. This yields six time integration schemes with different properties. Technically, these methods are not variational in the sense that they are not derived from the virtual action of the total domain. Interestingly, the most favorable of these methods --- denoted as $\mrp2$-scheme --- offers similar advantageous properties like variational integrators: 1)~a qualitatively accurate long-term behavior, and 2) symplecticity for simple linear systems.

We first investigate the properties of our schemes by considering a harmonic oscillator. We then demonstrate numerically that the $\mrp2$-scheme shows both conditional stability and convergence of order four. Afterwards, we examine both linear and nonlinear problems accounting either for inherently continuous or spatially discrete systems. Our results show that for an appropriate spatial discretization, the $\mrp2$-scheme provides both reasonable computational effort and remarkably higher accuracy than variational integrators based on linear interpolation in time. This scheme also well-preserves the energy of the system for long-time integration. It remains to be seen whether our method is symplectic for arbitrary nonlinear systems. This and other properties (such as the convergence behavior and the preservation of momentum maps in the presence of symmetries) should be addressed in an analytical study in the future.

Since the scope of this work is the construction of time integration methods, we here focus on one-dimensional elastodynamic problems. An elaborate study of elastodynamic problems in either two or three dimensions may be the scope of future work. In addition, since the numerical results obtained within this study for spatial Hermite finite elements look very promising, it would be interesting to apply our time integrator in combination with other spatially $C^1$-continuous FE discretizations and --- in particular --- to isogeometric approaches~\cite{cottrell09,hughes05}. In contrast to higher-order Lagrangian finite elements, isogeometric elements can prevent high-frequency errors~\cite{cottrell06,hughes08} that require unnecessarily small time steps.


\begin{acknowledgement}
We are grateful to the German Research Foundation (DFG) supporting J.~C.~Mergel and R.~A.~Sauer through grants SA1822/5-1 and GSC111. We further thank Prof.~Yuri B.~Suris (Technical University of Berlin), Prof.~Sachin S.~Gautam (IIT~Guwahati), and Marcus G.~Schmidt (formerly RWTH Aachen University) for their helpful comments.
\end{acknowledgement}


\appendix

\section{Implementation of our schemes}

This section contains helpful details for the implementation of our time integration schemes.


\subsection{Shape functions} \label{a:shpfct}

In both the spatial and temporal discretizations, we consider either linear Lagrange or cubic Hermite shape functions. Usually, these functions are defined on master domains, denoted e.g.~by $\tau \in [-1,1]$ for the discretization in time. The mapping from the master to the temporal domain, $\tau \mapsto t$ ($t \in [t_n,t_{n+1}]$), is characterized by the Jacobian determinant $J_\tau := \partial t/\partial\tau$; this determinant is given by $J_\tau = (t_{n+1}-t_n)/2$ for both the L1-integrator and our six schemes.
\begin{itemize}
	\item Linear Lagrange shape functions (for the L1-integrator):
				\begin{equation}
					R_1(\tau) = \frac{1}{2} \, (1 - \tau), \qquad
					R_2(\tau) = \frac{1}{2} \, (1 + \tau). \label{e:shpL}
				\end{equation}
	\item Cubic Hermite shape functions (for our schemes):
				\begin{equation}
					\begin{array}{rclcrcl}
						R_{1}(\tau)
							&	=	&	\ds\frac{1}{4} \, (2 + \tau)(1 - \tau)^2,	& \quad	&
						R_{2}(\tau)
							&	=	&	\ds\frac{1}{4} \, (2 - \tau)(1 + \tau)^2
										= 1 - R_1(\tau),	\\[2ex]
						H_{1}(\tau)
							&	=	&	\ds\frac{J_\tau}{4} \, (\tau + 1) (1 - \tau)^2,
							&	\quad	&
						H_{2}(\tau)
							&	=	&	\ds\frac{J_\tau}{4} \, (\tau - 1) (1 + \tau)^2.
					\end{array} \label{e:shpH}
				\end{equation}
\end{itemize}


\subsection{Discrete Lagrangian and action integral} \label{a:Aht}

The discrete (pseudo-)momenta \eqref{e:defp} and \eqref{e:defq} are obtained by taking the derivatives of the spatially and temporally discrete action
\begin{equation}
	\sS^{ht} = \sum_{n=0}^{N-1} \sS^{ht}_{n+1}, \qquad
	\sS^{ht}_{n+1} = \intTnpo L^{ht} \big(\mx^t,\dot{\mx}^t\big)~\mrd t.
\end{equation}
Here, the discrete Lagrangian, $L^{ht}\big(\mx^t,\dot{\mx}^t\big) = K^{ht}\big(\dot{\mx}^t\big) - \Pi^{ht}\big(\mx^t\big)$, can be determined from
\begin{equation}
	K^{ht} \big(\dot{\mx}^t\big) = \sum_{e=1}^\nel \, \left[{ \frac{1}{2} \, 
	{\left( {\dot\mx_e^t} \right)}^\mrT \, \mm_e \, \dot\mx_e^t }\right], \qquad
	\Pi^{ht} \big(\mx^t\big) = \sum_{e=1}^\nel \, \left[{ \int_{\Omega_0^e}
	W \big(\mx^t_e\big)~\mrd V -\ {\left( {\mx_e^t} \right)}^\mrT
	\, \fext^e \big(\mx_e^t\big) }\right].
\end{equation}
In \Eq{e:pnmp} and \eqref{e:qnmp} we additionally use the relation $\mf\big(\mx^t\big) = \partial\Pi^{ht}\big(\mx^t\big) / \partial\mx^t$.

If required, the total energy of the system can be computed from $E^{ht} \big(\mx^t,\dot{\mx}^t\big) = K^{ht} \big(\dot{\mx}^t\big) + \Pi^{ht} \big(\mx^t\big)$.


\subsection{Implementation of the integrals} \label{a:implempq}

As pointed out in \Sect{s:implement}, some of the integrals in the discrete (pseudo-)momenta can be computed analytically. For this purpose, we split the internal force into a linear and into a nonlinear part,
\begin{equation}
	\fint\big(\mx^t\big) = \mk \, \muu^t + \mf_\mathrm{nlin} \big( \mx^t \big),
\end{equation}
where $\muu^t = \mx^t - \mX$, and where $\mk$ is a linear stiffness matrix. The case $\mf_\mathrm{nlin} = \mathbf{0}$ corresponds to linear elasticity, while $\mk = \mathbf{0}$ represents the fully nonlinear case. Further let
\begin{equation}
	\widetilde{\mf}\big(\mx^t\big) := \mf_\mathrm{nlin}\big(\mx^t\big) - \fext.
\end{equation}
Inserting \Eq{e:xt} and \eqref{e:vt} into \Eq{e:pnmp} and \eqref{e:qnmp} yields with $\Dt := t_{n+1}-t_n$
\begin{align}
	\pnm		&	=	\frac{1}{10\,\Dt} \, \mm \, \Big[ 12 \, \unpo - 12 \, \un
						- \Dt \, \vnpo - \Dt \, \vn \Big] + \intTnpo
						R_1(t) \, \widetilde{\mf} \big(\mx^t\big)~\mrd t \nonumber \\
				&	\qquad + \frac{\Dt}{420} \, \mk \, \Big[ 54 \, \unpo
						+ 156 \, \un - 13\,\Dt\,\vnpo + 22\,\Dt\,\vn \Big],
						\label{e:pnman} \\
	\pnpop	&	=	\frac{1}{10\,\Dt} \, \mm \, \Big[ 12 \, \unpo - 12 \, \un
						- \Dt \, \vnpo - \Dt \, \vn \Big] - \intTnpo
						R_2(t) \, \widetilde{\mf} \big(\mx^t\big)~\mrd t \nonumber \\
				& \qquad + \frac{\Dt}{420} \, \mk \, \Big[ - 156 \, \unpo
						- 54 \, \un + 22\,\Dt\,\vnpo - 13\,\Dt\,\vn \Big]
\end{align}
and
\begin{align}
	\qnm		&	=	\frac{1}{30} \, \mm \, \Big[ 3 \, \unpo - 3 \, \un
						+ \Dt \, \vnpo - 4 \, \Dt \, \vn \Big] + \intTnpo
						H_1(t) \, \widetilde{\mf} \big(\mx^t\big)~\mrd t \nonumber \\
				&	\qquad + \frac{\Dt^2}{420} \, \mk \, \Big[ 13 \, \unpo
						+ 22 \, \un - 3\,\Dt\,\vnpo + 4\,\Dt\,\vn \Big], \\
	\qnpop	&	=	\frac{1}{30} \, \mm \, \Big[ -3 \, \unpo + 3 \, \un
						+ 4 \, \Dt \, \vnpo - \Dt \, \vn \Big] - \intTnpo
						H_2(t) \, \widetilde{\mf} \big(\mx^t\big)~\mrd t \nonumber \\
				&	\qquad + \frac{\Dt^2}{420} \, \mk \, \Big[ 22 \, \unpo
						+ 13 \, \un - 4\,\Dt\,\vnpo + 3\,\Dt\,\vn \Big]. \label{e:qnpopan}
\end{align}
The derivatives with respect to $\xnpo$ and $\vnpo$ (required for linearization, see \Sect{s:implement}) are given by
\begin{align}
	\pa{\pnm}{\xnpo}		&	=	\frac{6}{5\,\Dt} \, \mm + \frac{9\,\Dt}{70} \, \mk
												+ \intTnpo R_1(t) \, R_2(t) \,
												\pa{\widetilde{\mf}\big(\mx^t\big)}{\mx^t}~\mrd t, \\
	\pa{\pnpop}{\xnpo}	&	=	\frac{6}{5\,\Dt} \, \mm - \frac{13\,\Dt}{35} \, \mk
												- \intTnpo R_2^2(t) \,
												\pa{\widetilde{\mf}\big(\mx^t\big)}{\mx^t}~\mrd t, \\
	\pa{\qnm}{\xnpo}		&	=	\frac{1}{10} \, \mm + \frac{13\,\Dt^2}{420} \, \mk
												+ \intTnpo H_1(t) \, R_2(t) \,
												\pa{\widetilde{\mf}\big(\mx^t\big)}{\mx^t}~\mrd t, \\
	\pa{\qnpop}{\xnpo}	&	=	-\frac{1}{10} \, \mm + \frac{11\,\Dt^2}{210} \, \mk
												- \intTnpo H_2(t) \, R_2(t) \,
												\pa{\widetilde{\mf}\big(\mx^t\big)}{\mx^t}~\mrd t
\end{align}
and
\begin{align}
	\pa{\pnm}{\vnpo}		&	=	- \frac{1}{10} \, \mm - \frac{13\,\Dt^2}{420} \, \mk
												+ \intTnpo R_1(t) \, H_2(t) \,
												\pa{\widetilde{\mf}\big(\mx^t\big)}{\mx^t}~\mrd t, \\
	\pa{\pnpop}{\vnpo}	&	=	- \frac{1}{10} \, \mm + \frac{11\,\Dt^2}{210} \, \mk
												- \intTnpo R_2(t) \, H_2(t) \,
												\pa{\widetilde{\mf}\big(\mx^t\big)}{\mx^t}~\mrd t, \\
	\pa{\qnm}{\vnpo}		&	=	\frac{\Dt}{30} \, \mm - \frac{\Dt^3}{140} \, \mk
												+ \intTnpo H_1(t) \, H_2(t) \,
												\pa{\widetilde{\mf}\big(\mx^t\big)}{\mx^t}~\mrd t, \\
	\pa{\qnpop}{\vnpo}	&	=	\frac{2\,\Dt}{15} \, \mm - \frac{\Dt^3}{105} \, \mk
												- \intTnpo H_2^2(t) \,
												\pa{\widetilde{\mf}\big(\mx^t\big)}{\mx^t}~\mrd t.
\end{align}


\subsection{Alternative representation of the p2-scheme} \label{a:ppreform}

By integrating $\pnm$ and $\pnpop$ from \Eq{e:pnmp} by parts, we obtain
\begin{equation}
	\pnm = \mm \, \vn + \intTnpo R_1(t) \, \Big[ \ddot{\mx}^t
	+ \mf\big(\mx^t\big) \Big]~\mrd t, \qquad \pnpop = \mm \, \vnpo
	- \intTnpo R_2(t) \, \Big[ \ddot{\mx}^t + \mf\big(\mx^t\big) \Big]~\mrd t;
\end{equation}
the $\mrp2$-scheme~\eqref{e:pp} thus fulfills
\begin{equation}
	\intTnpo R_1(t) \, \Big[ \ddot{\mx}^t + \mf\big(\mx^t\big) \Big]~\mrd t = 0,
	\qquad
	\intTnpo R_2(t) \, \Big[ \ddot{\mx}^t + \mf\big(\mx^t\big) \Big]~\mrd t = 0.
\end{equation}
Since further $R_1(t), R_2(t) \ge 0$ and $R_1(t) + R_2(t) = 1$ (see \Eq{e:shpH}), we arrive at \Eq{e:ppreform}.


\section{Linear variational integrator (L1)} \label{a:varintlin}

This section outlines the linear variational integrator (L1) that we use for comparison in \Sect{s:props} and~\ref{s:results}. Here, the deformation and velocity are approximated by $\mx(t) \approx \mx^t(t)$ and $\dot{\mx}(t) \approx \dot{\mx}^t(t)$,
\begin{equation}
	\mx^t(t) = R_1(t) \, \xn + R_2(t) \, \xnpo, \qquad \dot{\mx}^t(t)
	= \dot{R}_1(t) \, \xn + \dot{R}_2(t) \, \xnpo, \qquad t \in [t_n,t_{n+1}],
	\label{e:xtL}
\end{equation}
where $R_1$ and $R_2$ are defined by \Eq{e:shpL}. In analogy, we introduce
\begin{equation}
	\delta\mx^t(t) = R_1(t) \, \delta\xn + R_2(t) \, \delta\xnpo, \qquad
	\delta\dot{\mx}^t(t) = \dot{R}_1(t) \, \delta\xn
	+ \dot{R}_2(t) \, \delta\xnpo.
\end{equation}
With this interpolation the discretized action becomes $\A^h \approx \A^{ht}$,
\begin{equation}
	\A^{ht} = \sum_{n=0}^{N-1} \A_{n+1}^{ht} \left({ \xn, \xnpo }\right).
\end{equation}
From the discretization of the virtual action then follows that
\begin{equation}
	\sum_{n=0}^{N-1} \left[{ \pa{\A_{n+1}^{ht}}{\xn} \, \delta\xn
	+ \pa{\A_{n+1}^{ht}}{\xnpo} \, \delta\xnpo }\right]
	= \delta{\mx}_N \cdot \mm
	\, \dot{\mx}^t(T) - \delta{\mx}_0 \cdot \mm \, \dot{\mx}^t(0).
\end{equation}
Since the $n^\mathrm{th}$ time step belongs to two time intervals, $[t_{n-1},t_n]$ and $[t_n,t_{n+1}]$, we obtain
\begin{align}
	\mm \, \dot{\mx}^t(0) + \pa{\A_1^{ht} \left({\hat{\mx}_0,\hat{\mx}_1}
	\right)} {\hat{\mx}_0}	&	=	\bzero, \\
	\pa{\A_n^{ht} \left({ \xnmo, \xn }\right)}{\xn}
	+ \pa{\A_{n+1}^{ht} \left({ \xn, \xnpo }\right)}{\xn}
		&	= \bzero, \qquad n = 1,\dots,N\!-\!1.
\end{align}
The integrals within the partial derivatives of $\A_{n+1}^{ht}$ are also evaluated by Gaussian quadrature.


\section{Amplification matrices for the simple harmonic oscillator} \label{a:ampmat}

The amplification matrices of our six schemes (\ref{e:pp}) -- (\ref{e:pmqp}) (required for \Eq{e:ampmat}) are given by
\begin{align*}
	\mA_{\mrp2}
		&	=	\frac{1}{8\gamma^4 \!+\! 132\gamma^2 \!+\! 2016}
				\begin{bmatrix}
						\big(  26\gamma^4 \!-\!  876\gamma^2 \!+\! 2016 \big)
					&	\big( 204\gamma^3 \!-\! 2016\gamma \big) \\
						\big(   3\gamma^5 \!-\!  204\gamma^3 \!+\! 2016\gamma \big)
					&	\big(  26\gamma^4 \!-\!  876\gamma^2 \!+\! 2016 \big) \\
				\end{bmatrix}, \\[1ex]
	\mA_{\mrq2}
		&	=	\frac{1}{2\gamma^4 \!+\! 18\gamma^2 \!+\! 420}
				\begin{bmatrix}
						\big(  7\gamma^4 \!-\! 192\gamma^2 \!+\! 420 \big)
					&	\big( 45\gamma^3 \!-\! 420\gamma \big) \\
						\big(   \gamma^5 \!-\!  52\gamma^3 \!+\! 420\gamma \big)
					&	\big(  7\gamma^4 \!-\!  192\gamma^2 \!+\! 420 \big) \\
				\end{bmatrix}, \\[1ex]
	\mA_{\mrp^+\mrq^-}
		&	=	\frac{1}{26\gamma^4 \!+\! 198\gamma^2 \!+\! 3780}
				\begin{bmatrix}
						\big(  65\gamma^4 \!-\! 1692\gamma^2 \!+\! 3780 \big)
					&	\big( 390\gamma^3 \!-\! 3780\gamma \big) \\
						\big(   7\gamma^5 \!-\!  432\gamma^3 \!+\! 3780\gamma \big)
					&	\big(  46\gamma^4 \!-\! 1692\gamma^2 \!+\! 3780 \big) \\
				\end{bmatrix}, \\[1ex]
	\mA_{\mrp^+\mrq^+}
		&	=	\frac{1}{10\gamma^4 \!+\! 24\gamma^2 \!+\! 630}
				\begin{bmatrix}
						\big( 13\gamma^4 \!-\! 291\gamma^2 \!+\! 630 \big)
					&	\big( 60\gamma^3 \!-\! 630\gamma \big) \\
						\big(   \gamma^5 \!-\!  81\gamma^3 \!+\! 630\gamma \big)
					&	\big(  5\gamma^4 \!-\! 291\gamma^2 \!+\! 630 \big) \\
				\end{bmatrix}, \\[1ex]
	\mA_{\mrp^-\mrq^-}
		&	=	\frac{1}{\gamma^4 \!+\! 48\gamma^2 \!+\! 1260}
				\begin{bmatrix}
						\big(  10\gamma^4 \!-\!  582\gamma^2 \!+\! 1260 \big)
					&	\big( 120\gamma^3 \!-\! 1260\gamma \big) \\
						\big(   2\gamma^5 \!-\!  162\gamma^3 \!+\! 1260\gamma \big)
					&	\big(  26\gamma^4 \!-\!  582\gamma^2 \!+\! 1260 \big) \\
				\end{bmatrix}, \\[1ex]
	\mA_{\mrp^-\mrq^+}
		&	=	\frac{1}{10\gamma^4 \!+\! 198\gamma^2 \!+\! 3780}
				\begin{bmatrix}
						\big(  46\gamma^4 \!-\! 1692\gamma^2 \!+\! 3780 \big)
					&	\big( 390\gamma^3 \!-\! 3780\gamma \big) \\
						\big(   7\gamma^5 \!-\!  432\gamma^3 \!+\! 3780\gamma \big)
					&	\big(  65\gamma^4 \!-\! 1692\gamma^2 \!+\! 3780 \big) \\
				\end{bmatrix}.
\end{align*}


\section{Equations for the 1D bar} \label{a:bar}

In \Sect{s:results} we investigate the axial deformation of a thin bar that is characterized by its length, $L$, density, $\rho_0$, cross section area, $A$, and Young's modulus,~$E$. We study both linear and nonlinear problems by considering different material behavior for the bar. Let $X$, $x$, and $v$ denote the 1D equivalents of the terms written in bold font. We further assume that $\bar B = \bar T = 0$. For the 1D case the virtual kinetic energy, given by \Eq{e:dK}, reduces to
\begin{equation}
	\delta K(v) = \rho_0 \, A \int_0^L \delta v \cdot v~\mrd X. \label{e:dKbar}
\end{equation}
For a linear elastic bar (discussed in \Sect{s:freevibr} and \ref{s:freevibr:conv}), the virtual work \eqref{e:dPi} is given by
\begin{equation}
	\delta\Pi(x) = A \int_0^L \delta\varepsilon \cdot E \,
	\varepsilon~\mrd X, \qquad \varepsilon := \lambda - 1, \label{e:siglin}
\end{equation}
where $\lambda := \partial x / \partial X$. In \Sect{s:nonlvibr} we consider a Neo-Hooke material with Poisson's ratio $\nu = 0$; see e.g.~Ref.~\cite{zienkiewicz05}:
\begin{equation}
	\delta\Pi(x) = A \int_0^L \delta\lambda \cdot P~\mrd X, \qquad
	P := \frac{E}{2} \left({\lambda - \lambda^{-1}}\right). \label{e:sigNH}
\end{equation}
For the spatial discretization of \Eq{e:dKbar} -- \eqref{e:sigNH} by means of 1D finite elements, we refer e.g.~to Ref.~\cite{wriggers08}.


\bibliographystyle{zamm-title}
\bibliography{bibliography}

\providecommand{\WileyBibTextsc}{}
\let\textsc\WileyBibTextsc
\providecommand{\othercit}{}
\providecommand{\jr}[1]{#1}
\providecommand{\etal}{~et~al.}


\begin{thebibliography}{[10]}

\bibitem{argyris69}
 \textsc{J.\,H. Argyris} and  \textsc{D.\,W. Scharpf},
Finite elements in time and space,
 \jr{Nucl. Eng. Des.} \textbf{\textbf{10}}(4), 456--464 (1969).


\bibitem{bailey75a}
 \textsc{C.\,D. Bailey},
Application of {H}amilton's law of varying action,
 \jr{AIAA J.} \textbf{\textbf{13}}(9), 1154--1157 (1975).


\bibitem{bailey75b}
 \textsc{C.\,D. Bailey},
A new look at {H}amilton's principle,
 \jr{Found. Phys.} \textbf{\textbf{5}}(3), 433--451 (1975).


\bibitem{bailey76a}
 \textsc{C.\,D. Bailey},
Hamilton, {R}itz, and elastodynamics,
 \jr{J. Appl. Mech.} \textbf{\textbf{43}}(4), 684--688 (1976).


\bibitem{bailey76b}
 \textsc{C.\,D. Bailey},
The method of {R}itz applied to the equation of {H}amilton,
 \jr{Comput. Methods Appl. Mech. Eng.} \textbf{\textbf{7}}, 235--247 (1976).


\bibitem{baruch82}
 \textsc{M.~Baruch} and  \textsc{R.~Riff},
{H}amilton's principle, {H}amilton's law -- {$6^n$} correct formulations,
 \jr{AIAA J.} \textbf{\textbf{20}}(5), 687--692 (1982).


\bibitem{bauchau99}
 \textsc{O.\,A. Bauchau} and  \textsc{T.~Joo},
Computational schemes for non-linear elasto-dynamics,
 \jr{Int. J. Numer. Methods Eng.} \textbf{\textbf{45}}, 693--719 (1999).


\othercit
\bibitem{belytschko00}
 \textsc{T.~Belytschko},  \textsc{W.\,K. Liu},  and  \textsc{B.~Moran},
Nonlinear Finite Elements for Continua and Structures (Wiley, Chichester,
  2000).


\bibitem{betsch16}
 \textsc{P.~Betsch} and  \textsc{A.~Janz},
An energy-momentum consistent method for transient simulations with mixed
  finite elements developed in the framework of geometrically exact shells,
 \jr{Int. J. Numer. Methods Eng.} (2016),
DOI: 10.1002/nme.5217.


\bibitem{betsch01}
 \textsc{P.~Betsch} and  \textsc{P.~Steinmann},
Conservation properties of a time {FE} method -- part {II}: {T}ime-stepping
  schemes for non-linear elastodynamics,
 \jr{Int. J. Numer. Methods Eng.} \textbf{\textbf{50}}, 1931--1955 (2001).


\bibitem{borri85}
 \textsc{M.~Borri},  \textsc{G.\,L. Ghiringhelli},  \textsc{M.~Lanz},
  \textsc{P.~Mantegazza},  and  \textsc{T.~Merlini},
Dynamic response of mechanical systems by a weak {H}amiltonian formulation,
 \jr{Comput. Struct.} \textbf{\textbf{20}}(1--3), 495--508 (1985).


\bibitem{bourabee08}
 \textsc{N.~Bou-Rabee} and  \textsc{H.~Owhadi},
Stochastic variational integrators,
 \jr{IMA J. Numer. Anal.} \textbf{\textbf{29}}, 421--443 (2008).


\othercit
\bibitem{chadwick99}
 \textsc{P.~Chadwick},
Continuum Mechanics: Concise Theory and Problems (Dover Publications, 1999).


\bibitem{chung93}
 \textsc{J.~Chung} and  \textsc{G.\,M. Hulbert},
A time integration algorithm for structural dynamics with improved numerical
  dissipation: {T}he generalized-{$\alpha$} method,
 \jr{J. Appl. Mech.} \textbf{\textbf{60}}, 371--375 (1993).


\bibitem{cortes01}
 \textsc{J.~Cort{\'e}s} and  \textsc{S.~Mart{\'\i}nez},
Nonholonomic integrators,
 \jr{Nonlinearity} \textbf{\textbf{14}}(5), 1365--1392 (2001).


\othercit
\bibitem{cottrell09}
 \textsc{J.\,A. Cottrell},  \textsc{T.\,J. Hughes},  and  \textsc{Y.~Bazilevs},
Isogeometric Analysis: Toward Integration of CAD and FEA (Wiley, 2009).


\bibitem{cottrell06}
 \textsc{J.\,A. Cottrell},  \textsc{A.~Reali},  \textsc{Y.~Bazilevs},  and
  \textsc{T.\,J.\,R. Hughes},
Isogeometric analysis of structural vibrations,
 \jr{Comput. Methods Appl. Mech. Eng.} \textbf{\textbf{195}}, 5257--5296
  (2006).


\bibitem{demoures15}
 \textsc{F.~Demoures},  \textsc{F.~Gay-Balmaz},  \textsc{S.~Leyendecker},
  \textsc{S.~Ober-Bl{\"o}baum},  \textsc{T.\,S. Ratiu},  and
  \textsc{Y.~Weinand},
Discrete variational {L}ie group formulation of geometrically exact beam
  dynamics,
 \jr{Numer. Math.} \textbf{\textbf{130}}(1), 73--123 (2015).


\bibitem{fetecau03}
 \textsc{R.\,C. Fetecau},  \textsc{J.\,E. Marsden},  \textsc{M.~Ortiz},  and
  \textsc{M.~West},
Nonsmooth {L}agrangian mechanics and variational collision integrators,
 \jr{SIAM J. Appl. Dyn. Syst.} \textbf{\textbf{2}}(3), 381--416 (2003).


\bibitem{fried69}
 \textsc{I.~Fried},
Finite-element analysis of time-dependent phenomena,
 \jr{AIAA J.} \textbf{\textbf{7}}(6), 1170--1173 (1969).


\bibitem{fung96}
 \textsc{T.\,C. Fung},
Unconditionally stable higher-order accurate {H}ermitian time finite elements,
 \jr{Int. J. Numer. Methods Eng.} \textbf{\textbf{39}}, 3475--3495 (1996).


\bibitem{gautam13}
 \textsc{S.\,S. Gautam} and  \textsc{R.\,A. Sauer},
An energy-momentum-conserving temporal discretization scheme for adhesive
  contact problems,
 \jr{Int. J. Numer. Methods Eng.} \textbf{\textbf{93}}(10), 1057--1081 (2013).


\bibitem{geradin74}
 \textsc{M.~G{\'e}radin},
On the variational method in the direct integration of the transient structural
  response,
 \jr{J. Sound Vib.} \textbf{\textbf{34}}(4), 479--487 (1974).


\bibitem{gonzalez00}
 \textsc{O.~Gonzalez},
Exact energy and momentum conserving algorithms for general models in nonlinear
  elasticity,
 \jr{Comput. Methods Appl. Mech. Eng.} \textbf{\textbf{190}}, 1763--1783
  (2000).


\bibitem{gross05}
 \textsc{M.~Gro{\ss}},  \textsc{P.~Betsch},  and  \textsc{P.~Steinmann},
Conservation properties of a time {FE} method. {P}art {IV}: {H}igher order
  energy and momentum conserving schemes,
 \jr{Int. J. Numer. Methods Eng.} \textbf{\textbf{63}}, 1849--1897 (2005).


\othercit
\bibitem{hairer06}
 \textsc{E.~Hairer},  \textsc{C.~Lubich},  and  \textsc{G.~Wanner},
Geometric Numerical Integration, $2^{\text{nd}}$ edition (Springer, Berlin
  Heidelberg, 2006).


\bibitem{hamilton34}
 \textsc{W.\,R. Hamilton},
On a general method in dynamics; by which the study of the motions of all free
  systems of attracting or repelling points is reduced to the search and
  differentiation of one central relation, or characteristic function,
 \jr{Phil. Trans. R. Soc. London} \textbf{\textbf{124}}, 247--308 (1834).


\bibitem{hamilton35}
 \textsc{W.\,R. Hamilton},
Second essay on a general method in dynamics,
 \jr{Phil. Trans. R. Soc. London} \textbf{\textbf{125}}, 95--144 (1835).


\bibitem{hesch09}
 \textsc{C.~Hesch} and  \textsc{P.~Betsch},
A mortar method for energy-momentum conserving schemes in frictionless dynamic
  contact problems,
 \jr{Int. J. Numer. Methods Eng.} \textbf{\textbf{77}}, 1468--1500 (2009).


\bibitem{hesch10}
 \textsc{C.~Hesch} and  \textsc{P.~Betsch},
Transient three-dimensional domain decomposition problems: {F}rame-indifferent
  mortar constraints and conserving integration,
 \jr{Int. J. Numer. Methods Eng.} \textbf{\textbf{82}}, 329--358 (2010).


\bibitem{hilber77}
 \textsc{H.\,M. Hilber},  \textsc{T.\,J.\,R. Hughes},  and  \textsc{R.\,L.
  Taylor},
Improved numerical dissipation for time integration algorithms in structural
  dynamics,
 \jr{Earthq. Eng. Struct. Dyn.} \textbf{\textbf{5}}, 283--292 (1977).


\bibitem{hodges91}
 \textsc{D.\,H. Hodges} and  \textsc{R.\,R. Bless},
Weak {H}amiltonian finite element method for optimal control problems,
 \jr{J. Guidance Control Dyn.} \textbf{\textbf{14}}(1), 148--156 (1991).


\othercit
\bibitem{holzapfel00}
 \textsc{G.\,A. Holzapfel},
Nonlinear Solid Mechanics: A Continuum Approach for Engineering (Wiley,
  Chichester, 2000).


\bibitem{hughes05}
 \textsc{T.\,J.\,R. Hughes}, ,  \textsc{J.\,A. Cottrell},  and
  \textsc{Y.~Bazilevs},
Isogeometric analysis: {CAD}, finite elements, {NURBS}, exact geometry and mesh
  refinement,
 \jr{Comput. Methods Appl. Mech. Eng.} \textbf{\textbf{194}}, 4135--4195
  (2005).


\bibitem{hughes76}
 \textsc{T.\,J.\,R. Hughes},
Stability, convergence and growth and decay of energy of the average
  acceleration method in nonlinear structural dynamics,
 \jr{Comput. Struct.} \textbf{\textbf{6}}(4--5), 313--324 (1976).


\bibitem{hughes88}
 \textsc{T.\,J.\,R. Hughes} and  \textsc{G.\,M. Hulbert},
Space-time finite element methods for elastodynamics: {F}ormulations and error
  estimates,
 \jr{Comput. Methods Appl. Mech. Eng.} \textbf{\textbf{66}}, 339--363 (1988).


\bibitem{hughes08}
 \textsc{T.\,J.\,R. Hughes},  \textsc{A.~Reali},  and  \textsc{G.~Sangalli},
Duality and unified analysis of discrete approximations in structural dynamics
  and wave propagation: {C}omparison of {$p$}-method finite elements with
  {$k$}-method {NURBS},
 \jr{Comput. Methods Appl. Mech. Eng.} \textbf{\textbf{197}}, 4104--4124
  (2008).


\bibitem{hulbert92}
 \textsc{G.\,M. Hulbert},
Time finite element methods for structural dynamics,
 \jr{Int. J. Numer. Methods Eng.} \textbf{\textbf{33}}, 307--331 (1992).


\bibitem{hulbert90}
 \textsc{G.\,M. Hulbert} and  \textsc{T.\,J.\,R. Hughes},
Space-time finite element methods for second-order hyperbolic equations,
 \jr{Comput. Methods Appl. Mech. Eng.} \textbf{\textbf{84}}, 327--348 (1990).


\bibitem{johnson14}
 \textsc{G.~Johnson},  \textsc{S.~Leyendecker},  and  \textsc{M.~Ortiz},
Discontinuous variational time integrators for complex multibody collisions,
 \jr{Int. J. Numer. Methods Eng.} \textbf{\textbf{100}}, 871--913 (2014).


\bibitem{kane00}
 \textsc{C.~Kane},  \textsc{J.\,E. Marsden},  \textsc{M.~Ortiz},  and
  \textsc{M.~West},
Variational integrators and the {N}ewmark algorithm for conservative and
  dissipative mechanical systems,
 \jr{Int. J. Numer. Methods Eng.} \textbf{\textbf{49}}, 1295--1325 (2000).


\bibitem{kobilarov10}
 \textsc{M.~Kobilarov},  \textsc{J.\,E. Marsden},  and  \textsc{G.\,S.
  Sukhatme},
Geometric discretization of nonholonomic systems with symmetries,
 \jr{Discrete Contin. Dyn. Syst. Ser. S} \textbf{\textbf{3}}(1), 61--84 (2010).


\bibitem{krenk06}
 \textsc{S.~Krenk},
Energy conservation in {N}ewmark based time integration algorithms,
 \jr{Comput. Methods Appl. Mech. Eng.} \textbf{\textbf{195}}, 6110--6124
  (2006).


\bibitem{krenk14}
 \textsc{S.~Krenk},
Global format for energy-momentum based time integration in nonlinear dynamics,
 \jr{Int. J. Numer. Methods Eng.} \textbf{\textbf{100}}, 458--476 (2014).


\bibitem{kuhl99a}
 \textsc{D.~Kuhl} and  \textsc{M.\,A. Crisfield},
Energy-conserving and decaying algorithms in non-linear structural dynamics,
 \jr{Int. J. Numer. Methods Eng.} \textbf{\textbf{45}}, 569--599 (1999).


\bibitem{kuhl99b}
 \textsc{D.~Kuhl} and  \textsc{E.~Ramm},
Generalized energy-momentum method for non-linear adaptive shell dynamics,
 \jr{Comput. Methods Appl. Mech. Eng.} \textbf{\textbf{178}}, 343--366 (1999).


\bibitem{kuhl00}
 \textsc{D.~Kuhl} and  \textsc{E.~Ramm},
Time integration in the context of energy control and locking free finite
  elements,
 \jr{Arch. Comput. Methods Eng.} \textbf{\textbf{7}}(3), 299--332 (2000).


\othercit
\bibitem{lanczos70}
 \textsc{C.~Lanczos},
The Variational Principles of Mechanics, $4^{\text{th}}$ edition (Dover
  Publications, New York, 1970).


\othercit
\bibitem{leimkuhler05}
 \textsc{B.~Leimkuhler} and  \textsc{S.~Reich},
Simulating Hamiltonian Dynamics (Cambridge University Press, Cambridge, 2005).


\bibitem{leok12}
 \textsc{M.~Leok} and  \textsc{T.~Shingel},
Prolongation-collocation variational integrators,
 \jr{IMA J. Numer. Anal.} \textbf{\textbf{32}}(3), 1194--1216 (2012).


\bibitem{lew03}
 \textsc{A.~Lew},  \textsc{J.\,E. Marsden},  \textsc{M.~Ortiz},  and
  \textsc{M.~West},
Asynchronous variational integrators,
 \jr{Arch. Rational Mech. Anal.} \textbf{\textbf{167}}(2), 85--146 (2003).


\othercit
\bibitem{lew04a}
 \textsc{A.~Lew},  \textsc{J.\,E. Marsden},  \textsc{M.~Ortiz},  and
  \textsc{M.~West},
An overview of variational integrators,
 in: Finite Element Methods: 1970's and Beyond, edited by L.\,P. Franca, T.\,E.
  Tezduyar,  and A.~Masud (CIMNE, Barcelona, 2004).


\bibitem{lew04b}
 \textsc{A.~Lew},  \textsc{J.\,E. Marsden},  \textsc{M.~Ortiz},  and
  \textsc{M.~West},
Variational time integrators,
 \jr{Int. J. Numer. Methods Eng.} \textbf{\textbf{60}}, 153--212 (2004).


\bibitem{leyendecker06}
 \textsc{S.~Leyendecker},  \textsc{P.~Betsch},  and  \textsc{P.~Steinmann},
Objective energy-momentum conserving integration for the constrained dynamics
  of geometrically exact beams,
 \jr{Comput. Methods Appl. Mech. Eng.} \textbf{\textbf{195}}, 2313--2333
  (2006).


\bibitem{leyendecker08}
 \textsc{S.~Leyendecker},  \textsc{J.\,E. Marsden},  and  \textsc{M.~Ortiz},
Variational integrators for constrained dynamical systems,
 \jr{Z. angew. Math. Mech. (ZAMM)} \textbf{\textbf{88}}(9), 677--708 (2008).


\othercit
\bibitem{leyendecker13}
 \textsc{S.~Leyendecker} and  \textsc{S.~Ober-Bl{\"o}baum},
A variational approach to multirate integration for constrained systems,
 in: Multibody Dynamics, edited by J.\,C. Samin and P.~Fisette, Computational
  Methods in Applied Sciences Vol.\,28 (Springer Netherlands, 2013),
  pp.\,97--121.


\bibitem{leyendecker10}
 \textsc{S.~Leyendecker},  \textsc{S.~Ober-Bl{\"o}baum},  \textsc{J.\,E.
  Marsden},  and  \textsc{M.~Ortiz},
Discrete mechanics and optimal control for constrained systems,
 \jr{Optim. Control Appl. Methods} \textbf{\textbf{31}}(6), 505--528 (2010).


\bibitem{marsden98}
 \textsc{J.\,E. Marsden},  \textsc{G.\,W. Patrick},  and  \textsc{S.~Shkoller},
Multisymplectic geometry, variational integrators, and nonlinear {PDE}s,
 \jr{Commun. Math. Phys.} \textbf{\textbf{199}}(2), 351--395 (1998).


\bibitem{marsden01}
 \textsc{J.\,E. Marsden} and  \textsc{M.~West},
Discrete mechanics and variational integrators,
 \jr{Acta Numer.} \textbf{\textbf{10}}, 357--514 (2001).


\bibitem{modak02}
 \textsc{S.~Modak} and  \textsc{E.\,D. Sotelino},
The generalized method for structural dynamics applications,
 \jr{Adv. Eng. Softw.} \textbf{\textbf{33}}, 565--575 (2002).


\bibitem{newmark59}
 \textsc{N.\,M. Newmark},
A method of computation for structural dynamics,
 \jr{ASCE J. Eng. Mech. Div.} \textbf{\textbf{85}}(EM3), 67--94 (1959).


\bibitem{oberbloebaum11}
 \textsc{S.~Ober-Bl{\"o}baum},  \textsc{O.~Junge},  and  \textsc{J.\,E.
  Marsden},
Discrete mechanics and optimal control: {A}n analysis,
 \jr{ESAIM Control Optim. Calc. Var.} \textbf{\textbf{17}}(2), 322--352 (2011).


\bibitem{oberbloebaum15}
 \textsc{S.~Ober-Bl{\"o}baum} and  \textsc{N.~Saake},
Construction and analysis of higher order {G}alerkin variational integrators,
 \jr{Adv. Comput. Math.} \textbf{\textbf{41}}, 955--986 (2015).


\bibitem{oberbloebaum13}
 \textsc{S.~Ober-Bl{\"o}baum},  \textsc{M.~Tao},  \textsc{M.~Cheng},
  \textsc{H.~Owhadi},  and  \textsc{J.\,E. Marsden},
Variational integrators for electric circuits,
 \jr{J. Comput. Phys.} \textbf{\textbf{242}}, 498--530 (2013).


\bibitem{oden69}
 \textsc{J.\,T. Oden},
A general theory of finite elements {II}. {A}pplications,
 \jr{Int. J. Numer. Methods Eng.} \textbf{\textbf{1}}, 247--259 (1969).


\bibitem{peters88}
 \textsc{D.\,A. Peters} and  \textsc{A.\,P. Izadpanah},
hp-version finite elements for the space-time domain,
 \jr{Comput. Mech.} \textbf{\textbf{3}}, 73--88 (1988).


\bibitem{reich94}
 \textsc{S.~Reich},
Momentum conserving symplectic integrators,
 \jr{Physica D} \textbf{\textbf{76}}, 375--383 (1994).


\bibitem{riff84a}
 \textsc{R.~Riff} and  \textsc{M.~Baruch},
Stability of time finite elements,
 \jr{AIAA J.} \textbf{\textbf{22}}(8), 1171--1173 (1984).


\bibitem{riff84b}
 \textsc{R.~Riff} and  \textsc{M.~Baruch},
Time finite element discretization of {H}amilton's law of varying action,
 \jr{AIAA J.} \textbf{\textbf{22}}(9), 1310--1318 (1984).


\bibitem{simkins78}
 \textsc{T.\,E. Simkins},
Unconstrained variational statements for initial and boundary-value problems,
 \jr{AIAA J.} \textbf{\textbf{16}}(6), 559--563 (1978).


\bibitem{simkins81}
 \textsc{T.\,E. Simkins},
Finite elements for initial value problems in dynamcs,
 \jr{AIAA J.} \textbf{\textbf{19}}(10), 1357--1362 (1981).


\bibitem{simo92a}
 \textsc{J.\,C. Simo} and  \textsc{N.~Tarnow},
The discrete energy-momentum method. {C}onserving algorithms for nonlinear
  elastodynamics,
 \jr{Z. angew. Math. Phys. (ZAMP)} \textbf{\textbf{43}}, 757--792 (1992).


\bibitem{simo92b}
 \textsc{J.\,C. Simo},  \textsc{N.~Tarnow},  and  \textsc{K.\,K. Wong},
Exact energy-momentum conserving algorithms and symplectic schemes for
  nonlinear dynamics,
 \jr{Comput. Methods Appl. Mech. Eng.} \textbf{\textbf{100}}, 63--116 (1992).


\bibitem{suris90}
 \textsc{Y.\,B. Suris},
Hamiltonian methods of {R}unge-{K}utta type and their variational
  interpretation,
 \jr{Math. Model.} \textbf{\textbf{2}}(4), 78--87 (1990).


\bibitem{tao10}
 \textsc{M.~Tao},  \textsc{H.~Owhadi},  and  \textsc{J.\,E. Marsden},
Non-intrusive and structure preserving multiscale integration of stiff {ODE}s,
  {SDE}s and {H}amiltonian systems with hidden slow dynamics via flow
  averaging,
 \jr{Multiscale Model. Simul.} \textbf{\textbf{8}}(4), 1269--1324 (2010).


\bibitem{wolff13}
 \textsc{S.~Wolff} and  \textsc{C.~Bucher},
Asynchronous variational integration using continuous assumed gradient
  elements,
 \jr{Comput. Methods Appl. Mech. Eng.} \textbf{\textbf{255}}, 158--166 (2013).


\bibitem{wood80}
 \textsc{W.\,L. Wood},  \textsc{M.~Bossak},  and  \textsc{O.\,C. Zienkiewicz},
An alpha modification of {N}ewmark's method,
 \jr{Int. J. Numer. Methods Eng.} \textbf{\textbf{15}}, 1562--1566 (1980).


\othercit
\bibitem{wriggers08}
 \textsc{P.~Wriggers},
Nonlinear Finite Element Methods (Springer, Berlin Heidelberg, 2008).


\othercit
\bibitem{zienkiewicz05}
 \textsc{O.\,C. Zienkiewicz} and  \textsc{R.\,L. Taylor},
The Finite Element Method for Solid and Structural Mechanics, $6^{\text{th}}$
  edition (Elsevier Ltd., Butterworth-Heinemann, 2005).


\bibitem{zienkiewicz77}
 \textsc{O.\,C. Zienkiewicz},
A new look at the {N}ewmark, {H}oubolt and other time stepping formulas. {A}
  weighted residual approach,
 \jr{Earthq. Eng. Struct. Dyn.} \textbf{\textbf{5}}, 413--418 (1977).


\end{thebibliography}

\end{document}